\documentclass[a4paper,12pt]{article}

\usepackage{epstopdf}
\usepackage {svg}
\usepackage{here}

\usepackage{amsmath}
\usepackage{amssymb}
\usepackage{amsthm}
\usepackage{mathrsfs}
\usepackage[all]{xy}
\usepackage{float}
\usepackage{graphicx}

\usepackage{hyperref}

\usepackage{tikz}
\usetikzlibrary{intersections,calc,arrows.meta}

\newtheorem{thm}{Theorem}[section]
\newtheorem{prop}[thm]{Proposition}
\newtheorem{lem}[thm]{Lemma}
\newtheorem{dfn}[thm]{Definition}

\newtheorem{corollary}[thm]{Corollary}
\newtheorem{remark}[thm]{Remark}

\newcommand{\C}{{\mathbb C}}

\newcommand{\PP}{{\mathcal P}}
\newcommand{\K}{{\mathcal K}}

\newcommand{\X}{{\mathcal X}}
\newcommand{\W}{{\mathcal W}}
\newcommand{\V}{{\mathcal V}}

\newcommand{\Z}{{\mathbb Z}}

\numberwithin{equation}{section}

\makeatletter % プリアンブルで定義開始

\@addtoreset{figure}{section}
\@addtoreset{table}{section}

\makeatother % プリアンブルで定義終了

\usepackage{lipsum} %連絡先

\usepackage{braket}

\newcommand\restr[2]{{% we make the whole thing an ordinary symbol
  \left.\kern-\nulldelimiterspace % automatically resize the bar with \right
  #1 % the function
  \vphantom{\big|} % pretend it's a little taller at normal size
  \right|_{#2} % this is the delimiter
  }}

\usepackage{comment}

\usepackage{bm}

\title{Projective covers of the simple modules for the triplet $W$-algebra $\mathcal{W}_{p_+,p_-}$}
\author{Hiromu Nakano}
\date{}

\begin{document}

\maketitle
\begin{abstract}
We study the structure of the abelian category of modules for the triplet $W$-algebra $\mathcal{W}_{p_+,p_-}$. Using the logarithmic deformation by Fjelstad et al., we construct logarithmic $\mathcal{W}_{p_+,p_-}$-modules that have $L_0$ nilpotent rank three or two. By using the structure of these logarithmic modules and the results on logarithmic Virasoro modules by Kyt\"{o}l\"{a} and Ridout, we compute ${\rm Ext}^1$ groups between certain indecomposable modules and simple modules. Based on these ${\rm Ext}^1$ groups we determine the structure of the projective covers of all $\mathcal{W}_{p_+,p_-}$-simple modules. 

\end{abstract}
%追加
\setcounter{tocdepth}{1}
\tableofcontents
\clearpage

%Arike FF22 quantum group
\section{Introduction}
In 2006, Feigin, Gainutdinov, Semikhatov and Tipunin introduced the triplet $W$-algebra $\W_{p_+,p_-}$ indexed by coprime integers $p_+$ and $p_-$ with $p_->p_+\geq 2$ \cite{FF2}. As in the case of the triplet $W$-algebra $\W_{p}$ (cf. \cite{AM,FF3,FF4,FHST,Ka}), which is a VOA extension of the Virasoro simple VOA $L(c_{1,p})$, $\mathcal{W}_{p_+,p_-}$ is an extension of the minimal Virasoro model
$
L(c_{p_+,p_-})
$, where $c_{k,l}$ is the Virasoro central charge defined by
\begin{align*}
c_{k,l}:=1-6\frac{(k-l)^2}{kl}.
\end{align*}
In the detailed studies \cite{AMW2p,AMW3p,TW}, it is shown that the triplet $W$-algebra $\W_{p_+,p_-}$ satisfies Zhu's $C_2$-cofinite condition, and the simple modules are classified.
Their results also show that the structure of $\W_{p_+,p_-}$-simple modules are similar to those of $\W_{p}$ (cf. \cite{AM,FF3,FF4}).
On the other hand, unlike $\W_{p}$, the triplet $W$-algebra $\W_{p_+,p_-}$ itself is a reducible module of length two with a simple submodule, usually denoted by $\mathcal{X}^+_{1,1}$.
This subtle difference is also closely related to the complexity of the structure of the module category.

Let $\mathcal{C}_{p_+,p_-}$ be the abelian category of $\W_{p_+,p_-}$-modules.
As shown in Theorem~\ref{prop:block}, $\mathcal{C}_{p_+,p_-}$ admits the block decomposition, similarly to the module category of $\mathcal{W}_p$ \cite[Theorem~4.4]{AM}. Each block is assigned to one of three groups: $\frac{(p_+-1)(p_--1)}{2}$ thick blocks, $p_++p_--2$ thin blocks and two semisimple blocks.
The thin blocks have counterparts in the module category of $\mathcal{W}_p$. As will be shown in this paper, their properties as full abelian subcategories are almost identical to those of the corresponding blocks in the module category of $\mathcal{W}_p$.
The most complex groups are the thick blocks and each thick block contains five simple modules, exactly one of which is the Virasoro minimal simple module. 
These blocks have no counterpart in the module category of $\mathcal{W}_p$, where no Virasoro simple module appears as a simple $\mathcal{W}_p$-module (cf. \cite{AM}).
Each thick block contains certain logarithmic modules whose $L_0$ nilpotent rank three and the detailed structure of these rank three logarithmic modules has been studied from logarithmic conformal field theories and vertex operator algebras.
In \cite{Ra}, Rasmussen examined the structure of indecomposable modules by using the method of solvable lattice models. In \cite{GRW0,GRW,W}, Gaberdiel, Runkel and Wood examined the structure of rank three logarithmic modules from the direction of boundary conformal field theories and showed that the tensor category on $\mathcal{C}_{p_+,p_-}$ is not rigid (see also \cite{MS1,Na3}). In \cite{AK}, Adamovi\'{c} and Milas constructed certain rank three logarithmic modules by using the method of the lattice construction.
However, unlike the case of the triplet $W$-algebra $\W_{p}$ \cite{McRae,NT}, there were few detailed studies of the structure of the abelian category $\mathcal{C}_{p_+,p_-}$, such as the properties of ${\rm Ext}^1$-groups.

The main goal of this paper is to determine the structure of logarithmic $\mathcal{W}_{p_+,p_-}$-modules in the thin and thick blocks, and show that the rank three logarithmic modules are projective covers of the corresponding top composition factors.
As stated in Theorem~\ref{logarithmic1}, any rank three logarithmic module whose top composition factor is not a minimal module can be realized by combining the logarithmic deformation method developed by Fjelstad et al.~\cite{FJ} with the construction introduced in \cite{AK}.
The main technical difficulties in proving the projectivity of these modules are the determination of the socle series structure and the computation of ${\rm Ext}^1$ groups within each thick block.
Although the logarithmic deformation provides an explicit construction of rank three logarithmic modules, its definition is rather complicated, and the structure of these logarithmic modules is not clear from the definition itself. This leads to the difficulties mentioned above.
To overcome these difficulties, in this paper we study the structure of these $\mathcal{W}_{p_+,p_-}$-modules and their ${\rm Ext}^1$ properties by viewing them as Virasoro modules, that is, by restricting the $\mathcal{W}_{p_+,p_-}$-action to the Virasoro action.
This restriction naturally leads us to consider the corresponding category of generalised Virasoro modules. As will be described in Section~\ref{KRsection}, the nontriviality of the logarithmic coupling established by Kyt\"{o}l\"{a} and Ridout \cite{KR} plays a crucial role in the analysis of the ${\rm Ext}^1$ and ${\rm Hom}$ spaces in this category.
This logarithmic coupling is an invariant determined by the isomorphism class of a rank two logarithmic Virasoro module, and it may be regarded as a quantity characterizing the quotient structure of logarithmic Virasoro modules.
As a consequence, by combining the detailed computations of the ${\rm Ext}^1$ and ${\rm Hom}$ spaces in Sections~\ref{KRsection} and~\ref{section proj}, we obtain the main result of this paper, Theorem~\ref{77h}, which states that every rank three $\mathcal{W}_{p_+,p_-}$-module whose top composition factor is not a minimal module is projective.
In the case of the minimal modules, the structure of the corresponding projective covers can be determined by using the description of the center of the Zhu-algebra $A(\W_{p_+,p_-})$ obtained in \cite{AMW2p,AMW3p,TW}, as will be explained in Subsection~\ref{Projminimal}.

This paper is organized as follows:
In Section \ref{Basic}, we review the structure of Fock modules and the Felder complex in accordance with \cite{TW}. The basic facts in this section are frequently used in later sections.

In Section \ref{tripletW}, we introduce the vertex operator algebra $\W_{p_+,p_-}$ and review some results in \cite{AMW2p,AMW3p,TW} briefly. In Subsection \ref{W-mod}, we introduce the block decomposition of $\mathcal{C}_{p_+,p_-}$. Each block of $\mathcal{C}_{p_+,p_-}$-mod is assigned to one of three groups: $\frac{(p_+-1)(p_--1)}{2}$ thick blocks $C^{thick}_{r,s}$, $p_++p_--2$ thin blocks $C^{thin}_{r,p_-},C^{thin}_{p_+,s}$ and two semisimple blocks $C^\pm_{p_+,p_-}$. The most complex groups are the thick blocks and each thick block $C^{thick}_{r,s}$ contains five simple modules $\X^+_{r,s},\X^+_{r^\vee,s^\vee},\X^-_{r^\vee,s},\X^-_{r,s^\vee}$ and $L(h_{r,s})$, where $L(h_{r,s})$ is the minimal simple module module of the Virasoro algebra. The thick blocks and the thin blocks contain logarithmic $\W_{p_+,p_-}$-modules on which the Virasoro zero-mode $L_0$ acts non-semisimply.

In Section \ref{logarithmic section}, by gluing lattice simple modules $\V^\pm_{r,s}$ using the logarithmic deformation by \cite{FJ}, we construct logarithmic $\W_{p_+,p_-}$-modules $\mathcal{P}^\pm_{r,s}$ and $\mathcal{Q}(\X^\pm_{r,s})_{\bullet,\bullet}$ whose $L_0$ nilpotent rank three and two, respectively.

In Section \ref{KRsection}, we investigate the ${\rm Ext}^1$-groups between certain Virasoro modules of finite length and the quotient structures of certain logarithmic Virasoro modules in preparation for Section~\ref{section proj}.
As described above, the nontriviality of the logarithmic coupling for rank-two logarithmic modules, established in \cite{KR}, plays an important role in this section.

In Section \ref{section proj}, we determine the structure of the projective covers of all $\W_{p_+,p_-}$-simple modules. 
In Subsections \ref{structureQ} and \ref{W-mod}, we determine the structre of the socle series of indecomposable modules $\mathcal{Q}(\X^\pm_{\bullet,\bullet})_{\bullet,\bullet}$ and the ${\rm Ext}^1$-groups between all simple modules.
In Subsections \ref{54} and \ref{53}, we study the structure of ${\rm Ext}^1$-groups in the thin blocks $C^{thin}_{r,p_-}$, $C^{thin}_{p_+,s}$ and the thick blocks $C^{thick}_{r,s}$. 
Based on the structure of the logarithmic Virasoro modules determined in Section \ref{KRsection}, we compute ${\rm Ext}^1$ groups between certain indecomposable $\W_{p_+,p_-}$-modules and the simple modules, and show that the logarithmic modules $\mathcal{Q}(\X^\pm_{\bullet,p_-})_{\bullet,p_-}$, $\mathcal{Q}(\X^\pm_{p_+,\bullet})_{p_+,\bullet}$ and $\mathcal{P}^\pm_{\bullet,\bullet}$ are the projective covers of the top composition factors.
In Subsection \ref{Projminimal}, we determine the structure of the projective covers of the minimal simple modules $L(h_{r,s})$ by using the structure of the center of the Zhu-algebra $A(\W_{p_+,p_-})$ determined in \cite{AMW2p,AMW3p,TW}.

\section{Bosonic Fock modules}
\label{Basic}
Recall that the Virasoro algebra $\mathcal{L}$ is the Lie algebra over $\mathbb{C}$ generated by $\mathfrak{L}_n(n\in \mathbb{Z})$ and $C$ (the central charge) with the relations
\begin{align*}
&[\mathfrak{L}_m,\mathfrak{L}_n]=(m-n)\mathfrak{L}_{m+n}+\frac{m^3-m}{12}C\delta_{m+n,0},
&[\mathfrak{L}_n,C]=0.
\end{align*}
Fix two coprime integers $p_+,p_-$ such that $p_->p_+\geq 2$, and let
\begin{align*}
c_{p_+,p_-}:=1-6\frac{(p_+-p_-)^2}{p_+p_-}
\end{align*}
be the central charge of the minimal model ${\rm M}(p_+,p_-)$.
%Fix two coprime integers $p_+,p_-$ such that $p_->p_+\geq 2$. 
In this section, we briefly review theories of Fock modules whose central charges are $C=c_{p_+,p_-}\cdot{\rm id}$
in accordance with \cite{TW}. As for the representation theory of the Virasoro algebra, see \cite{FF} and \cite{IK}. 
\subsection{Free field theory}
The Heisenberg Lie algebra 
\begin{align*}
\mathcal{H}=\bigoplus_{n\in\mathbb{Z}}\mathbb{C} a_{n}\oplus \mathbb{C} K_{\mathcal{H}}
\end{align*}
is the Lie algebra whose commutation is given by
\begin{align*} 
&[a_m,a_n]=m\delta_{m+n,0}K_{\mathcal{H}},
&[K_{\mathcal{H}},\mathcal{H}]=0.
\end{align*}
Let
\begin{align*}
\mathcal{H}^\pm&=\bigoplus_{n>0}\mathbb{C} a_{\pm n},&\mathcal{H}^0&=\mathbb{C} a_0\oplus \mathbb{C} K_{\mathcal{H}},\\
\mathcal{H}^{\geq}&=\mathcal{H}^+\oplus \mathcal{H}^0,&\mathcal{H}^{\leq }&=\mathcal{H}^-\oplus \mathcal{H}^0.
\end{align*}
For any $\alpha\in\mathbb{C}$, let $\mathbb{C}{\mid}\alpha\rangle$ be the one dimensional $\mathcal{H}^\geq$-module defined by
\begin{align*}
&a_n{\mid}\alpha\rangle=\delta_{n,0}\alpha{\mid}\alpha\rangle\ (n\geq 0),
&K_{\mathcal{H}}{\mid}\alpha\rangle={\mid}\alpha\rangle. 
\end{align*}
For any $\alpha\in \mathbb{C}$, the bosonic Fock module is defined by the induced module
\begin{align*}
F_{\alpha}={\rm Ind}_{\mathcal{H}^\geq}^{\mathcal{H}}\mathbb{C}{\mid}\alpha\rangle.
\end{align*} 
Let 
\begin{equation*}
a(z)=\sum_{n\in\mathbb{Z}}a_nz^{-n-1}
\end{equation*}
be the bosonic current. Then we have the following operator expansion
\begin{equation*}
 a(z)  a(w)=\frac{1}{(z-w)^2}+\cdots,
\end{equation*}
where $\cdots$ denotes the regular part in $z=w$.
We define the energy-momentum tensor
\begin{equation*}
T(z):=\frac{1}{2}:a(z)a(z):+\frac{\alpha_0}{2}\partial a(z),\ \ \ \ \ \alpha_0:=\sqrt{\frac{2p_-}{p_+}}-\sqrt{\frac{2p_+}{p_-}}.
\end{equation*}
where $:\ :$ is the normal ordered product. The energy-momentum tensor satisfies the following operator expansion
\begin{align*}
T(z)T(w)= \frac{c_{p_+,p_-}}{2(z-w)^4}+\frac{2T(w)}{(z-w)^2}+\frac{\partial T(w)}{z-w}+\cdots.
\end{align*}
We denote by $\sum_{n\in\mathbb{Z}}L_{n}z^{-n-2}$ the Fourier mode expansion of $T(z)$.
Then the modes 
$\{
L_{n}
\}$
generate the Virasoro algebra at central charge $c_{p_+,p_-}$.
We set
\begin{align}
\label{h_alpha}
h_{\alpha}:=\frac{1}{2}\alpha(\alpha-\alpha_0).
\end{align}
From the universality of the Virasoro Verma module $M(h_{\alpha},c_{p_+,p_-})$ with lowest weight $h_{\alpha}$ and central charge $c_{p_+,p_-}$,
there exists a homomorphism
$
\pi : M(h_{\alpha},c_{p_+,p_-})\to F_{\alpha}
$ 
such that 
\begin{align*}
&\pi(\ket{h_{\alpha}})=\ket{\alpha},
&\pi(\mathfrak{L}_n)=L_n,
\end{align*}
where $\ket{h_\alpha}$ is the lowest weight vector of $M(h_{\alpha},c_{p_+,p_-})$. 
Then, by the field $T(z)$, each Fock module $F_{\alpha}$ admits the structure of a Virasoro module with the central charge $c_{p_+,p_-}$, and has the $L_0$ weight decomposition
\begin{align*}
&F_{\alpha}=\bigoplus_{n\in \mathbb{Z}_{\geq 0}}F_{\alpha}[n],
&F_{\alpha}[n]:=\{v\in F_{\alpha}\setminus\{0\}\mid L_0v=(h_{\alpha}+n)v\}.
\end{align*}
Each weight space $F_{\alpha}[n]$ has a basis
\begin{align*}
\{a_{-\lambda}{\mid}\alpha\rangle\mid \lambda\vdash n\}
\end{align*}
with $a_{-\lambda}=a_{-\lambda_k}\cdots a_{-\lambda_1}$ for a partition $\lambda=(\lambda_1,\dots,\lambda_k)$.

The dual Fock module $F^\vee_{\alpha}$ is defined by
\begin{align*}
F^\vee_{\alpha}={\rm Ind}_{\mathcal{H}^\leq}^{\mathcal{H}}\mathbb{C}\langle\alpha{\mid},
\end{align*} 
where $\mathbb{C}\langle\alpha{\mid}$ is the one dimensional $\mathcal{H}^\leq$-module defined by
\begin{align*}
&\langle\alpha{\mid}a_n=\delta_{n,0}\alpha\langle\alpha{\mid},\qquad n\leq 0,
&K_{\mathcal{H}}\langle\alpha{\mid}=\langle\alpha{\mid}. 
\end{align*}
We see that the two Fock modules $F_{\alpha},F^\vee_{\alpha}$ are equipped with an inner product
\begin{align}
\label{inn-pro}
(\ \cdot\ ,\ \cdot\ ):\ F^\vee_{\alpha}\times F_{\alpha}\rightarrow \mathbb{C}
\end{align}
defined by 
\begin{align*}
(\bra{\alpha},\ket{\alpha})=1,\qquad (\bra{\alpha}u_1,u_2\ket{\alpha})=(\bra{\alpha}u_1u_2,\ket{\alpha})=(\bra{\alpha},u_1u_2\ket{\alpha})
\end{align*}
for $u_1,u_2\in U(\mathcal{H})$, where $U(\mathcal{H})$ is the universal enveloping algebra of $\mathcal{H}$.
In this paper, we use the notation $\bra{\alpha}u_1u_2\ket{\alpha}=(\bra{\alpha},u_1u_2\ket{\alpha})$.

We define the following conformal vector in $F_0$
\begin{align*}
T=\frac{1}{2}(a^2_{-1}+\alpha_0 a_{-2})\ket{0}\in F_0[2].
\end{align*}
\begin{dfn}
The Fock module ${F}_0$ carries the structure of a $\Z_{\geq 0}$-graded vertex operator algebra, with
\begin{align*}
&Y(\ket{0};z)={\rm id},\ \ \ \ \ \ Y(a_{-1}\ket{0};z)=a(z),\ \ \ \ \ \ Y(T;z)=T(z).
\end{align*}
We denote this vertex operator algebra by $\mathcal{F}^{}_{\alpha_0}$.
\end{dfn}

\subsection{Structure of Fock modules}
We set
\begin{align*}
&\alpha_+=\sqrt{\frac{2p_-}{p_+}},
&\alpha_-=-\sqrt{\frac{2p_+}{p_-}}.
\end{align*}
For $r,s,n\in\Z$ we introduce the symbols
\begin{align}
&\alpha_{r,s;n}=\frac{1-r}{2}\alpha_++\frac{1-s}{2}\alpha_-+\frac{\sqrt{2p_+p_-}}{2}n,
&\alpha_{r,s}=\alpha_{r,s;0}
\label{alpha_{r,s;n}}
\end{align}
and use the simple notation
\begin{align}
&F_{r,s;n}=F_{\alpha_{r,s;n}},
&F_{r,s}&=F_{\alpha_{r,s}},\nonumber\\
&h_{r,s;n}=h_{\alpha_{r,s;n}},
&h_{r,s}&=h_{\alpha_{r,s}}.\label{eq:not-hrs}
\end{align}
Note that
\begin{align*}
&h_{r,s;n}=h_{r-np_+,s}=h_{r,s+np_-},
&h_{r,s;n}=h_{-r,-s;-n}
\end{align*}
for $r,s,n\in\Z$. For each $r,s,n\in\Z$, let $L(h_{r,s;n})$ be the simple Virasoro module whose lowest weight is $h_{r,s;n}$ and the central charge $C=c_{p_+,p_-}\cdot{\rm id}$. 

Before describing the structure of Fock modules, let us introduce the notion of socle series. 
\begin{dfn}
Let $V$ be a vertex operator algebra or the Virasoro algebra, and let $M$ be a finite length $V$-module.
Let ${\rm S}(M)$ denote the maximal semisimple submodule of $M$, called the {\rm socle}. 
Furthermore, we call the following sequence the {\rm socle} {\rm series} of $M$:
\begin{align*}
{\rm S}_1(M)\subsetneq {\rm S}_2(M)\subsetneq \cdots \subsetneq {\rm S}_n(M)=M
\end{align*}
with ${\rm S}_1(M)={\rm S}(M)$ and ${\rm S}_{i+1}(M)/{\rm S}_{i}(M)={\rm S}(M/{\rm S}_{i}(M))$.
\end{dfn}
The following proposition is due to \cite{FF}.
\begin{prop}
\label{FockSocle}
Let $r$ and $s$ be integers satisfying $1\leq r\leq  p_+$, $1\leq s\leq p_-$, and let $n\in \mathbb{Z}$.
As the Virasoro modules, there are four cases of socle series for the Fock modules $F_{r,s;n}\in \mathcal{F}_{\alpha_0}\mathchar`-{\rm mod}$:
\begin{enumerate}
\item Let $1\leq r< p_+$, $1\leq s<p_-$. Then we have
\begin{align*}
{\rm S}_1(F_{r,s;n})\subsetneq {\rm S}_2(F_{r,s;n})\subsetneq {\rm S}_3(F_{r,s;n})=F_{r,s;n}
\end{align*} 
such that
\begin{align*}
&{\rm S}_1(F_{r,s;n})={\rm S}(F_{r,s;n})=\bigoplus_{k\geq 0}L(h_{r,p_--s;|n|+2k+1}),\\
&{\rm S}_2(F_{r,s;n})/{\rm S}_1(F_{r,s;n})=\bigoplus_{k\geq a}L(h_{r,s;|n|+2k})\oplus \bigoplus_{k\geq 1-a}L(h_{p_+-r,p_--s;|n|+2k}),\\
&{\rm S}_3(F_{r,s;n})/{\rm S}_2(F_{r,s;n})={\rm S}(F_{r,s;n}/{\rm S}_2(F_{r,s;n}))=\bigoplus_{k\geq 0}L(h_{p_+-r,s;|n|+2k+1}),
\end{align*}
where $a=0$ if $n\geq 0$ and $a=1$ if $n<0$.

\item Let $1\leq s<p_-$. Then we have
\begin{align*}
{\rm S}_1(F_{p_+,s;n})\subsetneq {\rm S}_2(F_{p_+,s;n})=F_{p_+,s;n}
\end{align*}
such that
\begin{align*}
&{\rm S}_1(F_{p_+,s;n})={\rm S}(F_{p_+,s;n})=\bigoplus_{k\geq 0}L(h_{p_+,p_--s;|n|+2k+1}),\\
&{\rm S}_2(F_{p_+,s;n})/{\rm S}_1(F_{p_+,s;n})=\bigoplus_{k\geq a}L(h_{p_+,s;|n|+2k})
\end{align*}
where $a=0$ if $n\geq 1$ and $a=1$ if $n<1$.

\item Let $1\leq r < p_+$. Then we have
\begin{align*}
{\rm S}_1(F_{r,p_-;n})\subsetneq {\rm S}_2(F_{r,p_-;n})=F_{r,p_-;n}
\end{align*}
such that
\begin{align*}
&{\rm S}_1(F_{r,p_-;n})={\rm S}(F_{r,p_-;n})=\bigoplus_{k\geq 0}L(h_{r,p_-;|n|+2k}),\\
&{\rm S}_2(F_{r,p_-;n})/{\rm S}_1(F_{r,p_-;n})=\bigoplus_{k\geq a}L(h_{p_+-r,p_-;|n|+2k-1})
\end{align*}
where $a=1$ if $n\geq 0$ and $a=0$ if $n<0$.

\item The Fock module $F_{p_+,p_-;n}$  is semi-simple as a Virasoro module:
\begin{align*}
F_{p_+,p_-;n}={\rm S}(F_{p_+,p_-;n})=\bigoplus_{k\geq 0}L(h_{p_+,p_-;|n|+2k}).
\end{align*}
\end{enumerate}
\end{prop}
In the above proposition, the first family of the Fock modules is referred to as braided type, the second and third as chain type, and the fourth as semisimple type.

\subsection{Screening currents and Felder complex}
We introduce a free scalar field $\phi(z)$, which is a formal primitive of $a(z)$
\begin{equation*}
\phi(z)=\hat{a}+a_0\log z-\sum_{n\neq 0}\frac{a_n}{n}z^{-n}
\end{equation*}
where $\hat{a}$ is defined by
\begin{equation}
\label{conjugate}
[a_m,\hat{a}]=\delta_{m,0}{\rm id}.
\end{equation}
\begin{remark}
For each $\gamma \in \mathbb{C}$, $e^{\gamma\hat{a}}$ defines a Heisenberg weight shifting map $e^{\gamma\hat{a}}: F_{\beta}\rightarrow F_{\beta+\gamma}$. From now on, we identify $e^{\gamma\hat{a}}\cdot F_{\beta}$ with $F_{\beta+\gamma}$.
\label{rem:beta-gamma}
\end{remark}
The scalar field $\phi(z)$ satisfies the operator product expansion
\begin{equation*}
\phi(z)\phi(w)=\log (z-w)+\cdots.
\end{equation*} 
For any $\alpha\in\C$ we introduce the field $V_{\alpha}(z)$
\begin{align*}
&V_{\alpha}(z)=:e^{\alpha\phi(z)}:=e^{\alpha\hat{a}}z^{\alpha a_{0}}\overline{V}_{\alpha}(z),\ z^{\alpha a_{0}}=e^{\alpha a_{0}\log z}\ ,\\
&\overline{V}_{\alpha}(z)=e^{\alpha\sum_{n\geq 1}\frac{a_{-n}}{n}z^{n}}e^{-\alpha\sum_{n \geq 1}\frac{a_n}{n}z^{-n}}.
\end{align*}
In the context of vertex operator algebras, this field is represented by $Y(\ket{\alpha};z)$. 
By definition, the following operator product expansion holds
\begin{align}
\label{eq:vavb}
V_{\alpha}(z)V_{\beta}(w)= (z-w)^{\alpha\beta}:{V}_{\alpha}(z){V}_{\beta}(w):.
\end{align}
We introduce the following two screening currents ${Q}_{+}(z),{Q}_-(z)$
\begin{align*}
{Q}_{\pm}(z)=V_{\alpha_\pm}(z)
%:e^{\alpha_{\pm}\phi(z)}:
\end{align*} 
whose conformal weights are $h_{\alpha_{\pm}}=1$. In fact,
\begin{align*}
T(z){Q}_{\pm}(w)&=\frac{{Q}_{\pm}(w)}{(z-w)^2}+\frac{\partial_{w}{Q}_{\pm}(w)}{z-w}+\cdots\\
&=\partial_{w}\bigl(\frac{{Q}_{\pm}(w)}{z-w}\bigr)+\cdots.
\end{align*}
Therefore the zero modes of the fields ${Q}_{\pm}(z)$
\begin{align*}
&{\rm Res}_{z=0}{Q}_{+}(z){\rm d}z={Q}_{+}\ :\ F_{1,k}\rightarrow F_{-1,k},\ \ k\in\Z\\
&{\rm Res}_{z=0}{Q}_{-}(z){\rm d}z={Q}_{-}\ :\ F_{k,1}\rightarrow F_{k,-1},\ \ k\in\Z
\end{align*}
commute with every Virasoro mode.

For $n\geq 2$, we introduce more complicated screening currents 
\begin{align*}
&{Q}^{[n]}_{+}(z)\in {\rm Hom}_{\C}(F_{n,k},F_{-n,k})[[z,z^{-1}]],\ \ k\in\Z,\\
&{Q}^{[n]}_{-}(z)\in {\rm Hom}_{\C}(F_{k,n},F_{k,-n})[[z,z^{-1}]],\ \ k\in\Z,
\end{align*}
as follows
\begin{equation}
\label{Tsuchiya-Kanie}
\begin{split}
&{Q}^{[n]}_{\pm}(z)=\int_{\Delta^{\pm}_{n-1}}Q_\pm(z)Q_\pm(zx_{1})Q_\pm(zx_2)\cdots Q_\pm(zx_{n-1})z^{n-1}{\rm d}\bm{x},
\end{split}
\end{equation}
where 
we have used the notation
\begin{align*}
{\rm d}\bm{x}={\rm d}x_1\wedge \cdots\wedge {\rm d}x_{n-1},
\end{align*}
and 
$\Delta^{\pm}_{n-1}$ denotes the twisted cycle constructed from the simplex 
$
\{\bm{x}\in \mathbb{R}^{n-1}|1>x_1>\cdots >x_{n-1}>0\}
$
according to the multivaluedness of the integrand by applying the construction of \cite{TK}.
These fields satisfy the following operator product expansion (cf. \cite{TW})
\begin{align*}
&T(z){Q}^{[n]}_{\pm}(w)=\frac{{Q}^{[n]}_{\pm}(w)}{(z-w)^2}+\frac{\partial_{w}{Q}^{[n]}_{\pm}(w)}{z-w}+\cdots.
\end{align*}
In particular the following proposition holds.
\begin{prop}
For $n\geq 2$, define
\begin{align*}
&Q^{[n]}_+={\rm Res}_{z=0}{Q}^{[n]}_{+}(z){\rm d}z,
&Q^{[n]}_-={\rm Res}_{z=0}{Q}^{[n]}_{-}(z){\rm d}z.
\end{align*}
Then these zero modes commute with every Virasoro mode of $\mathcal{F}_{\alpha_0}\mathchar`-{\rm Mod}$. These zero modes are called {\rm screening operators}.
\end{prop}
In what follows, we use the notation $Q^{[1]}_\pm:=Q_\pm$.

For integers $r$ and $s$ satisfying $1\leq r< p_+$, $1\leq s< p_-$, we set
\begin{align*}
&r^{\vee}(p_+):=p_+-r,
&s^{\vee}(p_-):=p_--s
\end{align*}
and use the abbreviations 
\begin{align}
\label{notrs}
&r^\vee=r^{\vee}(p_+),
&s^\vee=s^{\vee}(p_-).
\end{align}

We define the following Virasoro modules :
\begin{enumerate}
\item For $1\leq r<p_+,\ 1\leq s\leq p_-,\ n\in\Z$ 
\begin{align*}
K_{r,s;n;+}&={\rm ker}Q^{[r]}_+:F_{r,s;n}\rightarrow F_{r^\vee,s;n+1},\\
X_{r^\vee,s;n+1;+}&={\rm im}Q^{[r]}_+:F_{r,s;n}\rightarrow F_{r^\vee,s;n+1}.
\end{align*}
\item For $1\leq r\leq p_+,\ 1\leq s<p_-,\ n\in\Z$ 
\begin{align*}
K_{r,s;n;-}&={\rm ker}Q^{[s]}_-:F_{r,s;n}\rightarrow F_{r,s^\vee;n-1},\\
X_{r,s^\vee;n-1;-}&={\rm im}Q^{[s]}_-:F_{r,s;n}\rightarrow F_{r,s^\vee;n-1}.
\end{align*}
\end{enumerate}

The following propositions are due to \cite{Felder}.
\begin{prop}
\label{Felder complex}
Let $r$ and $s$ be integers satisfying $1\leq r\leq  p_+$, $1\leq s\leq p_-$, and let $n\in \mathbb{Z}$.
The socle series of $K_{r,s;n;\pm}$ and $X_{r,s;n;\pm}$ are given by :
\begin{enumerate}
\item Let $1\leq r< p_+$ and $1\leq s< p_-$. Then we have
\begin{align*}
&{\rm S}_1(K_{r,s;n;\pm})=\subsetneq K_{r,s;n;\pm},
&{\rm S}_1(X_{r,s;n;\pm})=\subsetneq X_{r,s;n;\pm}
\end{align*} 
such that
\begin{align*}
n&\geq 0&n&\leq -1\\
\hspace{-6mm}{\rm S}_1(K_{r,s;n;+})&=\bigoplus_{k\geq 1}L(h_{r,s^\vee;n+2k-1}),&{\rm S}_1(K_{r,s;n;+})&=\bigoplus_{k\geq 1}L(h_{r,s^\vee;-n+2k-1}),\\
\hspace{-6mm}K_{r,s;n;+}/{\rm S}_1&=\bigoplus_{k\geq 1}L(h_{r,s;n+2(k-1)}),&K_{r,s;n;+}/{\rm S}_1&=\bigoplus_{k\geq 1}L(h_{r,s;-n+2k}),\\
\hspace{-6mm}{\rm S}_1(X_{r,s;n+1;+})&=\bigoplus_{k\geq 1}L(h_{r,s^\vee;n+2k}),&{\rm S}_1(X_{r,s;n+1;+})&=\bigoplus_{k\geq 1}L(h_{r,s^\vee;-n+2(k-1)}),\\
\hspace{-6mm}X_{r,s;n+1;+}/{\rm S}_1&=\bigoplus_{k\geq 1}L(h_{r,s;n+2k-1}),&X_{r,s;n+1;+}/{\rm S}_1&=\bigoplus_{k\geq 1}L(h_{r,s;-n+2k-1}).
\end{align*}
\begin{align*}
n&\geq 1&n&\leq 0\\
\hspace{-6mm}{\rm S}_1(K_{r,s;n;-})&=\bigoplus_{k\geq 1}L(h_{r,s^\vee;n+2k-1}),&{\rm S}_1(K_{r,s;n;-})&=\bigoplus_{k\geq 1}L(h_{r,s^\vee;-n+2k-1}),\\
\hspace{-6mm}K_{r,s;n;-}/{\rm S}_1&=\bigoplus_{k\geq 1}L(h_{r,s;n+2(k-1)}),&K_{r,s;n;-}/{\rm S}_1&=\bigoplus_{k\geq 1}L(h_{r,s;-n+2k}),\\
\hspace{-6mm}{\rm S}_1(X_{r,s;n+1;-})&=\bigoplus_{k\geq 1}L(h_{r,s^\vee;n+2(k-1)}),&{\rm S}_1(X_{r,s;n+1;-})&=\bigoplus_{k\geq 1}L(h_{r,s^\vee;-n+2k}),\\
\hspace{-6mm}X_{r,s;n+1;-}/{\rm S}_1&=\bigoplus_{k\geq 1}L(h_{r^\vee,s^\vee;n+2k-1}),&X_{r,s;n+1;-}/{\rm S}_1&=\bigoplus_{k\geq 1}L(h_{r^\vee,s^\vee;-n+2k-1}).
\end{align*}
\item Let $1\leq r< p_+$ and $s=p_-$. Then we have
$
X_{r,p_-;n}={\rm S}(F_{r,p_-;n}).
$
\item Let $r=p_+$ and $1\leq s< p_-$. Then we have
$
X_{p_+,s;n}={\rm S}(F_{p_+,s;n}).
$
\end{enumerate}

\end{prop}

\begin{prop}[\cite{Felder}]
\label{Felder complex2}
Let $r$ and $s$ be integers satisfying $1\leq r\leq  p_+$, $1\leq s\leq p_-$.
\begin{enumerate}
\item 
The screening operators $Q^{[r]}_+$ and $Q^{[r^{\vee}]}_+$ define the {complex}
\begin{align*}
\cdots\xrightarrow{Q^{[r]}_+}F_{r^{\vee},s;n-1}\xrightarrow{Q^{[r^{\vee}]}_+}F_{r,s;n}\xrightarrow{Q^{[r]}_+}F_{r^{\vee},s;n+1}\xrightarrow{Q^{[r^{\vee}]}_+}\cdots.
\end{align*}
This complex is exact everywhere except in $F_{r,s}=F_{r,s;0}$ where the cohomology is given by
\begin{align*}
{\rm ker}Q^{[r]}_+/{\rm im}Q^{[r^\vee]}_+\simeq L(h_{r,s;0}).
\end{align*}
\item The screening operators $Q^{[s]}_-$ and $Q^{[s^{\vee}]}_-$ define the {complex}
\begin{align*}
\cdots\xrightarrow{Q^{[s]}_-}F_{r,s^{\vee};n+1}\xrightarrow{Q^{[s^{\vee}]}_-}F_{r,s;n}\xrightarrow{Q^{[s]}_-}F_{r,s^{\vee};n-1}\xrightarrow{Q^{[s^{\vee}]}_-}\cdots.
\end{align*}
This complex is exact everywhere except in $F_{r,s}=F_{r,s;0}$ where the cohomology is given by
\begin{align*}
{\rm ker}Q^{[s]}_-/{\rm im}Q^{[s^\vee]}_-\simeq L(h_{r,s;0}).
\end{align*}
\item The screening operators $Q^{[r]}_+$ and $Q^{[r^{\vee}]}_+$ define the {complex}
\begin{align*}
\cdots\xrightarrow{Q^{[r]}_+}F_{r^{\vee},p_-;n-1}\xrightarrow{Q^{[r^{\vee}]}_+}F_{r,p_-;n}\xrightarrow{Q^{[r]}_+}F_{r^{\vee},p_-;n+1}\xrightarrow{Q^{[r^{\vee}]}_+}\cdots
\end{align*}
and this complex is exact.
\item The screening operators $Q^{[s]}_-$ and $Q^{[s^{\vee}]}_-$ define the {complex}
\begin{align*}
\cdots\xrightarrow{Q^{[s]}_-}F_{p_+,s^{\vee};n+1}\xrightarrow{Q^{[s^{\vee}]}_-}F_{p_+,s;n}\xrightarrow{Q^{[s]}_-}F_{p_+,s^{\vee};n-1}\xrightarrow{Q^{[s^{\vee}]}_-}\cdots.
\end{align*}
and this complex is exact.

\end{enumerate}
\end{prop}
We call the complexes introduced above $Felder$ $complexes$.

\section{The triplet $W$-algebra $\W_{p_+,p_-}$}
\label{tripletW}
In this section, we introduce a vertex operator algebra $\W_{p_+,p_-}$ which is called the triplet $W$-algebra of type $(p_+,p_-)$ and review some results in \cite{AMW2p,AMW3p,TW} briefly. In Subsection \ref{W-modori}, we introduce the abelian category of $\W_{p_+,p_-}$-modules and the block decomposition of this abelian category.
\subsection{The lattice vertex operator algebra and the vertex operator algebra $\W_{p_+,p_-}$}

%In this subsection, we introduce the triplet $W$-algebra $\W_{p_+,p_-}$ given by \cite{FF2}. This vertex operator algebra is a generalization of the triplet $W$-algebra $\W_p$ \cite{AM,FF3,FF4,Ka}.

We define the following lattice vertex operator algebra associated to the integer lattice $\mathbb{Z}\sqrt{2p_+p_-}$. 
%whose $\langle,\langle$
\begin{dfn}
%\mbox{}\\
The {\rm lattice vertex operator algebra} $\mathcal{V}_{[p_+,p_-]}$ is the tuple 
\begin{align*}
(\V^{+}_{1,1}, \ket{0}, \frac{1}{2}(a^2_{-1}-\alpha_0a_{-2})\ket{0}, Y),
\end{align*}
where the underlying vector space of $\mathcal{V}_{[p_+,p_-]}$ is given by
\begin{align*}
\V^{+}_{1,1}=\bigoplus_{n\in\Z}F_{1,1;2n}=\bigoplus_{n\in\Z}F_{n\sqrt{2p_+p_-}},
\end{align*}
and the vertex operator $Y$ satisfies $Y(\ket{\alpha_{1,1;2n}},z)=V_{\alpha_{1,1;2n}}(z)$ for $n\in\Z$.
\end{dfn}
It is a known fact that simple $\mathcal{V}_{[p_+,p_-]}$-modules are given by the following $2p_+p_-$ direct sum of Fock modules
\begin{align*}
&\V^{+}_{r,s}=\bigoplus_{n\in\Z}F_{r,s;2n},
&\V^{-}_{r,s}=\bigoplus_{n\in\Z}F_{r,s;2n+1},
\end{align*}
where $1\leq r\leq p_+,1\leq s\leq p_-$.

Note that the two screening operators $Q_+$ and $Q_-$ act on $\V^+_{1,1}$. We define the following vector subspace of $\V^+_{1,1}$:
\begin{align*}
\K_{1,1}={\rm ker}Q_+\cap{\rm ker}Q_-\subset \V^+_{1,1}.
\end{align*}
\begin{dfn}
The {\rm triplet $W$-algebra} is defined by the following vertex operator algebra
\begin{align*}
\W_{p_+,p_-}:=(\K_{1,1},\ket{0},T,Y),
\end{align*}
where the vacuum vector, conformal vector and vertex operator map are those of $\mathcal{V}_{[p_+,p_-]}$.
\end{dfn}
We define
\begin{align*}
&W^+:=Q^{[p_--1]}_-\ket{\alpha_{1,p_--1;3}},
&W^-&:=Q^{[p_+-1]}_+\ket{\alpha_{p_+-1,1;-3}},\\ 
&W^0:=Q^{[2p_+-1]}_+\ket{\alpha_{p_+-1,1;-3}}.
\end{align*}
The $L_0$ weights of these vectors are $h_{4p_+-1,1}$.
The following two propositions are known to characterize $\mathcal{W}_{p_+,p_-}$.
\begin{prop}[\cite{AMW2p,AMW3p,TW}]
\label{genW}
$\W_{p_+,p_-}$ is strongly  generated by the fields $T(z),Y(W^{\pm},z),Y(W^0,z)$.
\end{prop}

\begin{prop}[\cite{AMW2p,AMW3p,TW}]
$\W_{p_+,p_-}$ is $C_2$-cofinite.
\end{prop}

\subsection{Simple $\W_{p_+,p_-}$-modules}
For $1\leq r\leq p_+,\ 1\leq s\leq p_-$, let $\X^{\pm}_{r,s}$ be the following vector subspace of $\V^\pm_{r,s}$: 
\begin{enumerate}
\item For $1\leq r< p_+$, $1\leq s< p_-$, we set
\begin{align*}
&\X^+_{r,s}=Q^{[r^\vee]}_+(\V^-_{r^\vee,s})\cap Q^{[s^\vee]}_-(\V^-_{r,s^\vee}),
&\X^-_{r,s}=Q^{[r^\vee]}_+(\V^+_{r^\vee,s})\cap Q^{[s^\vee]}_-(\V^+_{r,s^\vee}),
\end{align*}
where we use the notation $r^\vee=p_+-r$, $s^\vee=p_--s$ defined by (\ref{notrs}).
\item For $1\leq r< p_+,\ s=p_-$, we set
\begin{align*}
&\X^+_{r,p_-}=Q^{[r^\vee]}_+(\V^-_{r^\vee,p_-}),
&\X^-_{r,p_-}=Q^{[r^\vee]}_+(\V^+_{r^\vee,p_-}).
\end{align*}
\item For $r=p_+,\ 1\leq s< p_-$, we set
\begin{align*}
&\X^+_{p_+,s}=Q^{[s^\vee]}_-(\V^-_{p_+,s^\vee}),
&\X^-_{p_+,s}=Q^{[s^\vee]}_-(\V^+_{p_+,s^\vee}).
\end{align*}
\item For $r=p_+,\ s=p_-$, we set
\begin{align*}
&\X^+_{p_+,p_-}=\V^+_{p_+,p_-},
&\X^-_{p_+,p_-}=\V^-_{p_+,p_-}.
\end{align*}
\end{enumerate}

\begin{dfn}
\label{symbol}
\begin{enumerate}
\item We define the interior Kac table $\mathbb{K}_{p_+,p_-}$ as the following quotient set
\begin{align*}
\mathbb{K}_{p_+,p_-}=\{(r,s)|\ 1\leq r<p_+,1\leq s<p_-\}/\sim
\end{align*}
where $(r,s)\sim (r',s')$ if and only if $r'=p_+-r, s'=p_--s$. Note that $\#\mathbb{K}_{p_+,p_-}=\frac{(p_+-1)(p_--1)}{2}$.

\item For each $1\leq r\leq p_+$, $1\leq s\leq p_-$, $n\geq 0$, we define the following symbols
\begin{equation*}
\hspace{-10mm}\Delta^+_{r,s;n}{}=
\begin{cases}
h_{r^\vee,s;-2n-1}&r\neq p_+,s\neq p_-\\
h_{p_+,s;-2n}&r=p_+,s\neq p_-\\
h_{r,p_-;2n}&r\neq p_+,s=p_-\\
h_{p_+,p_-;-2n}&r=p_+,s=p_-
\end{cases}
,\ \ 
\Delta^-_{r,s;n}{}=
\begin{cases}
h_{r^\vee,s;-2n-2}&r\neq p_+,s\neq p_-\\
h_{p_+,s;-2n-1}&r=p_+,s\neq p_-\\
h_{r,p_-;2n+1}&r\neq p_+,s=p_-\\
h_{p_+,p_-;-2n-1}&r=p_+,s=p_-
\end{cases}
.\\
\end{equation*}
\end{enumerate}
\end{dfn}

\begin{prop}[\cite{AMW2p,AMW3p,TW}]
\label{simpleclass}
Each of the $\frac{(p_+-1)(p_--1)}{2}+2p_+p_-$ vector spaces
\begin{align*}
L(h_{k,l}),\ (k,l)\in\mathbb{K}_{p_+,p_-},\ \ \ \ \X^{\pm}_{r,s},\ 1\leq r\leq p_+,\ 1\leq s\leq p_-
\end{align*} 
has the structure of a simple $\W_{p_+,p_-}$-module. Furthermore, all simple $\W_{p_+,p_-}$-modules are given by these $\frac{(p_+-1)(p_--1)}{2}+2p_+p_-$ simple modules. 
\end{prop}

\begin{prop}[\cite{AMW2p,AMW3p,TW}]
\label{socleV}
Each of the $2p_+p_-$ simple $\mathcal{V}_{[p_+,p_-]}$-module becomes a $\W_{p_+,p_-}$-module by restriction, and has the following socle series:
\begin{enumerate}
\item For each $1\leq r<p_+,\ 1\leq s<p_-$, $\V^+_{r,s}$ has the following socle series
\begin{align*}
{\rm S}_1(\V^+_{r,s})\subsetneq {\rm S}_2(\V^{+}_{r,s})\subsetneq {\rm S}_3(\V^{+}_{r,s})=\V^{+}_{r,s}
\end{align*}
such that
\begin{align*}
{\rm S}_1=\X^+_{r,s},\quad {\rm S}_2/{\rm S}_1=\X^-_{r^{\vee},s}\oplus \X^{-}_{r,s^{\vee}}\oplus L(h_{r,s}),\quad \V^+_{r,s}/{\rm S}_2=\X^+_{r^{\vee},s^{\vee}}.
\end{align*}
\item
For each $1\leq r<p_+,\ 1\leq s<p_-$, $\V^-_{r,s}$ has the socle series
\begin{align*}
{\rm S}_1(\V^-_{r,s})\subsetneq {\rm S}_2(\V^{-}_{r,s})\subsetneq {\rm S}_3(\V^{-}_{r,s})=\V^{-}_{r,s}
\end{align*}
such that
\begin{align*}
{\rm S}_1=\X^-_{r,s},\quad {\rm S}_2/{\rm S}_1=\X^+_{r^{\vee},s}\oplus \X^{+}_{r,s^{\vee}},\quad 
\V^-_{r,s}/{\rm S}_2=\X^-_{r^{\vee},s^{\vee}}.
\end{align*}
\item 
For each $1\leq r<p_+$, $\V^+_{r,p_-}$ and $\V^-_{r^\vee,p_-}$ have the following socle series
\begin{align*}
\V^+_{r,p_-}/\X^+_{r,p_-}\simeq \X^-_{r^\vee,p_-},\ \ \ \ \V^-_{r^\vee,p_-}/\X^-_{r^\vee,p_-}\simeq \X^+_{r,p_-}.
\end{align*}
\item
For each $1\leq s<p_-$, $\V^+_{p_+,s}$ and $\V^-_{p_-,s^\vee}$ have the following socle series
\begin{align*}
\V^+_{p_+,s}/\X^+_{p_+,s}\simeq \X^-_{p_+,s^\vee},\ \ \ \ \V^-_{p_+,s^\vee}/\X^-_{p_+,s^\vee}\simeq \X^+_{p_+,s}.
\end{align*}
\item For $r=p_+$, $s=p_-$, we have
\begin{align*}
\V^+_{p_+,p_-}=\X^+_{p_+,p_-},\ \ \ \ \V^-_{p_+,p_-}=\X^-_{p_+,p_-}.
\end{align*}
\end{enumerate}
\end{prop}

Define the following Virasoro singular vectors.
\begin{dfn}
\label{dfn0830}
For $1\leq r\leq p_+$, $1\leq s\leq p_-$ and $n\geq 0$, we define the following vectors $\{u^{(n)}_{i}(r,s)\}_{i=-n}^{n}$ in $\mathcal{X}^+_{r,s}$:
\begin{equation*}
u^{(n)}_{i}(r,s):=
\begin{cases}
Q^{[(n+i)p_+]}_+\ket{\alpha_{p_+,p_-;-2n}}& r=p_+,\ s=p_-\\
Q^{[(n+i)p_+]}_+\ket{\alpha_{p_+,s;-2n}} & r=p_+,\ s\neq p_-\\
Q^{[(n+i+1)p_+-r]}_+\ket{\alpha_{r^\vee,p_-;-2n-1}} & r\neq p_+,\ s=p_-\\
Q^{[(n+i+1)p_+-r]}_+\ket{\alpha_{r^\vee,s;-2n-1}} & r\neq p_+,\ s\neq p_-
\end{cases}
\end{equation*}
where $Q^{[0]}_+={\rm id}$.
For $1\leq r\leq p_+$, $1\leq s\leq p_-$ and $n\geq 0$, we define the following vectors $\{v^{(n)}_{\frac{2i-1}{2}}(r,s)\}_{i=-n}^{n+1}$ in $\mathcal{X}^-_{r,s}$:
\begin{equation*}
v^{(n)}_{\frac{2i-1}{2}}(r,s):=
\begin{cases}
Q^{[(n+i)p_+]}_+\ket{\alpha_{p_+,p_-;-2n-1}}& r=p_+,\ s=p_-\\
Q^{[(n+i)p_+]}_+\ket{\alpha_{p_+,s;-2n-1}} & r=p_+,\ s\neq p_-\\
Q^{[(n+i+1)p_+-r]}_+\ket{\alpha_{r^\vee,p_-;-2n-2}} & r\neq p_+,\ s=p_-\\
Q^{[(n+i+1)p_+-r]}_+\ket{\alpha_{r^\vee,s;-2n-2}} & r\neq p_+,\ s\neq p_-
\end{cases}
.
\end{equation*}
\end{dfn}
In the following, we use the abbreviations $u^{(n)}_{i}=u^{(n)}_{i}(r,s)$ and $v^{(n)}_{\frac{2i-1}{2}}=v^{(n)}_{\frac{2i-1}{2}}(r,s)$.
From the structure of the images and kernels of the screening operators $Q^{[\bullet]}_\pm$ \cite{FF,Felder,TW}, we see that each $\X^\pm_{r,s}$ decomposes as the direct sum of  the simple Virasoro modules
\begin{align}
&\X^+_{r,s}=\bigoplus_{n\geq 0}(2n+1)L(\Delta^+_{r,s;n}),
&\X^-_{r,s}=\bigoplus_{n\geq 0}(2n+2)L(\Delta^-_{r,s;n}),
\label{decomp}
\end{align}
and the sets of vectors $\{u^{(n)}_{i}\}_{i=-n}^{n}$ and $\{v^{(n)}_{\frac{2i-1}{2}}\}_{i=-n}^{n+1}$ give a basis of the lowest weight space of $(2n+1)L(\Delta^+_{r,s;n})\subset \X^+_{r,s}$ and $(2n+2)L(\Delta^-_{r,s;n})\subset \X^-_{r,s}$, respectively. 
For $W=W^\pm,W^0$ and $n\in \mathbb{Z}$, we define the $n$-th mode of the field $Y(W,z)$ as follows
\begin{align*}
W[n]=\oint_{z=0}Y(W,z)z^{h_{4p_+-1,1}+n-1}{\rm d}z.
\end{align*}
The transitive $\W_{p_+,p_-}$-actions on the Virasoro lowest weight spaces of (\ref{decomp}) are given by the following proposition.
\begin{prop}[\cite{AMW2p,AMW3p,TW}]
\mbox{}
\begin{enumerate}
\item The Virasoro lowest weight vectors $u^{(n)}_{i}$ $(n\geq 0,-n\leq i\leq n)$ satisfy
\begin{align*}
W^\epsilon[-h]u^{(0)}_{0}=0,\ \ \ h<\Delta^+_{r,s;1}-\Delta^+_{r,s;0},\ \epsilon=\pm,
\end{align*}
and
\begin{align*}
W^\pm[0]u^{(n)}_{i}&\in \C^\times u^{(n)}_{i\pm1}+\sum_{k=0}^{n-1}U(\mathcal{L}).u^{(k)}_{i\pm 1},\\
u^{(n+1)}_{i\pm 1}&\in \C^\times W^\pm[\Delta^+_{r,s;n}-\Delta^+_{r,s;n+1}]u^{(n)}_{i} +\sum_{k=0}^{n}U(\mathcal{L}).u^{(k)}_{i\pm 1},
\end{align*}
where $u^{(n)}_{n+1}=u^{(n)}_{-n-1}=u^{(-1)}_{i}=0$, and $U(\mathcal{L})$ is the universal enveloping algebra of the Virasoro algebra $\mathcal{L}$.
\item  The Virasoro lowest weight vectors $v^{(n)}_{\frac{2i+1}{2}}$ $(n\geq 0,-n-1\leq i\leq n)$ satisfy
\begin{equation*}
W^{\epsilon}[-h]v^{(0)}_{\frac{\pm1}{2}}
\begin{cases}
=0 &\epsilon=\pm \\
\in  U(\mathcal{L}).v^{(0)}_{\frac{\pm1}{2}} &\epsilon=0\\
\in U(\mathcal{L}).v^{(0)}_{\frac{\mp1}{2}}& \epsilon=\mp
\end{cases}
\ \ \ h<\Delta^-_{r,s;1}-\Delta^-_{r,s;0},
\end{equation*}
and
\begin{align*}
W^\pm[0]v^{(n)}_{\frac{2i+1}{2}}&\in \C^\times v^{(n)}_{\frac{2i+1\pm 2}{2}}+\sum_{k=0}^{n-1}U(\mathcal{L}).v^{(k)}_{\frac{2i+1}{2}\pm 1},\\
v^{(n+1)}_{\frac{2i+1}{2}\pm 1}&\in \C^\times W^\pm[\Delta^-_{r,s;n}-\Delta^-_{r,s;n+1}]v^{(n)}_{\frac{2i+1}{2}}+\sum_{k=0}^{n}U(\mathcal{L}).v^{(k)}_{\frac{2i+1}{2}\pm 1},
\end{align*}
where $v^{(n)}_{\frac{-2n-3}{2}}=v^{(n)}_{\frac{2n+3}{2}}=v^{(-1)}_{\frac{2i+1}{2}}=0$.
\end{enumerate}
\label{sl2action2}
\end{prop}

Let $A(\W_{p_+,p_-})$ be the Zhu-algebra of $\W_{p_+,p_-}$ (see \cite{Zhu} for the definition and properties of the Zhu-algebra). 
\begin{prop}[\cite{AMW2p,AMW3p,TW}]
For any $1\leq r\leq p_+$, $1\leq s\leq p_-$, 
\begin{align*}
f(\Delta^-_{r,s;0})\neq 0.
\end{align*}
In particular, the lowest weight space of $\X^-_{r,s}$ admits the structure of a two dimensional simple $sl_2$-module with respect to the following elements
\begin{align*}
E=\frac{1}{\sqrt{2}f(\Delta^-_{r,s;0})}[W^+],\ \ \ \ F=-\frac{1}{\sqrt{2}f(\Delta^-_{r,s;0})}[W^-],\ \ \ \ H=\frac{1}{f(\Delta^-_{r,s;0})}[W^0].
\end{align*}
\label{sl2}
\end{prop}

\begin{thm}[\cite{AMW2p,AMW3p,AM Zhu,TW}]
\label{Zhucenter}
The center of the Zhu-algebra $A(\W_{p_+,p_-})$ is generated by $[T]$ and isomorphic to
\begin{align*}
\C[x]/f_{p_+,p_-}(x),
\end{align*}
where
\begin{align*}
f_{p_+,p_-}(x)=&\prod_{(i,j)\in\mathcal{T}}(x-h_{r,s})^3\\
&\times\prod_{i=1}^{p_+-1}\prod_{j=1}^{p_--1}(x-\Delta^+_{i,j;0})^2\prod_{i=1}^{p_+-1}\prod_{j=1}^{p_--1}(x-\Delta^-_{i,j;0})\\
&\times\prod_{i=1}^{p_+-1}(x-\Delta^+_{i,p_-;0})^2\prod_{i=1}^{p_+-1}(x-\Delta^-_{i,p_-;0})\\
&\times\prod_{j=1}^{p_--1}(x-\Delta^+_{p_+,j;0})^2\prod_{j=1}^{p_--1}(x-\Delta^-_{p_+,j;0})\\
&\times(x-\Delta^+_{p_+,p_-;0})(x-\Delta^-_{p_+,p_-;0}).
\end{align*}
\end{thm}

The following corollary follows from the above theorem.
\begin{corollary}
\label{prop202302110}
The Zhu algebra $A(\W_{p_+,p_-})$ has three dimensional indecomposable modules on which $[T]$ acts as
\begin{equation*}
\begin{pmatrix}
h_{r,s} & 1 & 0\\
0 & h_{r,s} & 1\\
0& 0 & h_{r,s}\\
\end{pmatrix}
,
\end{equation*}
where $(r,s)\in \mathbb{K}_{p_+,p_-}$.
\end{corollary}
The above corollary will be used to determine the structure of the projective covers of the minimal simple modules in Subsection \ref{Projminimal}.

\subsection{Block decomposition of $\mathcal{C}_{p_+,p_-}$}
\label{W-modori}
Let $\mathcal{C}_{p_+,p_-}$ be the abelian category of grading restricted generalised $\W_{p_+,p_-}$-modules. Since $\W_{p_+,p_-}$ is $C_2$-cofinite, it follows from \cite{H} that any object $M$ in $\mathcal{C}_{p_+,p_-}$ has finite length. For any $M$ in $\mathcal{C}_{p_+,p_-}$, let $M^*$ denote the contragredient of $M$. From Proposition~\ref{simpleclass}, $\mathcal{C}_{p_+,p_-}$ is closed under contragredient (see also \cite{MS}).

\begin{dfn}
In the following, we define $\frac{(p_+-1)(p_--1)}{2}$ thick blocks, $p_++p_--2$ thin blocks and two semi-simple blocks.
\begin{enumerate}
\item For each $(r,s)\in\mathbb{K}_{p_+,p_-}$, we denote by $C^{thick}_{r,s}=C^{thick}_{p_+-r,p_--s}$ the full abelian subcategory of $\mathcal{C}_{p_+,p_-}$ such that
\begin{align*}
\hspace{-6mm}M\in C^{thick}_{r,s}&\Leftrightarrow\ {\rm all\ composition\ factors\ of}\ M\ {\rm are\ given\ by}\ \mathcal{X}^{+}_{r,s},\mathcal{X}^{+}_{r^{\vee},s^{\vee}},\\
&\ \ \ \ \X^-_{r^{\vee},s},\X^-_{r,s^{\vee}}\ {\rm and}\ L(h_{r,s}),
\end{align*}
where we used the notation $r^\vee=p_+-r$, $s^\vee=p_--s$ defined by (\ref{notrs}).
\item For each $1\leq s< p_-$, we denote by $C^{thin}_{p_+,s}$ the full abelian subcategory of $\mathcal{C}_{p_+,p_-}$ such that
\begin{align*}
\hspace{-6mm}M\in C^{thin}_{p_+,s}\Leftrightarrow\ {\rm all\ composition\ factors\ of}\ M\ {\rm are\ given\ by}\ \mathcal{X}^{+}_{p_+,s}\ {\rm and }\ \mathcal{X}^{-}_{p_+,s^{\vee}}.
\end{align*}
\item For each $1\leq r< p_+$, we denote by $C^{thin}_{r,p_-}$ the full abelian subcategory of $\mathcal{C}_{p_+,p_-}$ such that
\begin{align*}
\hspace{-6mm}M\in C^{thin}_{r,p_-}\Leftrightarrow\ {\rm all\ composition\ factors\ of}\ M\ {\rm are\ given\ by}\ \mathcal{X}^{+}_{r,p_-}\ {\rm and }\ \mathcal{X}^{-}_{r^{\vee},p_-}.
\end{align*}
\item We denote by $C^{\pm}_{p_+,p_-}$ the full abelian subcategory of $\mathcal{C}_{p_+,p_-}$ such that
\begin{align*}
\hspace{-6mm}M\in C^{\pm}_{p_+,p_-}\Leftrightarrow\ {\rm all\ composition\ factors\ of}\ M\ {\rm are\ given\ by}\ \mathcal{X}^{\pm}_{p_+,p_-}.
\end{align*}
\end{enumerate}
\end{dfn}
By Proposition~\ref{simpleclass}, every $\mathcal{W}_{p_+,p_-}$-module is contained in a suitable sum of the above blocks.

In the following, when the central charge is fixed to be $c_{p_+,p_-}$, we write the action of $\mathfrak{L}_n$ simply as $L_n$.
\begin{dfn}
Let $\mathcal{L}_{c_{p_+,p_-}}\mathchar`-{\rm Mod}$ be the abelian category of left generalized $U(\mathcal{L})$-modules whose morphisms are Virasoro-homomorphisms and whose objects are left $U(\mathcal{L})$-modules that satisfy the following conditions:
\begin{enumerate}
\item For the central charge, $C=c_{p_+,p_-}\cdot {\rm id}$ on $M$.
\item Every object $M$ has the following decomposition $M=\sum_{h\in{H(M)}}M[h]${\rm :}
\begin{itemize}
\item For some finite subset $H_0(M)$ of $\C$, $H(M)=H_0(M)+\Z_{\geq 0}$.
\item For $h\in H(M)$, $M[h]=\{m\in M: \exists n\geq 0\ {\rm s.t.}\ (L_0-h)^nm=0 \}$.
\item $0<{\rm dim}_{\C}M[h]<\infty$.
\end{itemize}
\item For every object $M\in \mathcal{L}_{c_{p_+,p_-}}\mathchar`-{\rm Mod}$, there exists the contragredient object $M^*\in \mathcal{L}_{c_{p_+,p_-}}\mathchar`-{\rm Mod}$ on which the anti-involution $\sigma(L_n)=L_{-n}$ induces the structure of a left $U(\mathcal{L})$-module by
\begin{align*}
\langle L_n\phi,u\rangle=\langle\phi,\sigma(L_n)u\rangle,\ \ \phi\in M^*,\ u\in M,
\end{align*} 
where $\langle-,-\rangle$ is the canonical pairing between $M^*$ and $M$.
\end{enumerate}
\end{dfn}
We call a module $M\in \mathcal{L}_{c_{p_+,p_-}}\mathchar`-{\rm Mod}$ on which $L_0$ acts semisimply a non-logarithmic Virasoro module, and $N\in \mathcal{L}_{c_{p_+,p_-}}\mathchar`-{\rm Mod}$ on which $L_0$ acts non-semisimply a logarithmic Virasoro module.
Let $\mathcal{L}^s_{c_{p_+,p_-}}\mathchar`-{\rm Mod}$ be the full subcategory of $\mathcal{L}_{c_{p_+,p_-}}\mathchar`-{\rm Mod}$ consisting of non-logarithmic modules.

Let $M(h,c_{p_+,p_-})$ be the Verma module of the Virasoro algebra whose lowest weight is $h\in\mathbb{C}$ and the central charge $C=c_{p_+,p_-}\cdot {\rm id}$. Let ${\mid}h\rangle$ be the lowest weight vector of $M(h,c_{p_+,p_-})$.
It is known that, for $r,s\geq 1$, $M(h_{r,s},c_{p_+,p_-})$ has the singular vector whose $L_0$-weight is $h_{r,s}+rs$ (for the notation $h_{r,s}$, see (\ref{eq:not-hrs})). 
We denote this singular vector by $S_{r,s}{\mid}h_{r,s}\rangle$, where $S_{r,s}\in U(\mathcal{L})$ is the Shapovalov element normalized as
\begin{align}
\label{eq:norm-s}
S_{r,s}-L^{rs}_{-1}\in {\rm span}_{\mathbb{C}}\{L_{-\lambda}\ |\ \lambda\vdash rs,\ \lambda\neq (1^{rs})\}.
\end{align}

We also introduce the following definition.
\begin{dfn}
Given an object $M\in \mathcal{L}^s_{c_{p_+,p_-}}\mathchar`-{\rm Mod}$, any nonzeo vector $v\in M$ is called a {\rm primary vector} when the following satisfied
\begin{align*}
L_nv=0,\ \ \ n\geq 1.
\end{align*} 
\end{dfn}
\begin{prop}
\label{FusionInt}
Let $M_1,M_2,M^*_3\in \mathcal{L}^s_{c_{p_+,p_-}}\mathchar`-{\rm Mod}$. Let $v_1\in M_{1}$, $v_2\in M_2$ and $v^*_3\in M^*_3$ be primary vectors whose $L_0$ weights are $h$, $h_{r_2,s_2}$ and $h_{r_3,s_3}$, respectively, where $r_2\geq 1$ and $s_2\geq 1$.
Let $\mathcal{Y}$ be a Virasoro-intertwining operator of type {\scriptsize{$\begin{pmatrix}
   M_3  \\
   M_1\ M_2
\end{pmatrix}
$}}.
Then we have
\begin{align*}
&\langle v^*_3,\mathcal{Y}(v_1,z)S_{r_2,s_2}v_2\rangle=z^{-r_2s_2}\prod_{i=1}^{r_2}\prod_{j=1}^{s_2}(h-h_{r_2+r_3-2i+1,s_2+s_3-2j+1})\langle v^*_3,\mathcal{Y}(v_1,z)v_2\rangle.
\end{align*}
\end{prop}
\begin{proof}
Note that the correlation function $\langle v^*_3,\mathcal{Y}(v_1,z)v_2\rangle$ is a scalar multiple of $z^{h_{r_3,s_3}-h_{r_2,s_2}-h}$.
Let $\langle v^*_3,\mathcal{Y}(v_1,z)v_2\rangle=cz^{h_{r_3,s_3}-h_{r_2,s_2}-h}$ $(c\in \mathbb{C})$. Then by the relation
\begin{align*}
[L_n,\mathcal{Y}(v_1,z)]=z^n(z\partial_z+h(n+1))\mathcal{Y}(v_1,z),\qquad n\in \mathbb{Z}.
\end{align*}
(cf. \cite{FHL}), there exists an $r_2s_2$-th order differential operator $D_z(h)$ with coefficients in $\mathbb{C}[h]$ satisfying the following property
\begin{align*}
\langle v^*_3,\mathcal{Y}(v_1,z)S_{r_2,s_2}v_2\rangle=D_z(h)\langle v^*_3,\mathcal{Y}(v_1,z)v_2\rangle=cP(h)z^{h_{r_3,s_3}-h_{r_2,s_2}-h-r_2s_2},
\end{align*}
where $P(h)$ is a polynomial of degree $r_2s_2$. 
In what follows, based on the technique in \cite[Section 3.9]{Yam}, we determine the zeros of $P(h)$.
Following \cite{Na0,TW}, we define $t$-deformations $\alpha_\pm$, $\alpha_{r,s}$, $h_{r,s}$, and $Q_\pm(x)$ as follows
\begin{align*}
&\widetilde{\alpha}_+=\alpha_++t,
&&\widetilde{\alpha}_-=\frac{-2}{\widetilde{\alpha}_+},\\
&\widetilde{\alpha}_{r,s}=\frac{1-r}{2}\widetilde{\alpha}_++\frac{1-s}{2}\widetilde{\alpha}_-,&&\widetilde{Q}_\pm(x)=V_{\widetilde{\alpha}_\pm},\\
&\widetilde{h}_{r,s}=\frac{1}{2}\widetilde{\alpha}_{r,s}(\widetilde{\alpha}_{r,s}-\widetilde{\alpha}_+-\widetilde{\alpha}_-).
\end{align*}
In the following, we write $\widetilde{L}_n$ for the $n$-th mode of the Virasoro algebra at central charge
$
1-3(\widetilde{\alpha}_++\widetilde{\alpha}_-)^2
$. By definition, $1-3(\widetilde{\alpha}_++\widetilde{\alpha}_-)^2|_{t=0}=c_{p_+,p_-}$.
Then we may regard $\widetilde{L}_n$ as a $t$-deformation of $L_n$.
We can check that $\widetilde{Q}_+(x)$ and $\widetilde{Q}_-(x)$ have conformal weight $1$ with respect to the $t$-deformed energy-momentum tensor $\widetilde{T}(z):=\sum_{n\in \mathbb{Z}}\widetilde{L}_nz^{-n-2}$.
That is,
\begin{align}
\label{eq:t-defqt}
\widetilde{T}(z)\widetilde{Q}_{\pm}(x)=\partial_{x}\bigl(\frac{\widetilde{Q}_{\pm}(x)}{z-x}\bigr)+\cdots.
\end{align}
Let $z\in \mathbb{R}_{>0}$.
%Following \cite{DF1,DF2,FS}, 
For $i,j\in \mathbb{Z}_{\geq 1}$ and $r_1,s_1\in \mathbb{Z}$, we consider the following correlation function, which can be expressed by means of the Dotsenko-Fateev integral introduced in \cite{DF1,DF2} (see also \cite{FS}):
\begin{align*}
&\Psi^0_{i,j}(z)\\
&=\int_{\Box^z_{i,j}}\bra{\widetilde{\alpha}_{-r_1+r_2-1-2i,-s_1+s_2-1-2j}+t}\widetilde{Q}^i_+(\bm{x})\widetilde{Q}^j_-(\bm{x}')V_{\widetilde{\alpha}_{-r_1,-s_1}+t}(z)\ket{\widetilde{\alpha}_{r_2,s_2}}{\rm d}\bm{x}{\rm d}\bm{x}',
\end{align*}
where we set
\begin{align*}
&\widetilde{Q}^i_+(\bm{x})=\widetilde{Q}_+(x_1)\cdots \widetilde{Q}_+(x_i)&\widetilde{Q}^j_-(\bm{x}')=\widetilde{Q}_-(x'_1)\cdots \widetilde{Q}_-(x'_j),
\end{align*}
and $\Box^z_{i,j}$ denotes the regularized cycle of the product 
$
\lbrack0,z\rbrack^i_{x_1,\dots,x_i}\times \lbrack0,z\rbrack^j_{x'_1,\dots,x'_j}
$
constructed in \cite{S1}. 
As constructed in detail in \cite{S1}, the cycle $\Box^z_{i,j}$ is a compact boundaryless cycle in the monodromy cover of the manifold
\begin{align*}
\{(x_1,\dots,x_i,x'_1,\dots,x'_j)\in \mathbb{C}^{i+j}\ | x_k\neq x'_l,\ x_m\neq x_n,\ x'_m\neq x'_n,\ m\neq n \}.
\end{align*}
We also introduce the following correlation functions, which can be expressed by means of the Selberg integral \cite{Se}:
\begin{align*}
\Psi^+_{m}(z)&=\int_{\Delta^{z;+}_m}\bra{\widetilde{\alpha}_{-r_1+r_2-1-2m,-s_1+s_2-1}+t}\widetilde{Q}^m_+(\bm{x})V_{\widetilde{\alpha}_{-r_1,-s_1}+t}(z)\ket{\widetilde{\alpha}_{r_2,s_2}}{\rm d}\bm{x},\\
\Psi^-_{m}(z)&=\int_{\Delta^{z;-}_m}\bra{\widetilde{\alpha}_{-r_1+r_2-1,-s_1+s_2-1-2m}+t}\widetilde{Q}^m_-(\bm{x})V_{\widetilde{\alpha}_{-r_1,-s_1}+t}(z)\ket{\widetilde{\alpha}_{r_2,s_2}}{\rm d}\bm{x},
\end{align*}
where $\Delta^{z;\pm}_m$ denotes the twisted cycle constructed from the simplex 
$
\{\bm{x}\in \mathbb{R}^m|z>x_1>\cdots >x_m>0\}
$
according to the multivaluedness of the integrand by applying the construction of \cite{TK}.
We regard each of $F_{\widetilde{\alpha}_{r_2,s_2}}$, $F_{\widetilde{\alpha}_{-r_1,-s_1}+t}$ and $F^\vee_{\widetilde{\alpha}_{-r_3,-s_3}+t}$ as a representation of the Virasoro algebra generated by the modes $\widetilde{L}_n$ of $\widetilde{T}(z)$. Note that, due to the shift of $t$, $F_{\widetilde{\alpha}_{-r_1,-s_1}+t}$ is simple as a Virasoro module.
We denote by $\widetilde{\alpha}_{-r_3,-s_3}+t$ the $a_0$-weight in the bra vectors of $\Psi^0_{i,j}(z)$ and $\Psi^\pm_{m}(z)$.
Then, by using (\ref{eq:t-defqt}) and the Stokes' theorem, we have
\begin{align}
\label{eq:sol-1}
\begin{aligned}
&\int_{\Box^z_{i,j}}\bra{\widetilde{\alpha}_{-r_3,-s_3}+t}\bigl[\widetilde{L}_n,\widetilde{Q}^i_+(\bm{x})\widetilde{Q}^j_-(\bm{x}')V_{\widetilde{\alpha}_{-r_1,-s_1}+t}(z)\bigr]\ket{\widetilde{\alpha}_{r_2,s_2}}{\rm d}\bm{x}{\rm d}\bm{x}'\\
&=z^n(z\partial_z+\widehat{h}_{r_1,s_1}(n+1))\Psi_{i,j}(z),
\end{aligned}
\end{align}
where we set
\begin{align*}
\widehat{h}_{r_1,s_1}:=\frac{1}{2}(\widetilde{\alpha}_{-r_1,-s_1}+t)(\widetilde{\alpha}_{-r_1,-s_1}+t-\widetilde{\alpha}_+-\widetilde{\alpha}_-).
\end{align*}
In particular, we have
$
\Psi^0_{i,j}(z)=\Psi^0_{i,j}(1)z^{\widetilde{h}_{r_3,s_3}-\widetilde{h}_{r_2,s_2}-\widehat{h}_{r_1,s_1}}.
$
Also, we see that 
\begin{align}
\label{eq:sol-1-1}
\begin{aligned}
&\int_{\Delta^{z;\pm}_m}\bra{\widetilde{\alpha}_{-r_3,-s_3}+t}\bigl[\widetilde{L}_n,\widetilde{Q}^m_\pm(\bm{x})V_{\widetilde{\alpha}_{-r_1,-s_1}+t}(z)\bigr]\ket{\widetilde{\alpha}_{r_2,s_2}}{\rm d}\bm{x}\\
&=z^n(z\partial_z+\widehat{h}_{r_1,s_1}(n+1))\Psi^\pm_{m}(z).
\end{aligned}
\end{align}
In particular, $\Psi^\pm_{m}(z)=\Psi^\pm_{m}(1)z^{\widetilde{h}_{r_3,s_3}-\widetilde{h}_{r_2,s_2}-\widehat{h}_{r_1,s_1}}$.
For $i,j\in \mathbb{Z}_{\geq 1}$, the following explicit formula for the constant $\Psi^0_{i,j}(1)$ is known \cite{DF2,Fo}:
\begin{align*}
\frac{\Psi^0_{i,j}(1)}{i!j!\widetilde{\rho}^{2ij}_-}=&\prod^{i-1}_{k=0}\frac{\sin \pi \bigl(\widetilde{a}_++\widetilde{b}_++2+(i-1+k)\widetilde{\rho}_+\bigr)}{\sin \pi (\widetilde{a}_++1+k\widetilde{\rho}_+)}\prod^{i}_{k=1}\frac{\Gamma(k\widetilde{\rho}_+)}{\Gamma(\widetilde{\rho}_+)}\prod^j_{l=1}\frac{\Gamma(l\widetilde{\rho}_--i)}{\Gamma(\widetilde{\rho}_-)}\\
&\times \prod_{k=0}^{i-1}\frac{\Gamma(1+\widetilde{b}_++k\widetilde{\rho}_+)\Gamma(-1+2j-\widetilde{a}_+-\widetilde{b}_+-(i-1+k)\widetilde{\rho}_+)}{\Gamma(-\widetilde{a}_+-k\widetilde{\rho}_+)}\\
&\times \prod_{l=0}^{j-1}\frac{\Gamma(1-i+\widetilde{a}_-+l\widetilde{\rho}_-)\Gamma(1-i+\widetilde{b}_-+l\widetilde{\rho}_-)}{\Gamma(2-i+\widetilde{a}_-+\widetilde{b}_-+(j-1+l)\widetilde{\rho}_-)}
\end{align*}
where we set
\begin{align*}
\widetilde{a}_\pm=\widetilde{\alpha}_\pm\widetilde{\alpha}_{r_2,s_2},\qquad \widetilde{b}_\pm=\widetilde{\alpha}_\pm(\widetilde{\alpha}_{-r_1,-s_1}+t),\qquad \widetilde{\rho}_\pm=\frac{\widetilde{\alpha}^2_\pm}{2}.
\end{align*}
It is also known that \cite{Se,TK}
\begin{align*}
\Psi^\pm_{m}(1)=\frac{1}{m!}\prod_{k=1}^m\frac{\Gamma(1+k\widetilde{\rho}_\pm)\Gamma(1+\widetilde{a}_\pm+(k-1)\widetilde{\rho}_\pm)\Gamma(1+\widetilde{b}_\pm+(k-1)\widetilde{\rho}_\pm)}{\Gamma(1+\widetilde{\rho}_\pm)\Gamma(2+\widetilde{a}_\pm+\widetilde{b}_\pm+(m+k-2)\widetilde{\rho}_\pm)}.
\end{align*}
From these formulas, we see that 
\begin{align}
\label{eq:sol-2}
\begin{aligned}
&\Psi^0_{i,j}(z)\neq 0,\quad r_2-1\geq i\geq 1,\ s_2-1\geq j\geq 1,\\
&\Psi^+_i(z)\neq 0,\ \ \Psi^-_j(z)\neq 0, \quad r_2-1\geq i\geq 0,\ s_2-1\geq j\geq 0.
\end{aligned}
\end{align}
It is known that there exists $\widetilde{S}_{r_2,s_2}\in \mathbb{C}(\widetilde{\alpha}_+)\otimes U(\mathcal{L})$ such that
\begin{align}
\label{eq:sol-3}
&\widetilde{S}_{r_2,s_2}\ket{\widetilde{\alpha}_{r_2,s_2}}=0,
&\widetilde{S}_{r_2,s_2}|_{t=0}={S}_{r_2,s_2}
\end{align}
(cf. \cite{FF,IK}).
We again set
\begin{align*}
&-r_3=-r_1+r_2-1-2i,
&-s_3=-s_1+s_2-1-2j.
\end{align*}
Then, by (\ref{eq:sol-1})--(\ref{eq:sol-3}), there exists an $r_2s_2$-th order differential operator $\widetilde{D}_z$ with coefficients in $\mathbb{C}(t)$ satisfying
\begin{align}
\label{eq:dz}
\widetilde{D}_zz^{\widetilde{h}_{r_3,s_3}-\widetilde{h}_{r_2,s_2}-\widehat{h}_{r_1,s_1}}=0.
\end{align}
We see that all the coefficients of $\widetilde{D}_z$ are holomorphic at $t=0$.
Further, since $\widetilde{S}_{r_2,s_2}|_{t=0}={S}_{r_2,s_2}$, we have $\widetilde{D}_z|_{t=0}=D_z(h_{r_1,s_1})$.
Then, noting (\ref{eq:dz}), $\widetilde{h}_{r,s}|_{t=0}=h_{r,s}$ and $\widehat{h}_{r_1,s_1}|_{t=0}=h_{r_1,s_1}$, we see that
\begin{align*}
h=h_{r_2+r_3-2i-1,s_2+s_3-2j-1},\quad r_2-1\geq i\geq 0,\ s_2-1\geq j\geq 0
\end{align*}
are zeros of $P(h)$. Thus, we obtain the claimed equality.
\end{proof}

From the above proposition and \cite[Proposition~5.4.7,Theorem~5.5.1]{FHL}, we obtain the following proposition (see also \cite{IK,Milas,Lin}).
\begin{prop}
\label{VirasoroFusion}
For $h\in\C$, $1\leq r_1,r_2<p_+$, $1\leq s_1,s_2<p_-$ and $n_1,n_2\in\Z_{\geq 0}$, we have
\begin{equation*}
\mathcal{N}^{L(h)}_{L(h_{r_1,s_1;n_1}),L(h_{r_2,s_2;n_2})}\leq 1,
\end{equation*}
where $\mathcal{N}^{L(h_3)}_{L(h_2),L(h_1)}$ is the dimension of the space of Virasoro intertwining operators of type {\scriptsize{$\begin{pmatrix}
   L(h_3)  \\
   L(h_2)\ L(h_1)
\end{pmatrix}
$}}.
If $\mathcal{N}^{L(h)}_{L(h_{r_1,s_1;n_1}),L(h_{r_2,s_2;n_2})}\neq 0$, then $h$ is the common solution of the following equations
\begin{align*}
\prod_{i=1}^{r_1}\prod_{j=1}^{s_1+n_1p_-}&(h-h_{r_1+r_2-2i+1,s_1+s_2-2j+1;n_1+n_2})=0,\\
\prod_{i=1}^{(n_1+1)p_+-r_1}\prod_{j=1}^{p_--s_1}&(h-h_{2p_+-r_1-r_2-2i+1,2p_--s_1-s_2-2j+1;-n_1-n_2})=0,\\
\prod_{i=1}^{r_2}\prod_{j=1}^{s_2+n_2p_-}&(h-h_{r_1+r_2-2i+1,s_1+s_2-2j+1;n_1+n_2})=0,\\
\prod_{i=1}^{(n_2+1)p_+-r_2}\prod_{j=1}^{p_--s_2}&(h-h_{2p_+-r_1-r_2-2i+1,2p_--s_1-s_2-2j+1;-n_1-n_2})=0.
\end{align*}
\end{prop}

\begin{thm}
\label{prop:block}
The abelian category $\mathcal{C}_{p_+,p_-}$ has the following block decomposition
\begin{align*}
\mathcal{C}_{p_+,p_-}=\bigoplus_{(r,s)\in\mathbb{K}_{p_+,p_-}}C^{thick}_{r,s}\oplus \bigoplus_{r=1}^{p_+-1}C^{thin}_{r,p_-}\oplus\bigoplus_{s=1}^{p_--1}C^{thin}_{p_+,s}\oplus C^+_{p_+,p_-}\oplus C^-_{p_+,p_-}.
\end{align*}
\end{thm}
\begin{proof}
The proof is analogous to that in \cite[Theorem~4.4]{AM}. 
Here, we give a proof using Proposition~\ref{FusionInt}.

Let $M$ and $N$ be non-logarithmic objects belonging to one of the blocks. The $\mathcal{W}_{p_+,p_-}$-action defines
a Virasoro-intertwining operator $\mathcal{Y}$ of type {\scriptsize{$\begin{pmatrix}
   M  \\
   N\ L
\end{pmatrix}
$}}, where $L=L(h_{4p_+-1,1})$.
Take arbitrary primary vectors $u^*\in M^*$ and $v\in N$, and let $w$ be the lowest weight vector of $L(h_{4p_+-1,1})$.
Let $h_{r,s}$ and $h_v$ be the $L_0$-weight of $u^*$ and $v$, respectively. Note that $S_{4p_+-1,1}w=S_{1,4p_--1}w=0$. Then, by Proposition~\ref{FusionInt}, we have
\begin{align*}
0=&\langle u^*,\mathcal{Y}(v,z)S_{4p_+-1,1}w\rangle=z^{1-4p_+}\prod_{i=1}^{4p_+-1}(h_{v}-h_{r+4p_+-1-2i+1,s})\langle u^*,\mathcal{Y}(v,z)w\rangle,\\
0=&\langle u^*,\mathcal{Y}(v,z)S_{1,4p_--1}w\rangle=z^{1-4p_-}\prod_{j=1}^{4p_--1}(h_{v}-h_{r,s+4p_--1-2j+1})\langle u^*,\mathcal{Y}(v,z)w\rangle.
\end{align*}
Assume that $\langle u^*,\mathcal{Y}(v,z)w\rangle\neq 0$. Then from the above equalities, we have
\begin{align*}
h_v\in \{h_{r,s},h_{r,s;-2}, h_{r,s;2},h_{-r,s}, h_{r^\vee,s;-1}\}.
\end{align*}
This implies that $N$ must be an object in the same block as $M$.
\end{proof}
By this theorem, in order to compute the ${\rm Ext}^1$-groups in $\mathcal{C}_{p_+,p_-}$, it suffices to study each block separately.

\section{Logarithmic $\W_{p_+,p_-}$ modules}
\label{logarithmic section}
In what follows, we call a module $M\in \mathcal{C}_{{p_+,p_-}}$ on which $L_0$ acts semisimply a non-logarithmic $\mathcal{W}_{p_+,p_-}$-module, and $N\in \mathcal{C}_{{p_+,p_-}}$ on which $L_0$ acts non-semisimply a logarithmic $\mathcal{W}_{p_+,p_-}$-module.
In this section, by using the logarithmic deformation by \cite{FJ}, we construct certain logarithmic $\W_{p_+,p_-}$-modules which correspond to the projective covers of all simple $\W_{p_+,p_-}$-modules $\X^\pm_{\bullet,\bullet}$ in the thick blocks and the thin blocks, and we introduce indecomposable modules $\mathcal{Q}(\X^\pm_{\bullet,\bullet})_{\bullet,\bullet}$ which become important after this section.
These logarithmic modules are closely related to certain indecomposable modules of the quantum group $\mathfrak{g}_{p_+,p_-}$ at the root of unity $e^{\frac{\pi i}{2p_+p_-}}$ \cite{Arike,FF22}. In the one thick block $C^{thick}_{1,1}$, these logarithmic modules are constructed by \cite{AK}.
\subsection{Logarithmic deformation}
\label{Logarithmic deformation}
\begin{prop}
\label{Hom1}
For $r,s\geq 1$, we have the following relation
\begin{align*}
\alpha_-[Q^{[r]}_+,Q^{[s]}_-(z)]=\alpha_+[Q^{[s]}_-,Q^{[r]}_+(z)].
\end{align*}
\end{prop}
\begin{proof}
In this proof, we will use the argument in \cite[Proposition~2.1]{NT}.
Recall the defining equation (\ref{Tsuchiya-Kanie}) of $Q^{[\bullet]}_\pm(z)$.
Then we have
\begin{align}
&{\rm Res}_{z=w}Q^{[r]}_+(z)Q^{[s]}_-(w){\rm d}z\nonumber\\
&={\rm Res}_{z=w}\int_{\Delta^+_{r-1}}{\rm d}\bm{x}\int_{\Delta^-_{s-1}}{\rm d}\bm{y}\frac{1}{(z-w)^2}:e^{\alpha_+\phi(z)+\alpha_-\phi(w)}:\nonumber\\
&\qquad \qquad \qquad\times Q(zx_1)\cdots Q_+(zx_{r-1})Q_-(wy_1)\cdots Q_-(wy_{s-1})z^{r-1}w^{s-1}\nonumber\\
&=\int_{\Delta^+_{r-1}}{\rm d}\bm{x}\int_{\Delta^-_{s-1}}{\rm d}\bm{y}\frac{\alpha_+}{\alpha_++\alpha_-}\Bigl(\frac{\partial}{\partial w}V_{\alpha_++\alpha_-}(w)\Bigr)\nonumber\\
&\qquad \qquad \qquad\times Q(wx_1)\cdots Q_+(wx_{r-1})Q_-(wy_1)\cdots Q_-(wy_{s-1})w^{r+s-2}\nonumber\\
&+\int_{\Delta^+_{r-1}}{\rm d}\bm{x}\int_{\Delta^-_{s-1}}{\rm d}\bm{y}z^{r-1}V_{\alpha_++\alpha_-}(z)\nonumber\\
&\qquad \qquad \qquad\times \left.\frac{\partial}{\partial z}\bigl(Q_{+}(zx_1)\cdots Q_+(zx_{r-1})\bigr)Q_-(wy_1)\cdots Q_-(wy_{s-1}) w^{s-1}\right|_{z=w}\nonumber\\
&+(r-1)\int_{\Delta^+_{r-1}}{\rm d}\bm{x}\int_{\Delta^-_{s-1}}{\rm d}\bm{y}\nonumber\\
&\qquad \qquad \times V_{\alpha_++\alpha_-}(w)Q_+(wx_1)\cdots Q_+(wx_{r-1})Q_-(wy_1)\cdots Q_-(wy_{s-1})w^{r+s-3}.
\label{ast}
\end{align}
Since
\begin{align*}
&\frac{\partial}{\partial z}\bigl(Q_{+}(zx_1)\cdots Q_+(zx_{r-1})\\
&=\sum_{i=1}^{r-1}Q_+(zx_1)\cdots \frac{\partial}{\partial z}Q_+(zx_i)\cdots Q_{+}(zx_{r-1})\\
&=\sum_{i=1}^{r-1}Q_+(zx_1)\cdots \Bigl(\frac{1}{z}x_i\frac{\partial}{\partial x_i}Q_+(zx_i)\Bigr)\cdots Q_{+}(zx_{r-1})\\
&=\frac{1}{z}\sum_{i=1}^{r-1}Q_+(zx_1)\cdots \Bigl(\frac{\partial}{\partial x_i}x_iQ_+(zx_i)\Bigr)\cdots Q_{+}(zx_{r-1})\\
&\ \ \ -\frac{r-1}{z}Q_+(zx_1)\cdots Q_+(zx_{r-1}),
\end{align*}
the second term of the right-hand side of $(\ref{ast})$ becomes
\begin{equation}
\label{onaji2023}
\begin{split}
&-(r-1)\int_{\Delta^{+}_{r-1}}{\rm d}\bm{x}\int_{\Delta^{-}_{s-1}}{\rm d}\bm{y}\\
&\qquad \qquad \times V_{\alpha_++\alpha_-}(w)Q_+(wx_1)\cdots Q_+(wx_{r-1})Q_-(wy_1)\cdots Q_-(wy_{s-1})w^{r+s-3}\\
&+\int_{\Delta^{-}_{s-1}}{\rm d}\bm{y}\int_{\Delta^{+}_{r-1}}V_{\alpha_++\alpha_-}(w)\\
&\qquad \qquad\times {\rm d}_{\bm{x}}\Bigl(\sum_{i=1}^{r-1}(-1)^{i+1}x_iQ_+(wx_1)\cdots Q_+(wx_{r-1}) {\rm d}x_1\wedge \cdots\widehat{{\rm d}x_i}\cdots\wedge{\rm d}x_{r-1}\Bigr)\\
&\qquad \qquad\times Q_-(wy_1)\cdots Q_-(wy_{s-1}),
\end{split}
\end{equation}
where ${\rm d}_{\bm{x}}$ is the total derivative with respect to the $x_i$ variables.
Since $\Delta^+_{r-1}$ is a twisted cycle, it follows from the generalized Stokes' theorem (cf.~\cite{AoK}) that the second term of (\ref{onaji2023}) vanishes. 
Note that the first term of (\ref{onaji2023}) cancels with the third term of $(\ref{ast})$.
Thus $[Q^{[r]}_+,Q^{[s]}_-(w)]$ becomes 
\begin{align*}
&\frac{\alpha_+}{\alpha_++\alpha_-}\int_{\Delta^{+}_{r-1}}{\rm d}\bm{x}\int_{\Delta^{-}_{s-1}}{\rm d}\bm{y}\\
&\qquad \times\Bigl(\frac{\partial}{\partial w}V_{\alpha_++\alpha_-}(w)\Bigr)Q(wx_1)\cdots Q_+(wx_{r-1})Q_-(wy_1)\cdots Q_-(wy_{s-1})w^{r+s-2}.
\end{align*}
In the same way, we have
\begin{align*}
&[Q^{[s]}_-,Q^{[r]}_+(w)]\\
&=\frac{\alpha_-}{\alpha_++\alpha_-}\int_{\Delta^{+}_{r-1}}{\rm d}\bm{x}\int_{\Delta^{-}_{s-1}}{\rm d}\bm{y}\\
&\qquad \times\Bigl(\frac{\partial}{\partial w}V_{\alpha_++\alpha_-}(w)\Bigr)Q_+(wx_1)\cdots Q_+(wx_{r-1})Q_-(wy_1)\cdots Q_-(wy_{s-1})w^{r+s-2}.
\end{align*}
Therefore we obtain
$
\alpha_-[Q^{[r]}_+,Q^{[s]}_-(z)]=\alpha_+[Q^{[s]}_-,Q^{[r]}_+(z)].
$
\end{proof}

\begin{prop}
For $r,s\geq 1$ the screening operators $Q^{[r]}_+$ and $Q^{[s]}_-$ are $\W_{p_+,p_-}$-homomorphisms. 
\begin{proof}
For each of the triplet generators of $\W_{p_+,p_-}$, we have the following two expressions
\begin{align*}
&W^+=Q^{[p_--1]}_-\ket{\alpha_{1,p_--1;3}}=Q^{[3p_+-1]}_+\ket{\alpha_{p_+-1,1;-3}},\\
&W^-=Q^{[p_+-1]}_+\ket{\alpha_{p_+-1,1;-3}}=Q^{[3p_--1]}_-\ket{\alpha_{1,p_--1;3}},\\
&W^0=Q^{[2p_+-1]}_+\ket{\alpha_{p_+-1,1;-3}}=Q^{[2p_--1]}_-\ket{\alpha_{1,p_--1;3}},
\end{align*}
up to scalar multiples.
Note that the products $Q_+(z)Y(\ket{\alpha_{1,p_--1;3}};w)$ and $Q_-(z)Y(\ket{\alpha_{p_+-1,1;-3}};w)$ have no poles at $z=w$ (see (\ref{eq:vavb})).
By the proof of Proposition \ref{Hom1}, we see that $[Q^{k}_+,Q^{[l]}_-]=0$ for $k,l\geq 1$ (cf. \cite[Proposition~2.10]{Na0}).
Thus, from the above two expressions of each of the triplet generators, we obtain 
\begin{align*}
[Q^{[r]}_+,Y(W^a;z)]=[Q^{[s]}_-,Y(W^a;z)]=0
\end{align*}
for $a=+,-,0$. This implies that $Q^{[r]}_+$ and $Q^{[s]}_-$ are $\W_{p_+,p_-}$-homomorphisms.
\end{proof}
\label{Hom2}
\end{prop}

We introduce the following logarithmic deformation by \cite{FJ}.
\begin{dfn}[\cite{FJ}]
\begin{enumerate}
\item Given mutually local fields $E(z)$ and $A(z)$, we define the logarithmic deformation of $A(z)$ by $E(z)$ as follows
\begin{align*}
\Delta_{E}(A(z))=\log z(E[0]A)(z)+\sum_{n\geq 1}\frac{(-1)^{n+1}}{n}\frac{(E[n]A)(z)}{z^n},
\end{align*}
where
\begin{align*}
(E[n]A)(w)=\oint_{z=w}(z-w)^nE(z)A(w){\rm d}z.
\end{align*}
\item Given mutually local fields $E(z)$, $A(z)$ and $B(z)$, we define 
\begin{align*}
\Delta_{E}(A(z)B(w))=\sum_{n\in\Z}\frac{\Delta_E((A[n]B)(w))}{(z-w)^{n+1}}.
\end{align*}
\end{enumerate}
\end{dfn}
\begin{thm}[\cite{FJ}]
\label{FJ}
Let $E(z)$, $A(z)$ and $B(z)$ be any mutually local fields.
Then the operator $\Delta_E$ satisfies the following derivation property
\begin{align*}
\Delta_E(A(z)B(w))=\Delta_E(A(z))B(w)+A(z)\Delta_E(B(w)).
\end{align*}
\end{thm}
In the next subsection, we consider the logarithmic deformations by the screening currents $Q^{[\bullet]}_{\pm}(z)$. We set
\begin{align*}
&\Delta^{[r]}_{+}:=\Delta_{Q^{[r]}_+},
&\Delta^{[s]}_{-}:=\Delta_{Q^{[s]}_-}.
\end{align*}
For the energy-momentum tensor, we have
\begin{align}
\label{202210130}
&\Delta^{[r]}_{+}(T(z))=T(z)+\frac{Q^{[r]}_+(z)}{z},
&\Delta^{[s]}_{-}(T(z))=T(z)+\frac{Q^{[s]}_-(z)}{z}.
\end{align}
By Proposition \ref{Hom2}, for any $A\in \W_{p_+,p_-}$, the field $\Delta^{[\bullet]}_{\pm}(Y(A;z))$ contains no $\log z$ terms. 

\subsection{Logarithmic modules in the thick block}
Let $r$ and $s$ be integers satisfying $1\leq r< p_+$, $1\leq s< p_-$. We set
\begin{align*}
\mathcal{P}_{r,s}=\V^+_{r,s}\oplus \V^{+}_{r^{\vee},s^{\vee}}\oplus \V^{-}_{r,s^{\vee}}\oplus \V^-_{r^{\vee},s},
\end{align*}
where $r^\vee=p_+-r$ and $s^\vee=p_--s$. Note that $\V^+_{r,s}, \V^{+}_{r^{\vee},s^{\vee}}, \V^{-}_{r,s^{\vee}}, \V^-_{r^{\vee},s}\in C^{thick}_{r,s}$. Let $(\mathcal{P}_{r,s},Y_{\mathcal{P}_{r,s}})$ be the ordinary $\W_{p_+,p_-}$-module. 
Fix an arbitrary $\tau=(a,b,\epsilon)$ in
\begin{align*}
\{(a,b,\epsilon)\}=\{(r,s,+),(r^\vee,s^\vee,+),(r^\vee,s,-),(r,s^\vee,-)\}.
\end{align*}
For $A\in \W_{p_+,p_-}$, we define the following operators on $\mathcal{P}_{r,s}$ :
\begin{align*}
\begin{aligned}
&\widetilde{\Delta}_{\tau}(Y_{\mathcal{P}_{r,s}}(A;z))\\
&=
\begin{cases}
(\alpha_--\alpha_+)\bigl(\Delta^{[a]}_++\Delta^{[b]}_-\bigr)(Y_{\mathcal{P}_{r,s}}(A;z))\\
+\bigl(-\alpha_+\Delta^{[b]}_-\circ\Delta^{[a]}_++\alpha_-\Delta^{[a]}_+\circ\Delta^{[b]}_-\bigr)(Y_{\mathcal{P}_{r,s}}(A;z))& on\ \V^{\epsilon}_{a,b}\\
0 & on\ \mathcal{P}_{r,s}\setminus \V^\epsilon_{a,b},
\end{cases}
\end{aligned}
\end{align*}

\begin{equation*}
\Delta^{-}_{\tau}(Y_{\mathcal{P}_{r,s}}(A;z))
=
\begin{cases}
\Delta^{[b]}_-(Y_{\mathcal{P}_{r,s}}(A;z))& on\ \V^{-\epsilon}_{a^\vee,b}\\
0 & on\ \mathcal{P}_{r,s}\setminus \V^{-\epsilon}_{a^\vee,b},
\end{cases}
\end{equation*}

\begin{equation*}
\Delta^{+}_{\tau}(Y_{\mathcal{P}_{r,s}}(A;z))
=
\begin{cases}
\Delta^{[a]}_+(Y_{\mathcal{P}_{r,s}}(A;z))& on\ \V^{-\epsilon}_{a,b^\vee}\\
0 & on\ \mathcal{P}_{r,s}\setminus \V^{-\epsilon}_{a,b^\vee}.
\end{cases}
\end{equation*}

By the following lemma, we see that the operator $\widetilde{\Delta}_\tau$ contains no $\log z$ terms.
\begin{lem}
\label{logz}
Let $r$ and $s$ be integers satisfying $1\leq r< p_+$, $1\leq s< p_-$.
For any $A\in \W_{p_+,p_-}$, the operator
\begin{align*}
-\alpha_+\Delta^{[s]}_-\bigl(\Delta^{[r]}_+(Y(A;z))\bigr)+\alpha_-\Delta^{[r]}_+\bigl(\Delta^{[s]}_-(Y(A;z))\bigr)
\end{align*}
contains no $\log z$ terms.
\begin{proof}
The $\log z$ terms of $\Delta^{[s]}_-\bigl(\Delta^{[r]}_+(Y(A;z))\bigr)$ and $\Delta^{[r]}_+\bigl(\Delta^{[s]}_-(Y(A;z))$ are given by
\begin{align*}
[Q^{[s]}_-,\Delta^{[r]}_+(Y(A;z))]\log z,\ \ \ \ \ \ [Q^{[r]}_+,\Delta^{[s]}_-(Y(A;z))]\log z.
\end{align*}
By using Proposition \ref{Hom1} we have
\begin{align*}
&[Q^{[s]}_-,\Delta^{[r]}_+(Y(A;z))]\\
&=\Bigl[Q^{[s]}_-,\sum_{n\geq 1}\frac{(-1)^{n+1}}{n}\oint_{w=z}\frac{(w-z)^n}{z^n}Q^{[r]}_+(w)Y(A;z){\rm d}w\Bigr]\\
&=\frac{\alpha_-}{\alpha_+}\Bigl[Q^{[r]}_+,\sum_{n\geq 1}\frac{(-1)^{n+1}}{n}\oint_{w=z}\frac{(w-z)^n}{z^n}Q^{[s]}_-(w)Y(A;z){\rm d}w\Bigr]\\
&=\frac{\alpha_-}{\alpha_+}[Q^{[r]}_+,\Delta^{[s]}_-(Y(A;z))].
\end{align*}
Thus, the claim follows.
\end{proof}
\end{lem}

Motivated by the construction in \cite{AK}, we introduce rank three logarithmic $\mathcal{W}_{p_+,p_-}$-modules.
\begin{thm}
\label{logarithmic1}
Fix an arbitrary $\tau=(a,b,\epsilon)$ in 
\begin{align*}
\{(a,b,\epsilon)\}=\{(r,s,+),(r^\vee,s^\vee,+),(r^\vee,s,-),(r,s^\vee,-)\}.
%\{(a,b,c,d,\epsilon)\}=\{(r,s,r^\vee,s^\vee;+),(r^\vee,s^\vee,r,s,+),(r^\vee,s,r,s^\vee;-),(r,s^\vee,r^\vee,s,-)\}.
\end{align*}
We define a logarithmic $\W_{p_+,p_-}$-module $\mathcal{P}^\epsilon_{a^\vee,b^\vee}$ as follows.
As the vector space $\mathcal{P}^\epsilon_{a^\vee,b^\vee}=\mathcal{P}_{r,s}$ and the module actions is defined by
\begin{align*}
J^\epsilon_{a^\vee,b^\vee}(A;z)&=Y_{\mathcal{P}_{r,s}}(A;z)+\bigl(\widetilde{\Delta}_{\tau}+\Delta^{+}_{\tau}+\Delta^{-}_{\tau}\bigr)(Y_{\mathcal{P}_{r,s}}(A;z)),
\end{align*}
for any $A\in \W_{p_+,p_-}$.
\begin{proof}
By Lemma \ref{logz}, we have $J^\epsilon_{a^\vee,b^\vee}: \W_{p_+,p_-}\rightarrow {\rm End}\mathcal{P}^\epsilon_{a^\vee,b^\vee}[[z,z^{-1}]]$. The condition that $J^\epsilon_{a^\vee,b^\vee}(\ket{0};z)={\rm id}_{\mathcal{P}^\epsilon_{a^\vee,b^\vee}}$ is trivial from the definition of logarithmic deformation. In the following we prove the identity
\begin{align}
\label{eq:iden-yab}
J^{\epsilon}_{a^\vee,b^\vee}(A;z)J^\epsilon_{a^\vee,b^\vee}(B;w)=J^\epsilon_{a^\vee,b^\vee}(Y(A;z-w)B;w)
\end{align}
for $A,B\in \W_{p_+,p_-}$. Fix an arbitrary vector $v\in \mathcal{P}_{r,s}$ and write $v$ be as follows
\begin{align*}
v=k_1v^\epsilon_{a,b}+k_2v^{-\epsilon}_{a^\vee,b}+k_3v^{-\epsilon}_{a,b^\vee}+k_4v^\epsilon_{a^\vee,b^\vee},
\end{align*} 
where $v^\epsilon_{a,b}\in \V^\epsilon_{a,b},\ v^{-\epsilon}_{a^\vee,b}\in \V^{-\epsilon}_{a^\vee,b},\ v^{-\epsilon}_{a,b^\vee}\in \V^{-\epsilon}_{a,b^\vee},\ v^\epsilon_{a^\vee,b^\vee}\in \V^\epsilon_{a^\vee,b^\vee}$, and $k_i(i=1,\dots,4)$ are constants. 
By the definition of $J^{\epsilon}_{a^\vee,b^\vee}$,  we have
\begin{align*}
&J^{\epsilon}_{a^\vee,b^\vee}(A;z)J^{\epsilon}_{a^\vee,b^\vee}(B;z)v^\epsilon_{a,b}\\
&=Y(A;z)Y(B;w)v^\epsilon_{a,b}\\
&\ \ +(\alpha_--\alpha_+)\bigl[\Delta^{[a]}_+(Y(A;z))+\Delta^{[b]}_-(Y(A;z))\bigr]Y(B;w)v^\epsilon_{a,b}\\
&\ \ +\bigl[-\alpha_+\Delta^{[b]}_-(\Delta^{[a]}_+(Y(A;z)))+\alpha_-\Delta^{[a]}_+(\Delta^{[b]}_-(Y(A;z)))\bigr]Y(B;w)v^\epsilon_{a,b}\\
&\ \ +(\alpha_--\alpha_+)Y(A;z)\bigl[\Delta^{[a]}_+(Y(B;w))+\Delta^{[b]}_-(Y(B;w))\bigr]v^\epsilon_{a,b}\\
&\ \ +(\alpha_--\alpha_+)\bigl[\Delta^{[b]}_-(Y(A;z))\Delta^{[a]}_+(Y(B;w))+\Delta^{[a]}_+(Y(A;z))\Delta^{[b]}_-(Y(B;w))\bigr]v^\epsilon_{a,b}\\
&\ \ +Y(A;z)\bigl[-\alpha_+\Delta^{[b]}_-(\Delta^{[a]}_+(Y(B;w)))+\alpha_-\Delta^{[a]}_+(\Delta^{[b]}_-(Y(B;w)))\bigr]v^\epsilon_{a,b}
\end{align*}
Applying Theorem \ref{FJ} to the right-hand side of the above equality, we obtain
\begin{align*}
&J^{\epsilon}_{a^\vee,b^\vee}(A;z)J^{\epsilon}_{a^\vee,b^\vee}(B;z)v^\epsilon_{a,b}\\
&=Y(A;z)Y(B;w)v^\epsilon_{a,b}\\
&\ \ +(\alpha_--\alpha_+)\bigl[\Delta^{[a]}_+(Y(A;z)Y(B;w))+\Delta^{[b]}_-(Y(A;z)Y(B;w))\bigr]v^\epsilon_{a,b}\\
&\ \ -\alpha_+\Delta^{[b]}_-\circ\Delta^{[a]}_+(Y(A;z)Y(B;w))v^\epsilon_{a,b}+\alpha_-\Delta^{[a]}_+\circ\Delta^{[b]}_-(Y(A;z)Y(B;w))v^\epsilon_{a,b}\\
&=Y(Y(A;z-w)B;w)v^\epsilon_{a,b}+\bigl[\bigl(\widetilde{\Delta}_{\tau}+\Delta^{+}_{\tau}+\Delta^{-}_{\tau}\bigr)(Y(Y(A;z-w)B;w))\bigr]v^\epsilon_{a,b}\\
&=J^{\epsilon}_{a^\vee,b^\vee}(Y(A;z-w)B;w)v^\epsilon_{a,b}.
\end{align*}
Thus (\ref{eq:iden-yab}) holds on $v^\epsilon_{a,b}$.
In the same way, we can prove the relations
\begin{align*}
&J^{\epsilon}_{a^\vee,b^\vee}(A;z)J^{\epsilon}_{a^\vee,b^\vee}(B;z)v^{-\epsilon}_{a^\vee,b}=J^{\epsilon}_{a^\vee,b^\vee}(Y(A;z-w)B;w)v^{-\epsilon}_{a^\vee,b},\\
&J^{\epsilon}_{a^\vee,b^\vee}(A;z)J^{\epsilon}_{a^\vee,b^\vee}(B;z)v^{-\epsilon}_{a,b^\vee}=J^{\epsilon}_{a^\vee,b^\vee}(Y(A;z-w)B;w)v^{-\epsilon}_{a,b^\vee},\\
&J^{\epsilon}_{a^\vee,b^\vee}(A;z)J^{\epsilon}_{a^\vee,b^\vee}(B;z)v^\epsilon_{a^\vee,b^\vee}=J^{\epsilon}_{a^\vee,b^\vee}(Y(A;z-w)B;w)v^\epsilon_{a^\vee,b^\vee}.
\end{align*} 
Therefore we obtain (\ref{eq:iden-yab}).
\end{proof}
\end{thm}

As in the argument of \cite[Theorem 6.6]{AK}, by using Proposition \ref{Felder complex} and (\ref{202210130}), we see that the four logarithmic modules $\mathcal{P}^\pm_{\bullet,\bullet}\in C^{thick}_{r,s}$ have $L_0$ nilpotent rank three.

\begin{remark}
We will give the socle series of the logarithmic modules $\mathcal{P}^\pm_{\bullet,\bullet}$ in Subsection \ref{53} (see Proposition \ref{logarithmic2}).
These logarithmic modules $\mathcal{P}^+_{r,s}$, $\mathcal{P}^+_{r^\vee,s^\vee}$, $\mathcal{P}^-_{r^\vee,s}$ and $\mathcal{P}^-_{r,s^\vee}$ correspond to the projective covers of $\X^+_{r,s}$, $\X^+_{r^\vee,s^\vee}$, $\X^-_{r^\vee,s}$ and $\X^-_{r,s^\vee}$, respectively.
\end{remark}

Hereafter, we use the following notation.
\begin{dfn}
\label{not0519}
For any simple $\mathcal{W}_{p_+,p_-}$-module $\mathcal{X}$ and any $n\geq 1$, we define
\begin{align*}
n\mathcal{X}:=\overbrace{\mathcal{X} \oplus \cdots \oplus \mathcal{X}}^{n}.
\end{align*}
\end{dfn}
\begin{prop}
\label{logarithmic11}
By taking quotients of $\mathcal{P}^+_{r,s}$, $\mathcal{P}^+_{r^\vee,s^\vee}$, $\mathcal{P}^-_{r^\vee,s}$ and $\mathcal{P}^-_{r,s^\vee}$, we define eight logarithmic modules $\mathcal{Q}(\X^\epsilon_{a,b})_{b,c}$ where
\begin{align*}
\{(\epsilon,a,b,c,d)\}=&\bigl\{(+,r,s,r^\vee,s),(+,r,s,r,s^\vee),(+,r^\vee,s^\vee,r^\vee,s),(+,r^\vee,s^\vee,r,s^\vee),\\
&(-,r^\vee,s,r,s),(-,r^\vee,s,r^\vee,s^\vee),(-,r,s^\vee,r,s),(-,r,s^\vee,r^\vee,s^\vee)\bigr\},
\end{align*}
and each composition series is given by:
\begin{enumerate}
\item For $\mathcal{Q}(\X^+_{a,b})_{c,d}$, we have the series
$
G_1\subsetneq G_2\subsetneq G_3\subsetneq G_4=\mathcal{Q}(\X^+_{a,b})_{c,d}
$
such that
\begin{align*}
G_1=\X^+_{a,b},\quad G_2/G_1\oplus G_3/G_2=2\X^{-}_{c,d}\oplus L(h_{a,b}),\quad \mathcal{Q}(\X^+_{a,b})_{c,d}/G_3=\X^+_{a,b}.
\end{align*}
\item For $\mathcal{Q}(\X^-_{a,b})_{c,d}$, we have the series
$
G_1\subsetneq G_2\subsetneq G_3\subsetneq G_4=\mathcal{Q}(\X^-_{a,b})_{c,d}
$
such that 
\begin{align*}
G_1=\X^-_{a,b},\quad G_2/G_1\oplus G_3/G_2=2\X^{+}_{c,d},\quad \mathcal{Q}(\X^-_{a,b})_{c,d}/G_3=\X^-_{a,b}.\\
\end{align*}
\end{enumerate}
\end{prop}
\begin{remark}
We will give the socle series of these logarithmic modules in Subsection \ref{structureQ}.
\end{remark}
\subsection{Logarithmic modules in the thin blocks}
For each $(r,s)\in \mathbb{Z}^2_{\geq 1}$ satisfying $1\leq r< p_+$, $1\leq s< p_-$, we set
\begin{align*}
\mathcal{P}_{r,p_-}=\V^+_{r,p_-}\oplus \V^-_{r^\vee,p_-}\in C^{thin}_{r,p_-},\ \ \ \ \mathcal{P}_{p_+,s}=\V^+_{p_+,s}\oplus \V^-_{p_+,s^\vee}\in C^{thin}_{p_+,s}.
\end{align*}
Let $(\mathcal{P}_{r,p_-},Y_{\mathcal{P}_{r,p_-}})$ and $(\mathcal{P}_{p_+,s},Y_{\mathcal{P}_{p_+,s}})$ be the ordinary $\W_{p_+,p_-}$-module. 
As in \cite{NT}, we can construct the following logarithmic modules.
\begin{prop}
\label{action0303}
Let $r$ and $s$ be integers satisfying $1\leq r< p_+$, $1\leq s< p_-$. 
\begin{enumerate}
\item Let $(a,\epsilon)\in \{(r,+),(r^\vee,-)\}$. We define the logarithmic module 
\begin{align*}
(\mathcal{Q}(\X^\epsilon_{a,p_-})_{a^\vee,p_-},J^\epsilon_{a,p_-}) 
\end{align*}
as follows. As the vector space
\begin{align*}
\mathcal{Q}(\X^\epsilon_{a,p_-})_{a^\vee,p_-}=\mathcal{P}_{r,p_-}
\end{align*}
and the module action is defined by
\begin{equation*}
J^\epsilon_{a,p_-}(A;z)
=
\begin{cases}
Y_{\mathcal{P}_{r,p_-}}(A;z)+\Delta^{[a^{\vee}]}_+(Y_{\mathcal{P}_{r,p_-}}(A;z))& on\ \V^{-\epsilon}_{a^\vee,p_-}\\
Y_{\mathcal{P}_{r,p_-}}(A;z) & on\ \V^\epsilon_{a,p_-},
\end{cases}
\end{equation*}
for $A\in \W_{p_+,p_-}$.
\item Let $(b,\epsilon)\in \{(s,+),(s^\vee,-)\}$. We define the logarithmic module 
\begin{align*}
(\mathcal{Q}(\X^\epsilon_{p_+,b})_{p_+,b^\vee},J^\epsilon_{p_+,b})
\end{align*}
as follows. As the vector space 
\begin{align*}
\mathcal{Q}(\X^\epsilon_{p_+,b})_{p_+,b^\vee}=\mathcal{P}_{p_+,s}
\end{align*}
and the module action is defined by
\begin{equation*}
J^\epsilon_{p_+,b}(A;z)
=
\begin{cases}
Y_{\mathcal{P}_{p_+,s}}(A;z)+\Delta^{[b^{\vee}]}_-(Y_{\mathcal{P}_{p_+,s}}(A;z))& on\ \V^{-\epsilon}_{p_+,b^\vee}\\
Y_{\mathcal{P}_{p_+,s}}(A;z) & on\ \V^\epsilon_{p_+,b},
\end{cases}
\end{equation*}
for $A\in \W_{p_+,p_-}$.
\end{enumerate}
\end{prop}

\begin{prop}
Let $r$ and $s$ be integers satisfying $1\leq r< p_+$, $1\leq s< p_-$.
The length of the composition series of the logarithmic modules $\mathcal{Q}(\X^+_{r,p_-})_{r^\vee,p_-}$, $\mathcal{Q}(\X^-_{r^\vee,p_-})_{r,p_-}$, $\mathcal{Q}(\X^+_{p_+,s})_{p_+,s^\vee}$ and $\mathcal{Q}(\X^-_{p_+,s^\vee})_{p_+,s}$ are four, and the composition series of each module is given by:
\begin{enumerate}
\item Let $(a,\epsilon)\in \{(r,+),(r^\vee,-)\}$. The composition series $
G_1\subsetneq G_2\subsetneq G_3\subsetneq G_4
$
of $\mathcal{Q}(\X^\epsilon_{a,p_-})_{a^\vee,p_-}$ is given by
\begin{align*}
\hspace{-6mm}G_1=\X^\epsilon_{a,p_-},\quad G_2/G_1\oplus G_3/G_2=2\X^{-\epsilon}_{a^\vee,p_-},\quad\mathcal{Q}(\X^\epsilon_{a,p_-})_{a^\vee,p_-}/G_3=\X^\epsilon_{a,p_-}.
\end{align*}
\item Let $(b,\epsilon)\in \{(s,+),(s^\vee,-)\}$. The  composition series $
G_1\subsetneq G_2\subsetneq G_3\subsetneq G_4
$ of $\mathcal{Q}(\X^\epsilon_{p_+,b})_{p_+,b^\vee}$ is given by
\begin{align*}
\hspace{-6mm}G_1=\X^\epsilon_{p_+,b},\quad G_2/G_1\oplus G_3/G_2=2\X^{-\epsilon}_{p_+,b^\vee},\quad \mathcal{Q}(\X^\epsilon_{p_+,b})_{p_+,b^\vee}/G_3=\X^\epsilon_{p_+,b}.
\end{align*}
\end{enumerate}
\end{prop}

\begin{remark}
We will give the socle series of these logarithmic modules in Subsection \ref{structureQ}.
These logarithmic modules are projective modules in the thin blocks (see Subsection \ref{54}). The structure of these projective modules are similar to the projective modules of the triplet $W$-algebra $\W_{p}$ determined by \cite{McRae,NT}.
\end{remark}

\section{Logarithmic extension of Virasoro modules}
\label{KRsection}
In this section, by using the results in \cite{KR}, we will determine some ${\rm Ext}^1$ and ${\rm Hom}$ groups between certain logarithmic modules in the abelian category of generalized Virasoro modules.
We note that the logarithmic Virasoro modules considered in this section can also be realized as quotients of the $\mathcal{Q}(\mathcal{X}^\pm_{\bullet,\bullet})_{\bullet,\bullet}$ and $\mathcal{P}^\pm_{\bullet,\bullet}$-type $\mathcal{W}_{p_+,p_-}$-modules defined in the previous section.
From this section, we identify any Virasoro modules that are isomorphic among each other.

\subsection{${\rm Ext}^1$-groups between simple Virasoro modules}
We set
\begin{align*}
A_{p_+,p_-}&:=\{\ \alpha_{r,s;n}\ |\ r,s,n\in\Z\ \},\\
H_{p_+,p_-}&:=\Bigl\{\ h_{\alpha}\ \Bigr|\ \alpha\in A_{p_+,p_-}\ \Bigr\}
\end{align*}
(for the definition of symbols $\alpha_{r,s;n}$ and $h_\alpha$, see (\ref{alpha_{r,s;n}}) and (\ref{h_alpha}), respectively).
\begin{dfn}
We define $\mathcal{L}_{c_{p_+,p_-}}\mathchar`-{\rm mod}$ to be the full subcategory of $\mathcal{L}_{c_{p_+,p_-}}\mathchar`-{\rm Mod}$ such that all objects in $\mathcal{L}_{c_{p_+,p_-}}\mathchar`-{\rm mod}$ satisfy the following conditions:
\begin{enumerate}
\item The socle series of $M$ has finite length.
\item The lowest weights $h$ of the simple modules $L(h)$, appearing in the composition factors of $M$, are elements of $H_{p_+,p_-}$. 
\end{enumerate}
\end{dfn}
For $M_1,M_2\in \mathcal{L}_{c_{p_+,p_-}}\mathchar`-{\rm mod}$, we denote by ${\rm Ext}^1_{\mathcal{L}}(M_1,M_2)$ the ${\rm Ext}^1$-group in $\mathcal{L}_{c_{p_+,p_-}}\mathchar`-{\rm mod}$ consisting of $E\in \mathcal{L}_{c_{p_+,p_-}}\mathchar`-{\rm mod}$ satisfying the exact sequence
\begin{align*}
0\rightarrow M_2\rightarrow E\rightarrow M_1\rightarrow 0.
\end{align*}
For $M_1,M_2\in \mathcal{L}_{c_{p_+,p_-}}\mathchar`-{\rm Mod}$, we denote by ${\rm Hom}_{U(\mathcal{L})}(M_1,M_2)$ the set of the $U(\mathcal{L})\mathchar`-$homomorphisms from $M_1$ to $M_2$. Note that each of ${\rm Ext}^1_{\mathcal{L}}(\bullet,\bullet)$ and ${\rm Hom}_{U(\mathcal{L})}(\bullet,\bullet)$ has the structure of a $\C$-vector space.
%, in what follows.
\vspace{2mm}

Recall that, for $r,s\geq 1$, the Verma module $M(h_{r,s},c_{p_+,p_-})$ has the singular vector $S_{r,s}{\mid}h_{r,s}\rangle$ normalized as (\ref{eq:norm-s}).
For $r,s\geq 1$ and $h\in\mathbb{C}$, let us consider the value
$
\langle h{\mid} \sigma(S_{r,s})S_{r,s}{\mid}h\rangle,
$
where we have chosen a norm of the lowest weight vector ${\mid} h\rangle \in M(h,c_{p_+,p_-})$ as $\langle h{\mid} h\rangle=1$. We can see that this value is a polynomial of $h$ and is divisible by $(h-h_{r,s})$.
A more detailed value is given by the following proposition.
\begin{prop}[\cite{Yanagida}]
For $r,s\geq 1$ and $h\in\mathbb{C}$,
\begin{align*}
\langle h{\mid} \sigma(S_{r,s})S_{r,s}{\mid}h\rangle=R_{r,s}(h-h_{r,s})+O((h-h_{r,s})^2),
%&f(h,c_{p_+,p_-})=R_{r,s}(h-h_{r,s})+O((h-h_{r,s})^2),
\end{align*}
where $R_{r,s}$ is given by
\begin{align*}
R_{r,s}= 2\prod_{\substack{(k,l)\in \mathbb{Z}^2,\\ 1-r\leq k\leq r, 1-s\leq l\leq s,\\(k,l)\neq (0,0),(r,s)}}\Bigl(k\bigl(\frac{p_+}{p_-}\bigr)^{-\frac{1}{2}}+l\bigl(\frac{p_+}{p_-}\bigr)^{\frac{1}{2}}\Bigr).
\end{align*}
\label{Y}
\end{prop}

\begin{remark}
In this paper, it is important that $R_{r,s}$ be non-zero, specific value is not necessary. In fact, the non-triviality of $R_{r,s}$ can be shown using the Jantzen-filtration of the Fock module $F_{r,s}$.
\end{remark}
By using Proposition \ref{Y}, we obtain the following proposition (cf.\cite{GK}).
\begin{prop}
\label{Vir0}
For $h\in H_{p_+,p_-}$, we have
\begin{align*}
{\rm Ext}^1_{\mathcal{L}}(L(h),L(h))=0.
\end{align*}
\begin{proof}
We prove only in the case of $h=h_{r,s} (1\leq r<p_+,1\leq s<p_- )$. The other cases can be proved in the same way.

Assume ${\rm Ext}^1_{\mathcal{L}}(L(h_{r,s}),L(h_{r,s}))\neq 0$. Fix a non-trivial extension 
\begin{align*}
0\rightarrow L(h_{r,s})\xrightarrow{\iota} E\xrightarrow{\pi}L(h_{r,s})\rightarrow 0.
\end{align*}
Let $\{u_0,u_1\}$ be a basis of the lowest weight space of $E$ such that
\begin{align*}
\pi(u_0)=\ket{h_{r,s}},\qquad \iota(\ket{h_{r,s}})=u_1,\qquad 
(L_0-h_{r,s})u_0=cu_1,
\end{align*}
where $c$ is a non-zero constant and $\ket{h_{r,s}}$ is the lowest weight vector of $L(h_{r,s})$. Then, by Proposition \ref{Y}, we have
\begin{align*}
&\sigma(S_{r,s})S_{r,s}u_0=f(c)u_1,
&\frac{f(c)}{c}\Bigr|_{c=0}\neq 0,
\end{align*}
where $f(c)$ is a polynomial of $c$. 
Thus, we see that $S_{r,s}u_0$ is non-zero and 
\begin{align*}
S_{r,s}u_0\in \iota(L(h_{r,s})).
\end{align*} 
On the other hand, by the simplicity of $L(h_{r,s})$, we have
\begin{align*}
\sigma(S_{r,s})S_{r,s}u_0=0.
\end{align*}
But this is a contradiction.
\end{proof}
\end{prop}

The following proposition is due to \cite{BNW,FF,IK}.
\begin{prop}[\cite{BNW,FF,IK}]
For any $r,s,\in \mathbb{Z}$ and $h\in \C$, the ${\rm Ext}^1$-group ${\rm Ext}^1_{\mathcal{L}}(L(h_{r,s;n}),L(h))$ is given by:
\begin{enumerate}
\item For $1\leq r<p_+,1\leq s<p_-$ and $n=0$, we have
\begin{equation*}
{\rm Ext}^1_{\mathcal{L}}(L(h_{r,s}),L(h))=
\begin{cases}
\C& for\ h=h_{r^\vee,s;-1}\ or\ h_{r^\vee,s;1}\\
0& otherwise 
\end{cases}
.
\end{equation*}
\item For $1\leq r<p_+,1\leq s<p_-$ and $n\geq 1$, we have
\begin{equation*}
{\rm Ext}^1_{\mathcal{L}}(L(h_{r,s;n}),L(h))=
\begin{cases}
\C& for\ h=h_{r^\vee,s;n\pm1}\ or\ h_{r,s^\vee;n\pm1}\\
0& otherwise 
\end{cases}
.
\end{equation*}
\item For $1\leq r<p_+,s=p_-$ and $n=0$, we have
\begin{equation*}
{\rm Ext}^1_{\mathcal{L}}(L(h_{r,p_-}),L(h))=
\begin{cases}
\C& for\ h=h_{r^\vee,p_-;1}\\
0& otherwise 
\end{cases}
.
\end{equation*}
\item For $1\leq r<p_+,s=p_-$ and $n\geq 1$, we have
\begin{equation*}
{\rm Ext}^1_{\mathcal{L}}(L(h_{r,p_-;n}),L(h))=
\begin{cases}
\C& for\ h=h_{r^\vee,p_-;n+1}\ or\ h_{r^\vee,p_-;n-1}\\
0& otherwise 
\end{cases}
.
\end{equation*}
\item For $r=p_+,1\leq s<p_-$ and $n=0$, we have
\begin{equation*}
{\rm Ext}^1_{\mathcal{L}}(L(h_{p_+,s}),L(h))=
\begin{cases}
\C& for\ h=h_{p_+,s^\vee;-1}\\
0& otherwise 
\end{cases}
.
\end{equation*}
\item For $r=p_+,1\leq s<p_-$ and $n\leq -1$, we have
\begin{equation*}
{\rm Ext}^1_{\mathcal{L}}(L(h_{p_+,s;n}),L(h))=
\begin{cases}
\C& for\ h=h_{p_+,s^\vee;n-1}\ or\ h_{p_+,s^\vee;n+1}\\
0& otherwise 
\end{cases}
.
\end{equation*}
\item For $r=p_+,s=p_-,n\in\Z$, we have
\begin{align*}
{\rm Ext}^1_{\mathcal{L}}(L(h_{p_+,p_-;n}),L(h))=0.
\end{align*}
\end{enumerate}
\label{VirC}
\end{prop}

\subsection{Logarithmic extensions I}

In the following, we identify any Virasoro modules that are isomorphic among each other.
Let us define the following indecomposable modules in $\mathcal{L}_{c_{p_+,p_-}}\mathchar`-{\rm mod}$ as quotient modules of  certain Virasoro Verma modules. 

\begin{dfn}
We define $\mathcal{T}_{p_+,p_-}$ to be the subset of $A^3_{p_+,p_-}$ such that every element $(\alpha_1,\alpha_2,\alpha_3)\in A^3_{p_+,p_-}$ satisfies the following conditions:
\begin{enumerate}
\item $h_{\alpha_1}\leq h_{\alpha_2}< h_{\alpha_3}$.
\item The three Fock modules $F_{\alpha_1}$, $F_{\alpha_2}$ and $F_{\alpha_3}$ are contained in the same Felder complex in Proposition \ref{Felder complex2} and are adjacent to each other as
\begin{align*}
\cdots\xrightarrow{}F_{\alpha_1}\xrightarrow{Q^{[\bullet]}_\epsilon}F_{\alpha_2}\xrightarrow{Q^{[\bullet]}_\epsilon}F_{\alpha_3}\xrightarrow{}\cdots.
\end{align*}
\end{enumerate}
\end{dfn}

\begin{dfn}
We define the following subsets of $\mathcal{T}_{p_+,p_-}$:
\begin{align*}
\mathcal{T}^{\rm ch}_{p_+,p_-}&:=\{(\alpha_1,\alpha_2,\alpha_3)\in \mathcal{T}_{p_+,p_-}\ |\ F_{\alpha_1}\ {\rm is}\ {\rm chain}\ {\rm type}\ \},\\
\mathcal{T}^{\rm br}_{p_+,p_-}&:=\{(\alpha_1,\alpha_2,\alpha_3)\in \mathcal{T}_{p_+,p_-}\ |\ F_{\alpha_1}\ {\rm is}\ {\rm braided}\ {\rm type}\ \},\\
\mathcal{T}^0_{p_+,p_-}&:=\{(\alpha_1,\alpha_2,\alpha_3)\in \mathcal{T}^{\rm ch}_{p_+,p_-}\ |\ h_{\alpha_1}=h_{\alpha_2}\},\\
\mathcal{T}^{{\rm Min}}_{p_+,p_-}&:=\{\ (\alpha_{r,s},\alpha_{2},\alpha_{3})\in \mathcal{T}^{\rm br}_{p_+,p_-}\ |\ \alpha_1=\alpha_{r,s},\ 1\leq r<p_+,\ 1\leq s<p_-\}.
\end{align*}
\end{dfn}

\begin{remark}
{}\mbox{}
\begin{itemize}
\item Note that $\# \mathcal{T}^0_{p_+,p_-}=p_++p_--2$ and every element of $\mathcal{T}^0_{p_+,p_-}$ is given by
$
(\alpha_{r,p_-;-1},\alpha_{r,p_-;0},\alpha_{r,p_-;1}),
(\alpha_{p_+,s^\vee;1},\alpha_{p_+,s;0},\alpha_{p_+,s^\vee;-1})
$
for $1\leq r<p_+$, $1\leq s<p_-$.
\item Note that $\# \mathcal{T}^{\rm Min}_{p_+,p_-}=(p_+-1)(p_--1)$ and every element of $\mathcal{T}^{\rm Min}_{p_+,p_-}$ is given by
$
(\alpha_{r,s},\alpha_{r^\vee,s;1},\alpha_{r,s;2}),
(\alpha_{r,s},\alpha_{r,s^\vee;-1},\alpha_{r,s;-2})
$
for $1\leq r<p_+$, $1\leq s<p_-$.
\end{itemize}
\end{remark}
\begin{dfn}
For $h,h'\in H_{p_+,p_-}$ such that ${\rm Ext}^1_{\mathcal{L}}(L(h),L(h'))\simeq \C$ and $h<h'$, 
let $L(h,h')$ be an indecomposable module representing the nonzero class in ${\rm Ext}^1_{\mathcal{L}}(L(h),L(h'))$.
\end{dfn}

\begin{dfn}
\begin{enumerate}
\item For any $\tau=(\alpha_1,\alpha_2,\alpha_3)\in\mathcal{T}^0_{p_+,p_-}$, we define $K(\tau)=L(h_{\alpha_2},h_{\alpha_3})$.
\item Given $\tau=(\alpha_1,\alpha_2,\alpha_3)\in\mathcal{T}_{p_+,p_-}$ such that $h_{\alpha_1}\neq h_{\alpha_2}$, from Proposition \ref{VirC}, we have
\begin{align*}
{\rm Ext}^1_{\mathcal{L}}(L(h_{\alpha_2},h_{\alpha_3}), L(h_{\alpha_1}))\simeq \C.
\end{align*} 
Let $K(\tau)$ be an indecomposable module representing the nonzero class in the above ${\rm Ext}^1$-group.
\end{enumerate}
\end{dfn}
Let $\tau=(\alpha_1,\alpha_2,\alpha_3)\in \mathcal{T}_{p_+,p_-}\setminus\mathcal{T}^0_{p_+,p_-}$.
From the Virasoro subquotient structure of the logarithmic $\W_{p_+,p_-}$-modules $\mathcal{Q}(\X^\pm_{r,s})_{\bullet,\bullet}$ defined in Section \ref{logarithmic section}, we see that there exists a logarithmic module corresponding to an extension class in
\begin{align}
\label{eq:kr-ext2026}
{\rm Ext}^1_{\mathcal{L}}(L(h_{\alpha_2},h_{\alpha_3}),L(h_{\alpha_1},h_{\alpha_2})).
\end{align}
Fix an arbitrary logarithmic module $E$ corresponding to an extension class in the ${\rm Ext}^1$-group (\ref{eq:kr-ext2026}). Let $v_{\alpha_1}$ be a nonzero vector in the one dimensional subspace $E[h_{\alpha_1}]$. 
Let $v_{\alpha_2}$ be a vector in $E[h_{\alpha_2}]$ such that $(L_0-h_{\alpha_2})v_{\alpha_2}\neq 0$.
Note that the lowest weight vector of the submodule $L(h_{\alpha_2})\subset L(h_{\alpha_1},h_{\alpha_2})$ can be written in terms of the lowest weight vector of $L(h_{\alpha_1},h_{\alpha_2})$ and a Shapovalov element of $U(\mathcal{L})$. We denote by $\eta(h_{\alpha_1},h_{\alpha_2})$ such a Shapovalov element.
Then, there exist constants $c_E(v_{\alpha_1},v_{\alpha_2})$, $d_E(v_{\alpha_1},v_{\alpha_2})$ such that 
\begin{align*}
\sigma(\eta(h_{\alpha_1},h_{\alpha_2}))v_{\alpha_2}&=c_E(v_{\alpha_1},v_{\alpha_2})v_{\alpha_1},\\
(L_0-h_{\alpha_2})v_{\alpha_2}&=d_E(v_{\alpha_1},v_{\alpha_2})\eta(h_{\alpha_1},h_{\alpha_2})v_{\alpha_1}.
\end{align*}
We see that the quotient $c_E(v_{\alpha_1},v_{\alpha_2})/d_E(v_{\alpha_1},v_{\alpha_2})$ is independent of the choice of $v_{\alpha_1}$ and $v_{\alpha_2}$. Then we set
\begin{align*}
\beta_E=\frac{c_E(v_{\alpha_1},v_{\alpha_2})}{d_E(v_{\alpha_1},v_{\alpha_2})}.
\end{align*}
This value is called the logarithmic coupling of $E$.
The following theorem is due to \cite[Theorem~6.15]{KR} (see also \cite{CR}).
\begin{thm}[\cite{KR}]
\label{KR0}
Let $\tau=(\alpha_1,\alpha_2,\alpha_3)\in \mathcal{T}_{p_+,p_-}\setminus\mathcal{T}^0_{p_+,p_-}$. Let $E$ be an arbitrary logarithmic module  corresponding to an extension class in the ${\rm Ext}^1$-group (\ref{eq:kr-ext2026}). Then we have $\beta_E\neq 0$.
In particular, the quotient $E/L(h_{\alpha_2})$ is indecomposable.
\end{thm}
Fix an arbitrary $\tau=(\alpha_1,\alpha_2,\alpha_3)\in \mathcal{T}_{p_+,p_-}\setminus\mathcal{T}^0_{p_+,p_-}$. 
By Theorem~\ref{KR0}, we have a logarithmic Virasoro module representing classes in both (\ref{eq:kr-ext2026}) and 
\begin{align}
\label{eq:kr-ext2026-2}
{\rm Ext}^1_{\mathcal{L}}(K(\tau),L(h_{\alpha_2})).
\end{align}
On the other-hand, from the Virasoro subquotient structure of the logarithmic $\W_{p_+,p_-}$-modules $\mathcal{Q}(\X^\pm_{r,s})_{\bullet,\bullet}$ defined in Section \ref{logarithmic section} (see also logarithmic modules defined in Definition \ref{seealso}), we have a logarithmic Virasoro module corresponding to a nonzero extension class in (\ref{eq:kr-ext2026-2}).
The following proposition shows that they are isomorphic.
\begin{prop}
\label{sankakuthick}
For any $\tau=(\alpha_1,\alpha_2,\alpha_3)\in \mathcal{T}_{p_+,p_-}$, we have
\begin{align}
\label{ExtC0516}
{\rm Ext}^1_{\mathcal{L}}(K(\tau),L(h_{\alpha_2}))\simeq \C.
\end{align}
\begin{proof}
Fix a logarithmic extension $[E_1]$ in the ${\rm Ext}^1$-group (\ref{eq:kr-ext2026}).
% and denote it by $P(\tau)$. 
If $\tau\in \mathcal{T}^0_{p_+,p_-}$, then we obtain the assertion of this theorem by Propositions \ref{Vir0} and \ref{VirC}, and thus
let $\tau\notin \mathcal{T}^0_{p_+,p_-}$.
It is sufficient to show that ${\rm Ext}^1_{\mathcal{L}}(E_1,L(h_{\alpha_2}))=0$. Assume that 
\begin{align}
\label{assumpt20221203}
{\rm Ext}^1_{\mathcal{L}}(E_1,L(h_{\alpha_2}))\neq 0.
\end{align}
Note that, by Theorem \ref{KR0}, $E_1$ has $L(h_{\alpha_1},h_{\alpha_2})$ as a submodule. Then, by the exact sequence
$
0\rightarrow L(h_{\alpha_1},h_{\alpha_2})\rightarrow E_1\rightarrow L(h_{\alpha_2},h_{\alpha_3})\rightarrow 0
$
and by the assumption (\ref{assumpt20221203}), we have
\begin{align*}
{\rm Ext}^1_{\mathcal{L}}(L(h_{\alpha_2},h_{\alpha_3}),L(h_{\alpha_2}))\neq 0.
\end{align*}
Let $[E_2]$ be any non-trivial extension in ${\rm Ext}^1_{\mathcal{L}}(L(h_{\alpha_2},h_{\alpha_3}),L(h_{\alpha_2}))$. Then, by Propositions~\ref{Y} and \ref{Vir0}, $E_2$ is a logarithmic Virasoro module and must have $L(h_{\alpha_2},h_{\alpha_3})^*$ as a submodule. By the exact sequence
$
0\to L(h_{\alpha_2},h_{\alpha_3})^*\to E_2\rightarrow L(h_{\alpha_2})\to 0,
$
we have the exact sequence
$
0\rightarrow \C\rightarrow {\rm Ext}^1_{\mathcal{L}}(E_2,L(h_{\alpha_1})).
$
Thus we have ${\rm Ext}^1_{\mathcal{L}}(E_2,L(h_{\alpha_1}))\neq 0$.
Let $[E_3]\in {\rm Ext}^1_{\mathcal{L}}(L(h_{\alpha_1}),E^*_2)\setminus\{0\}$. 
Note that since $E_2$ is logarithmic, $E_3$ is also a logarithmic Virasoro module.
By definition, 
\begin{align*}
E_3/L(h_{\alpha_2})=L(h_{\alpha_1})\oplus L(h_{\alpha_2},h_{\alpha_3}).
\end{align*}
But this contradicts Theorem \ref{KR0}.
\end{proof}
\end{prop}
We define the following indecomposable module  (see also the schematic diagram in Figure~\ref{51figk}).
\begin{dfn}
\label{sankakuthickdf}
Given $\tau=(\alpha_1,\alpha_2,\alpha_3)\in\mathcal{T}_{p_+,p_-}$, 
let $P(\tau)$ be an indecomposable module representing the nonzero class in the ${\rm Ext}^1$-group (\ref{ExtC0516}).
\end{dfn}

\begin{remark}
By Propositions \ref{Vir0} and \ref{sankakuthick}, we see that for any $\tau\in \mathcal{T}^0_{p_+,p_-}$, the logarithmic Virasoro mocule $P(\tau)$ is self-contragredient. 
\label{rem:log-tau0}
\end{remark}
\begin{figure}[htbp]
  \centering
\begin{tikzpicture}[scale=0.85]
\node[inner sep=0.7pt] (a1) at (0, 0) {${\ket{h_{\alpha_1}}}$};
\node[inner sep=0.7pt] (a2) at (3, 0) {${\ket{h_{\alpha_2}}}$};
\node[inner sep=0.7pt] (b1) at (1.5, -3) {${\ket{h_{\alpha_3}}}$};
\draw[arrows = {-Stealth[scale=0.9]}] (a1) to (b1);
\draw[arrows = {-Stealth[scale=0.9]},dashed] (a1) to (a2);
\draw[arrows = {-Stealth[scale=0.9]}] (b1) to (a2);

\node[inner sep=0.7pt,shift={(7,0)}] (p1) at (0, -1.5) {${\ket{h_{\alpha_2}}}$};
\node[inner sep=0.7pt,shift={(7,0)}] (p2) at (1.5, 0) {${\ket{h_{\alpha_1}}}$};
\node[inner sep=0.7pt,shift={(7,0)}] (q1) at (1.5, -3) {${\ket{h_{\alpha_3}}}$};
\node[inner sep=0.7pt,shift={(7,0)}] (t) at (3, -1.5) {${\ket{h_{\alpha_2}}}$};
\draw[arrows = {-Stealth[scale=0.9]},shift={(7,0)}] (p1) to (q1);
\draw[arrows = {-Stealth[scale=0.9]},shift={(7,0)}] (p1) to (q1);
\draw[arrows = {-Stealth[scale=0.9]},shift={(7,0)}] (p1) to (p2);
\draw[arrows = {-Stealth[scale=0.9]},dotted,shift={(7,0)}] (q1) to (t);
\draw[arrows = {-Stealth[scale=0.9]},shift={(7,0)}] (p2) to (t);
\draw[arrows = {-Stealth[scale=0.9]},dashed,shift={(7,0)}] (p1) to (t);
\end{tikzpicture}
\caption{The schematic diagrams of the logarithmic modules ${P}(\tau)$ defined in Definition \ref{sankakuthickdf}. The left diagram corresponds to the case $\tau=(\alpha_1,\alpha_2,\alpha_3)\in \mathcal{T}^0_{p_+,p_-}$, and the right to the case $\tau\in \mathcal{T}_{p_+,p_-}\setminus\mathcal{T}^0_{p_+,p_-}$. Each ket vector represents the lowest weight vector of the simple module in some quotient. Each arrow between the ket vectors represents that in some quotient, multiplying by some element of $U(\mathcal{L})$ produces the ket vector at the arrow end. The dashed arrows represent the action $L_0-h_{\alpha_2}$. The dotted arrow in the right figure is correct if ${\rm S}({P}(\tau))=L(h_{\alpha_2})$, but this cannot be determined from the above arguments. \label{51figk}}
\mbox{}\\
\end{figure}

The following indecomposable modules $\widetilde{K}(h_{r,s;n})$ can be realized as the quotients of Fock modules.
\begin{dfn}
\label{dfnC2}
\mbox{}
Let $r$ and $s$ be integers satisfying $1\leq r< p_+$, $1\leq s< p_-$.
\begin{enumerate}
\item By Proposition \ref{VirC}, the indecomposable module module representing a nonzero class in
\begin{align*}
%{\rm Ext}^1_{\mathcal{L}}(L(h_{r,s;1}),L(h_{r,s;0})\oplus L(h_{r^\vee,s;2})\oplus L(h_{r,s^\vee;2}))\setminus\{0\}
{\rm Ext}^1_{\mathcal{L}}(L(h_{r,s;1}),L(h_{r^\vee,s;0})\oplus L(h_{r^\vee,s;2})\oplus L(h_{r,s^\vee;2}))\setminus\{0\}
\end{align*}
is unique up to isomorphism. 
Let $\widetilde{K}(h_{r,s;1})$ be such a module.
\item For $n\geq 2$, by Proposition \ref{VirC}, the indecomposable module module representing a nonzero class in
\begin{align*}
{\rm Ext}^1_{\mathcal{L}}(L(h_{r,s;n}),\bigoplus_{i=-1,1}L(h_{r,s^\vee;n+i})\oplus \bigoplus_{j=-1,1} L(h_{r^\vee,s;n+j}))\setminus\{0\}
\end{align*}
is unique up to isomorphism. 
Let $\widetilde{K}(h_{r,s;n})$ be such a module.
\end{enumerate}
\end{dfn}

\begin{dfn}
Given $1\leq r<p_+$, $1\leq s<p_-$, $n\geq 0$, by Proposition \ref{VirC}, we have 
\begin{align*}
{\rm Ext}^1(L(h_{r,s;n},h_{r^\vee,s;n+1}),L(h_{r,s^\vee;n+1}))\simeq \C.
\end{align*}
Let $\widetilde{L}(h_{r,s;n})$ be an indecomposable module representing the nonzero class in the above ${\rm Ext}^1$-group.
\end{dfn}
Let $1\leq r<p_+$, $1\leq s<p_-$, $n\geq 2$.
From the Virasoro subquotient structure of the logarithmic $\W_{p_+,p_-}$-modules $\mathcal{P}^\pm_{r,s}$ defined in Section \ref{logarithmic section}, we see that there exists a logarithmic module corresponding to an extension class in
\begin{align} 
\label{eq:d-beta}
{\rm Ext}^1_{\mathcal{L}}(\widetilde{L}(h_{r,s;n}),L(h_{a,b;n-1},h_{r,s;n})),
\end{align} 
where $(a,b)=(r^\vee,s)$ or $(r,s^\vee)$.
The following theorem is due to \cite[Theorem~{6.15}]{KR}.
\begin{thm}[\cite{KR}]
\label{KR02}
Let $r$, $s$ and $n$ be integers satisfying $1\leq r<p_+$, $1\leq s<p_-$ and $n\geq 2$. Let $(a,b)=(r^\vee,s)$ or $(r,s^\vee)$.
Then, for any logarithmic module $E$ representing a class in the ${\rm Ext}^1$-group (\ref{eq:d-beta}), the quotient module $E/L(h_{r,s;n})$ is indecomposable.
\end{thm}
The following proposition is a consequence of Theorem \ref{KR02} (see also the schematic diagram in Figure~\ref{51figkk}).
\begin{prop}
\label{sankakuthick2}
For $(r,s,n)\in \mathbb{Z}^3$ satisfying $1\leq r< p_+$, $1\leq s< p_-$ and $n\geq 2$, we have
\begin{align}
\label{20520}
{\rm Ext}^1_{\mathcal{L}}(\widetilde{K}(h_{r,s;n}),L(h_{r,s;n}))\simeq\C^2
\end{align}
and there exists an indecomposable module representing a class in 
\begin{align*}
{\rm Ext}^1_{\mathcal{L}}(\widetilde{K}(h_{r,s;n}),L(h_{r,s;n})\oplus L(h_{r,s;n}))
\end{align*}
such that has quotients isomorphic to $P(\tau_1)$ and $P(\tau_2)$, where 
\begin{align*}
&\tau_1=(\alpha_{r^\vee,s;n-1},\alpha_{r,s;n},\alpha_{r^\vee,s;n+1}),\\
&\tau_2=(\alpha_{r^\vee,s;-n+1},\alpha_{r^\vee,s^\vee;-n},\alpha_{r^\vee,s;-n-1}).
\end{align*}

\begin{proof}
Noting Proposition \ref{sankakuthick}, from the Virasoro subquotient structure of the $\W_{p_+,p_-}$-modules $\mathcal{P}^{\pm}_{\bullet,\bullet}$, we obtain an indecomposable module corresponding to an extension class in
\begin{align*}
{\rm Ext}^1_{\mathcal{L}}(\widetilde{K}(h_{r,s;n}),L(h_{r,s;n})\oplus L(h_{r,s;n})).
\end{align*}
We denote by ${\widetilde{P}}(h_{r,s;n})$ this indecomposable module. 
Note that, by Proposition \ref{sankakuthick}, ${\widetilde{P}}(h_{r,s;n})$ has $P(\tau_1)$ and $P(\tau_2)$ as quotients, where 
\begin{align*}
&\tau_1=(\alpha_{r^\vee,s;n-1},\alpha_{r,s;n},\alpha_{r^\vee,s;n+1}),\\
&\tau_2=(\alpha_{r^\vee,s;-n+1},\alpha_{r^\vee,s^\vee;-n},\alpha_{r^\vee,s;-n-1}).
\end{align*}
As in the proof of Proposition \ref{sankakuthick}, we can show 
\begin{align*}
{\rm Ext}^1_{\mathcal{L}}({\widetilde{P}}(h_{r,s;n}),L(h_{r,s;n}))=0
\end{align*}
by using Theorems~\ref{KR0} and~\ref{KR02}. Thus, by Proposition \ref{Vir0}, we obtain (\ref{20520}).
\end{proof}
\end{prop}

\mbox{}\\
\begin{figure}[htbp]
  \centering
\begin{tikzpicture}[scale=0.85]
\node[inner sep=0.7pt] (p1) at (0, -1.5) {${\ket{h_{r,s;n}}}$};
\node[inner sep=0.7pt] (d2) at (-1.5, 0) {${\ket{h_{r^\vee,s;n-1}}}$};
\node[inner sep=0.7pt] (f1) at (-1.5, -3) {${\ket{h_{r^\vee,s;n+1}}}$};
\node[inner sep=0.7pt] (s) at (-3, -1.5) {${\ket{h_{r,s;n}}}$};

\draw[arrows = {-Stealth[scale=0.9]}] (p1) to (f1);
\draw[arrows = {-Stealth[scale=0.9]}] (p1) to (d2);
\draw[arrows = {-Stealth[scale=0.9]},dotted] (f1) to (s);
\draw[arrows = {-Stealth[scale=0.9]}] (d2) to (s);
\draw[arrows = {-Stealth[scale=0.9]},dashed] (p1) to (s);

\node[inner sep=0.7pt] (p2) at (1.5, 0) {${\ket{h_{r^\vee,s;-n+1}}}$};
\node[inner sep=0.7pt] (q1) at (1.5, -3) {${\ket{h_{r^\vee,s;-n-1}}}$};
\node[inner sep=0.7pt] (t) at (3, -1.5) {${\ket{h_{r,s;n}}}$};

\draw[arrows = {-Stealth[scale=0.9]}] (p1) to (q1);
\draw[arrows = {-Stealth[scale=0.9]}] (p1) to (q1);
\draw[arrows = {-Stealth[scale=0.9]}] (p1) to (p2);
\draw[arrows = {-Stealth[scale=0.9]},dotted] (q1) to (t);
\draw[arrows = {-Stealth[scale=0.9]}] (p2) to (t);
\draw[arrows = {-Stealth[scale=0.9]},dashed] (p1) to (t);

\end{tikzpicture}
\caption{The schematic diagrams of the logarithmic modules ${\widetilde{P}}(h_{r,s;n})$ defined in the proof of Proposition \ref{sankakuthick2}. Each ket vector $\ket{h_{\bullet,\bullet;\bullet}}$ represents the lowest weight vector of the simple module $L(h_{\bullet,\bullet;\bullet})$ in some quotient. Each arrow between the ket vectors represents that in some quotient, multiplying by some element of $U(\mathcal{L})$ yields the ket vector at the arrow end. The dashed arrows represent the action $L_0-h_{r,s;n}$. The dotted arrows are correct if ${\rm S}({\widetilde{P}}(h_{r,s;n}))=L(h_{r,s;n})\oplus L(h_{r,s;n})$, but this cannot be determined from the above arguments. \label{51figkk}}
\end{figure}

\subsection{Logarithmic extensions II}
So far we have studied properties for ${\rm Ext}^1$-groups. In the following, we will investigate the subquotient structure of logarithmic extensions in the ${\rm Ext}^1$-groups
\begin{align*}
{\rm Ext}^1_{\mathcal{L}}(K(\tau),L(h_{\alpha_2})),\ \ \tau=(\alpha_1,\alpha_2,\alpha_3)\in \mathcal{T}^{\rm Min}_{p_+,p_-}
\end{align*}
and
\begin{align*}
{\rm Ext}^1_{\mathcal{L}}(K(\tau),L(h_{\alpha_2},h_{\alpha_3})^*),\ \ \tau=(\alpha_1,\alpha_2,\alpha_3)\in \mathcal{T}^0_{p_+,p_-}\sqcup\mathcal{T}^{\rm Min}_{p_+,p_-}.
\end{align*}
For that purpose, we will construct a family of infinite length logarithmic Virasoro modules by using certain limit operations. A similar construction is given by \cite{C,Nakano}.

For any $\alpha,\beta\in\mathbb{C}$, we identify
$
e^{\beta\hat{a}}{\mid}\alpha\rangle
$
with 
$
{\mid}\alpha+\beta\rangle
$
(see Remark~\ref{rem:beta-gamma}).
For $\alpha\in A_{p_+,p_-}$, let $F_{\alpha}$ be a braided or chain type Fock module and let $Q^{[k]}_{\epsilon}$ be a screening operator acting on $F_{\alpha}$, where $\epsilon=+$ or $-$. Let $v$ be any $L_0$-homogeneous vector of $F_{\alpha}$ and let $A\in U(\mathcal{L})$ be any $L_0$-homogeneous element. 
Let $n_1$ and $n_2$ be the integers such that $L_0v=(n_1+h_{\alpha})v$ and $[L_0,A]=-n_2A$.
For any $x\in \mathbb{C}^\times$, consider the commutator
\begin{align*}
[Q^{[k]}_\epsilon,e^{-x\hat{a}}Ae^{x\hat{a}}]= Q^{[k]}_\epsilon e^{-x\hat{a}}Ae^{x\hat{a}}-e^{-x\hat{a}}Ae^{x\hat{a}}Q^{[k]}_\epsilon
\end{align*}
on $v$, where $F_{\alpha+x},F_{\alpha+k\alpha_{\epsilon}+x}\in \mathcal{F}_{\alpha_0}\mathchar`-{\rm mod}$.
Write $[Q^{[k]}_\epsilon,e^{-x\hat{a}}Ae^{x\hat{a}}]v$ as 
\begin{align*}
[Q^{[k]}_\epsilon,e^{-x\hat{a}}Ae^{x\hat{a}}]v=\sum_{\lambda\vdash N}f_{\lambda}(x)a_{-\lambda}{\mid}\alpha+k\alpha_\epsilon\rangle,
\end{align*}
where $N=n_1+n_2+h_{\alpha}-h_{\alpha+k\alpha_{\epsilon}}$, and  $f_{\lambda}(x)$ are some polynomials of $x$. Since $[Q^{[k]}_\epsilon,A]=0$, we see that every $f_{\lambda}(x)$ is divisible by $x$. Then we define
\begin{align*}
\lim_{\substack{x\rightarrow 0\\ x\neq 0}}\frac{1}{x}[Q^{[k]}_\epsilon,e^{-x\hat{a}}Ae^{x\hat{a}}]v:=\sum_{\lambda\vdash N}(x^{-1}f_{\lambda}(x){\mid}_{x=0})a_{-\lambda}{\mid}\alpha+k\alpha_\epsilon\rangle.
%=\sum _{i=1}^Nt^i v_i,
%F_{r^\vee,k;n+1}.
\end{align*}
We extend this operation linearly to all elements of $U(\mathcal{L})$.
Then, we introduce the following $\mathbb{C}$-linear operators.
\begin{dfn}
Let $F_{\alpha}(\alpha\in A_{p_+,p_-})$ be any braided or chain type Fock module and let $Q^{[k]}_{\epsilon}$ be a screening operator acting on $F_{\alpha}$, where $\epsilon=+$ or $-$. We define the $\C$-linear operator
$
\Lambda_{Q^{[k]}_\epsilon}:\ U(\mathcal{L})\rightarrow {\rm Hom}_\C(F_{\alpha},F_{\alpha+k\alpha_{\epsilon}})
$
as follows
\begin{align*}
\Lambda_{Q^{[k]}_\epsilon}(A)=\lim_{\substack{x\rightarrow 0\\ x\neq 0}}\frac{1}{x}[Q^{[k]}_\epsilon,e^{-x\hat{a}}Ae^{x\hat{a}}],\ \ \ \ {\rm for}\ A\in U(\mathcal{L}),
\end{align*}
where $F_{\alpha+x},F_{\alpha+k\alpha_{\epsilon}+x}\in \mathcal{F}_{\alpha_0}\mathchar`-{\rm mod}$.
\end{dfn}

\begin{prop}
\label{deri}
The two operators $\Lambda_{Q^{[k]}_\epsilon}$ satisfiy the following derivation property
\begin{align*}
\Lambda_{Q^{[k]}_\epsilon}(AB)=\Lambda_{Q^{[k]}_\epsilon}(A)B+A\Lambda_{Q^{[k]}_\epsilon}(B),\ \ \ A,B\in U(\mathcal{L}).
\end{align*}
\begin{proof}
For any $A,B\in U(\mathcal{L})$, we have
\begin{align*}
&[Q^{[k]}_\epsilon,e^{-x\hat{a}}ABe^{x\hat{a}}]\\
&=[Q^{[k]}_\epsilon,e^{-x\hat{a}}Ae^{x\hat{a}}\cdot e^{-x\hat{a}}Be^{x\hat{a}}]\\
&=[Q^{[k]}_\epsilon,e^{-x\hat{a}}Ae^{x\hat{a}}]e^{-x\hat{a}}Be^{x\hat{a}}+e^{-x\hat{a}}Ae^{x\hat{a}}[Q^{[k]}_\epsilon,e^{-x\hat{a}}Be^{x\hat{a}}]\\
&=[Q^{[k]}_\epsilon,e^{-x\hat{a}}Ae^{x\hat{a}}]B+A[Q^{[k]}_\epsilon,e^{-x\hat{a}}Be^{x\hat{a}}]\\
&\ \ \ \ +[Q^{[k]}_\epsilon,e^{-x\hat{a}}Ae^{x\hat{a}}](e^{-x\hat{a}}Be^{x\hat{a}}-B)+(e^{-x\hat{a}}Ae^{x\hat{a}}-A)[Q^{[k]}_\epsilon,e^{-x\hat{a}}Be^{x\hat{a}}].
\end{align*}
Dividing both sides by $x$ and taking the limit, we have the derivation property.
\end{proof}
\end{prop}

We define the following symbols for some screening operators.
\begin{dfn}
\begin{enumerate}
\item For $\tau=(\alpha_1,\alpha_2,\alpha_3)\in \mathcal{T}_{p_+,p_-}$, we denote by $Q_{\tau}$ the screening operator from $F_{\alpha_1}$ to $F_{\alpha_2}$, and we set $\Lambda_{\tau}=\Lambda_{Q_{\tau}}$. 

\item Let $\tau=(\alpha_1,\alpha_2,\alpha_3)\in \mathcal{T}^{\rm br}_{p_+,p_-}$. Let $(r_\tau,s_\tau)\in\Z^2$ be a unique integer pair satisfying
\begin{align*}
(1\leq r_\tau<p_+\land 1\leq s_\tau<p_-) \land (Q^{[r_\tau]}_+\ {\rm and }\ Q^{[s_\tau]}_-\ {\rm are}\ {\rm acting}\ {\rm on}\ F_{\alpha_1})
\end{align*}
(note that $F_{\alpha_2}=F_{\alpha_1+r_\tau\alpha_+}$ or $F_{\alpha_1+s_\tau\alpha_-}$). We define $\widetilde{Q}_\tau=Q^{[r_\tau]}_++Q^{[s_\tau]}_-$ and $\widetilde{\Lambda}_{\tau}=\Lambda_{Q^{[r_\tau]}_+}+\Lambda_{Q^{[s_\tau]}_-}$.
\end{enumerate}
\end{dfn}

We define the following indecomposable Virasoro modules by gluing ceratin Fock modules.
\begin{dfn}
\label{seealso}
\begin{enumerate}
\item Let $\tau=(\alpha_1,\alpha_2,\alpha_3)\in\mathcal{T}_{p_+,p_-}$. We set 
\begin{equation*}
F(\tau)=F_{\alpha_1}\oplus F_{\alpha_2}.
\end{equation*}
For $A\in U(\mathcal{L})$, we define the following operator $J_{\tau}(A)$ on ${F}(\tau)$ 
\begin{equation*}
J_{\tau}(A)=
\begin{cases}
A+\Lambda_\tau(A)\ \ & {\rm on}\ \ F_{\alpha_1},\\
A\ \ & {\rm on}\ \ F_{\alpha_2}.
\end{cases}
\end{equation*}
Then, by Proposition \ref{deri}, we have 
\begin{align*}
J_{\tau}(AB)&=J_{\tau}(A)J_{\tau}(B),\ \ \ {\rm for}\ {\rm any}\ A,B\in U(\mathcal{L}).
\end{align*}
Thus we see that $J_{\tau}$ defines a structure of Virasoro module on ${F}(\tau)$. 
We denote this Virasoro module by $(F(\tau),J_{\tau})$.
\item Let $\tau=(\alpha_1,\alpha_2,\alpha_3)\in\mathcal{T}^{\rm br}_{p_+,p_-}$. We set 
\begin{equation*}
\widetilde{F}(\tau)=F_{\alpha_1}\oplus F_{\alpha_1+r_\tau\alpha_{+}}\oplus F_{\alpha_1+s_\tau\alpha_-}.
\end{equation*}
For $A\in U(\mathcal{L})$, we define the following operator $\widetilde{J}_{\tau}(A)$ on $\widetilde{F}(\tau)$ 
\begin{equation*}
\widetilde{J}_{\tau}(A)=
\begin{cases}
A+\widetilde{\Lambda}_\tau(A)\ \ & {\rm on}\ \ F_{\alpha_1},\\
A\ \ & {\rm on}\ \ F_{\alpha_1+r_\tau\alpha_{+}}\oplus F_{\alpha_1+s_\tau\alpha_-}.
\end{cases}
\end{equation*}
Then, by Proposition \ref{deri}, we have 
\begin{align*}
\widetilde{J}_{\tau}(AB)&=\widetilde{J}_{\tau}(A)\widetilde{J}_{\tau}(B),\ \ \ {\rm for}\ {\rm any}\ A,B\in U(\mathcal{L}).
\end{align*}
Thus we see that $\widetilde{J}_{\tau}$ defines a structure of Virasoro module on $\widetilde{F}(\tau)$. 
We denote this Virasoro module by $(\widetilde{F}(\tau),\widetilde{J}_{\tau})$.
\end{enumerate}
\end{dfn}

The Virasoro modules $(F(\tau),J_{\tau})$ and $(\widetilde{F}(\tau),\widetilde{J}_{\tau})$ have $L_0$-nilpotent rank two. In fact we have the following lemma.
\begin{lem}
\label{L0nilpotent2}
\begin{enumerate}
\item Let $\tau=(\alpha_1,\alpha_2,\alpha_3)\in \mathcal{T}_{p_+,p_-}$. 
The $J_{\tau}(L_n)$ action on the subspace $F_{\alpha_1}\subset F({\tau})$ is given by
\begin{align*}
J_{\tau}(L_n)=L_n+[Q_\tau,a_n]
\end{align*}
\item Let $\tau=(\alpha_1,\alpha_2,\alpha_3)\in \mathcal{T}^{\rm br}_{p_+,p_-}$. 
The $\widetilde{J}_{\tau}(L_n)$ action on the subspace $F_{\alpha_1}\subset \widetilde{F}({\tau})$ is given by
\begin{align*}
\widetilde{J}_{\tau}(L_n)=L_n+[\widetilde{Q}_\tau,a_n].
\end{align*}
\end{enumerate}
\begin{proof}
We only prove the first case. The second case can be proved in the same way.

Note that the ordinary action of $L_n$ on the Fock modules in $\mathcal{F}_{\alpha_0}\mathchar`-{\rm mod}$ is given by
\begin{align}
\label{L0}
L_n=\frac{1}{2}\sum_{m\in\mathbb{Z}}:a_ma_{n-m}:-\frac{1}{2}\alpha_0(n+1)a_n.
\end{align}
Then, by (\ref{conjugate}), (\ref{L0}) and $[Q_\tau,L_n]=0$, we have
\begin{align*}
J_{\tau}(L_n)-L_n
&=\Lambda_\tau(L_n)\\
&=\lim_{x\rightarrow 0}\frac{1}{x}Q_\tau e^{-x\hat{a}}L_ne^{x\hat{a}}-\lim_{x\rightarrow 0}\frac{1}{x}e^{-x\hat{a}}L_ne^{x\hat{a}}Q_\tau\\
&=\lim_{x\rightarrow 0}\frac{1}{x}Q_\tau xa_n-\lim_{x\rightarrow 0}\frac{1}{x}xa_nQ_\tau\\
&=[Q_\tau,a_n].
\end{align*}
\end{proof}
\end{lem}

\begin{prop}
\label{2023P}
For any $\tau=(\alpha_1,\alpha_2,\alpha_3)\in \mathcal{T}^{{\rm Min}}_{p_+,p_-}$, we have
\begin{align*}
{\rm S}(P(\tau))=L(h_{\alpha_2}).
\end{align*}
In particular $P(\tau)$ is self-contragredient.
\begin{proof}
Note that $\mathcal{T}^{{\rm Min}}_{p_+,p_-}$ is given by the following set
\begin{align*}
\{(\alpha_{r^\vee,s^\vee},\alpha_{r,s^\vee;1},\alpha_{r^\vee,s^\vee;2}),
(\alpha_{r^\vee,s^\vee},\alpha_{r^\vee,s;-1},\alpha_{r^\vee,s^\vee;-2})|1\leq r<p_+,1\leq s<p_-\}.
\end{align*}
Fix integers $r,s$ satisfying $1\leq r<p_+,1\leq s<p_-$ and set
\begin{align*}
&\tau=(\alpha_{r^\vee,s^\vee},\alpha_{r,s^\vee;1},\alpha_{r^\vee,s^\vee;2}),
&\tau'=(\alpha_{r^\vee,s^\vee},\alpha_{r^\vee,s;-1},\alpha_{r^\vee,s^\vee;-2}).
\end{align*}
Let us show
\begin{align}
\label{show0517}
{\rm S}(P(\tau))={\rm S}(P(\tau'))=L(\Delta^+_{r,s;0})
\end{align}
(for the definition of $\Delta^+_{r,s;0}$, see Definition~\ref{symbol}).
Consider the vector 
\begin{align*}
e^{-x\hat{a}}S_{r^\vee,s^\vee}e^{x\hat{a}}{\mid}\alpha_{r^\vee,s^\vee}\rangle\in F_{r^\vee,s^\vee},\quad x\in\mathbb{C}^\times,
\end{align*}
where we regard $F_{\alpha_{r^\vee,s^\vee}+x}$ as an object in $\mathcal{F}_{\alpha_0}\mathchar`-{\rm mod}$. Write
\begin{align*}
e^{-x\hat{a}}S_{r^\vee,s^\vee}e^{x\hat{a}}{\mid}\alpha_{r^\vee,s^\vee}\rangle=\sum_{\lambda\vdash r^\vee s^\vee}f_{\lambda}(x)a_{-\lambda}{\mid}\alpha_{r^\vee,s^\vee}\rangle,
\end{align*}
where $f_{\lambda}(x)$ are some polynomials of $x$. 
Since $S_{r^\vee,s^\vee}{\mid}\alpha_{r^\vee,s^\vee}\rangle=0$, we see that every $f_{\lambda}(x)$ is divisible by $x$.
Then we define
\begin{align*}
u_{\tau}:=\sum_{\lambda\vdash r^\vee s^\vee}(x^{-1}f_{\lambda}(x)\mid_{x=0})a_{-\lambda}{\mid}\alpha_{r^\vee,s^\vee}\rangle\in F_{r^\vee,s^\vee}[r^\vee s^\vee].
\end{align*}
By the Jantzen filtration of the Fock module $F_{r^\vee,s^\vee}$ (cf. \cite{FF,IK}), this vector is non-zero.
By the argument in \cite{Felder}, we see that
\begin{align*}
&Q^{[r^\vee]}_+(u_\tau)\in \C^\times\ket{\alpha_{r,s^\vee;1}},
&Q^{[s^\vee]}_-(u_\tau)\in \C^\times\ket{\alpha_{r^\vee,s;-1}}.
\end{align*} 
Consider the logarithmic module $(\widetilde{F}(\tau),\widetilde{J}_{\tau})$.
By Lemma \ref{L0nilpotent2}, we have
\begin{align}
\label{202211221}
(\widetilde{J}_{\tau}(L_0)-h_{r,s^\vee;1})u_\tau=-r^\vee Q^{[r^\vee]}_+(u_\tau)-s^\vee Q^{[s^\vee]}_-(u_\tau).
\end{align}
in $(\widetilde{F}(\tau),\widetilde{J}_{\tau})$. By using Proposition \ref{Y} we have
\begin{align}
\label{202211222}
\widetilde{J}_{\tau}(S_{r^\vee,s^\vee}\sigma(S_{r^\vee,s^\vee}))u_\tau=\frac{1}{2}(2\alpha_{r^\vee,s^\vee}-\alpha_0)R_{r^\vee,s^\vee}(Q^{[r^\vee]}_+(u_\tau)+Q^{[s^\vee]}_-(u_\tau))
\end{align}
in $(\widetilde{F}(\tau),\widetilde{J}_{\tau})$.
Define the following finite length submodule of $(\widetilde{F}(\tau),\widetilde{J}_{\tau})$
\begin{align*}
E_\tau:=\widetilde{J}_{\tau}(U(\mathcal{L})).u_\tau.
\end{align*}
By (\ref{202211221}) and (\ref{202211222}), $E_\tau$ has a indecomposable quotient in
\begin{align*}
{\rm Ext}^1_{\mathcal{L}}(\widetilde{K}(\Delta^+_{r,s;0}), L(\Delta^+_{r,s;0})\oplus L(\Delta^+_{r,s;0}))
\end{align*}
(for the definition of $\widetilde{K}(\Delta^+_{r,s;0})$, see Definition \ref{dfnC2}).
We denote this quotient by $\widetilde{P}(\Delta^+_{r,s;0})$. Note that by Proposition \ref{sankakuthick}, $\widetilde{P}(\Delta^+_{r,s;0})$ has $P(\tau)$ and $P(\tau')$,
as quotients.
Consider the quotient module
$
E_\tau/\widetilde{J}_{\tau}(U(\mathcal{L})).\ket{\alpha_{r^\vee,s^\vee}}.
$
By (\ref{202211221}) and (\ref{202211222}), this quotient module has $L_0$ nilpotent rank two and has a logarithmic quotient isomorphic to $\widetilde{P}(\Delta^+_{r,s;0})/L(h_{r,s},\Delta^+_{r,s;0})$. We set $R(\Delta^+_{r,s;0})=\widetilde{P}(\Delta^+_{r,s;0})/L(h_{r,s},\Delta^+_{r,s;0})$.
Then, by Propositions \ref{Vir0} and \ref{sankakuthick}, we have
$
{\rm S}(R(\Delta^+_{r,s;0}))=L(\Delta^+_{r,s;0}).
$
Therefore we obtain (\ref{show0517}).
\end{proof}
\end{prop}

\begin{remark}
By Proposition~\ref{2023P} and Remark~\ref{rem:log-tau0}, $P(\tau)$ is self-contragredient for $\tau\in \mathcal{T}^0_{p_+,p_-}\sqcup\mathcal{T}^{\rm Min}_{p_+,p_-}$. We expect this to hold for all $\tau \in \mathcal{T}_{p_+,p_-}$.
\end{remark}

\begin{dfn}
\label{indP1}
Let $\tau=(\alpha_1,\alpha_2,\alpha_3)\in \mathcal{T}^0_{p_+,p_-}\sqcup\mathcal{T}^{\rm Min}_{p_+,p_-}$. By Proposition \ref{sankakuthick}, we have
%2023P
\begin{align*}
{\rm Ext}^1_{\mathcal{L}}(K(\tau),L(h_{\alpha_2},h_{\alpha_3})^*)\simeq \C.
\end{align*}
Let $\widehat{P}(\tau)$ be an indecomposable module representing the nonzero class in the above ${\rm Ext}^1$-group.
\end{dfn}
From the definition of $\widehat{P}(\tau)$, we can easily see that
\begin{equation}
\label{agi}
\begin{split}
&{\rm dim}_{\C}{\rm Hom}_{U(\mathcal{L})}(L(h_{\alpha_2},h_{\alpha_3})^*,\widehat{P}(\tau))= 2,\\
&{\rm dim}_{\C}{\rm Hom}_{U(\mathcal{L})}(\widehat{P}(\tau),K(\tau))= 2.
\end{split}
\end{equation}
Note that $\widehat{P}(\tau)$ has a quotient isomorphic to ${P}(\tau)$. In fact, for any $L_0$-homogeneous vector $u$ not annihilated by the action of $L_0-h_{\alpha_2}$, the module $U(\mathcal{L}).u$ is isomorphic to ${P}(\tau)$.
We also note that this logarithmic module $\widehat{P}(\tau)$ can be realized as the quotient of $(F(\tau),J_\tau)$ or $(\widetilde{F}(\tau),\widetilde{J}_{\tau})$.
The following proposition is clear from the self-dualities $P(\tau)^*\simeq P(\tau)$ and $\widehat{P}(\tau)^*\simeq \widehat{P}(\tau)$ for $\tau\in \mathcal{T}^0_{p_+,p_-}\sqcup\mathcal{T}^{\rm Min}_{p_+,p_-}$, but we give a direct proof using the structure of $(F(\tau),J_\tau)$ and $(\widetilde{F}(\tau),\widetilde{J}_\tau)$.
\begin{prop}
\label{sankaku0000}
Let $\tau=(\alpha_1,\alpha_2,\alpha_3)\in \mathcal{T}^0_{p_+,p_-}\sqcup\mathcal{T}^{\rm Min}_{p_+,p_-}$.
\begin{enumerate}
\item Two injections can be chosen as a basis of the Hom-space
\begin{align}
\label{2basi}
{\rm Hom}_{U(\mathcal{L})}(L(h_{\alpha_2},h_{\alpha_3})^*,\widehat{P}(\tau))\quad (\simeq \mathbb{C}^2).
\end{align}
\item Two surjections can be chosen as a basis of the Hom-space
\begin{align}
\label{2bas}
{\rm Hom}_{U(\mathcal{L})}(\widehat{P}(\tau),K(\tau))\quad (\simeq \mathbb{C}^2).
\end{align}
\end{enumerate}

\begin{proof}
We only prove in the case of $\tau=(\alpha_{r^\vee,p_-;-1},\alpha_{r,p_-},\alpha_{r^\vee,p_-;1})\in \mathcal{T}^0_{p_+,p_-}$. The case where $\tau\in \mathcal{T}^{\rm Min}_{p_+,p_-}$ can be proved similarly by considering the structure of $(\widetilde{F}(\tau),\widetilde{J}_\tau)$.

Consider the logarithmic module $(F(\tau),J_\tau)$.
Set $u_0=\ket{\alpha_{r^\vee,p_-;-1}}$ and $u_1=\ket{\alpha_{r,p_-}}$. Fix an arbitrary $L_0$-homogeneous vector $v_1\in F_{r^\vee,p_-;-1}[rp_-]$ such that $Q^{[r]}_+(v_1)\in \C^\times \ket{\alpha_{r^\vee,p_-;1}}$, and let $v_2=Q^{[r]}_+(\ket{\alpha_{r,p_-;-2}})$. Note that 
\begin{align}
\label{20230401-1eq}
\sigma(S_{r,p_-})v_1\in \C^\times \ket{\alpha_{r,p_-}}.
\end{align}
By the Jantzen filtration of the Fock module $F_{r^\vee,p_-;-1}$ (cf. \cite{IK}), we have
\begin{align*}
e^{-x\hat{a}}\sigma(S_{r,p_-})e^{x\hat{a}}v_2=xf(x) u_0,
\end{align*}
where $f(x)$ is a nonzero polynomial in $\mathbb{C}[x]\setminus x\mathbb{C}[x]$.
Then, by noting the definition of $J_\tau$, we obtain
\begin{equation}
\label{202304010eq-0}
\begin{split}
J_\tau(\sigma(S_{r,p_-}))v_2=\lim_{\substack{x\rightarrow 0\\ x\neq 0}}\frac{1}{x}Q^{[r^\vee]}_+e^{-x\hat{a}}\sigma(S_{r,p_-})e^{x\hat{a}}v_2\in \C^\times u_1
\end{split}
\end{equation}
Similarly, from the Jantzen filtration of the Fock module $F_{r,p_-}$, we have
\begin{equation}
\label{202304010eq}
\begin{split}
J_\tau(S_{r,p_-})u_0\in \mathbb{C}^\times v_2+\mathbb{C}^\times v_1+U(\mathcal{L}).u_1.
\end{split}
\end{equation}
(see Figure~\ref{51fig} for these relations).
Consider the submodule
\begin{align*}
I=J_\tau(U(\mathcal{L})).u_0+J_\tau(U(\mathcal{L})).v_1.
\end{align*}
By definition, $I$ has a quotient isomorphic to $\widehat{P}(\tau)$. Let $\pi$ be the surjection from $I$ to $\widehat{P}(\tau)$, and define the following injections and surjections 
\begin{align*}
&\iota_i: L(h_{r,p_-},h_{r^\vee,p_-;1})^*\stackrel{\simeq}{\longrightarrow}\pi(J_\tau(U(\mathcal{L})).v_i)\subset \widehat{P}(\tau),\ \ i=1,2,\\
&p_j:\widehat{P}(\tau)\rightarrow \widehat{P}(\tau)/\pi(J_\tau(U(\mathcal{L})).v_j),\ \ j=1,2.
\end{align*}
Then, from (\ref{agi}), (\ref{20230401-1eq}), (\ref{202304010eq-0}), and (\ref{202304010eq}), we see that $\{\iota_1,\iota_2\}$ and $\{p_1,p_2\}$ give a basis of (\ref{2basi}) and (\ref{2bas}), respectively.
\end{proof}
\end{prop}

\begin{figure}[htbp]
  \centering
\begin{tikzpicture}[scale=0.85]
\node[inner sep=0.7pt] (a1) at (0, 0) {${\pi(u_0)}$};
\node[inner sep=0.7pt] (a2) at (3, 0) {${\pi(u_1)}$};
\node[inner sep=0.7pt] (b1) at (0, -3) {${\pi(v_2)}$};
\node[inner sep=0.7pt] (b2) at (3, -3) {${\pi(v_1)}$};

\draw[arrows = {-Stealth[scale=0.9]}] (a1) to (b1);
\draw[arrows = {-Stealth[scale=0.9]},very thick] (a1) to (b2);
\draw[arrows = {-Stealth[scale=0.9]},very thick] (a1) to (b1);
\draw[arrows = {-Stealth[scale=0.9]},dashed] (a1) to (a2);
\draw[arrows = {-Stealth[scale=0.9]}] (b1) to (a2);
\draw[arrows = {-Stealth[scale=0.9]}] (b2) to (a2);
\end{tikzpicture}
\caption{\label{51fig}The schematic diagram of the logarithmic modules $\widehat{P}(\tau)$ $(\tau=(\alpha_{r^\vee,p_-;-1},\alpha_{r,p_-},\alpha_{r^\vee,p_-;1})\in \mathcal{T}^0_{p_+,p_-})$ described in the proof of Proposition \ref{sankaku0000}. 
The two thick arrows represent that the vectors $\pi (v_i)$ $(i=1,2)$ cannot be obtained independently from the vector $\pi(u_0)$ at the starting point of the arrows by the action of $U(\mathcal{L})$.
The meaning of the remaining arrows is the same as in Figures \ref{51figk}-\ref{51figkk}.}
\mbox{}\\
\end{figure}

As in Definition \ref{indP1}, we define the following indecomposable modules.
\begin{dfn}
\begin{enumerate}
\item 
Fix $\tau=(\alpha_1,\alpha_2,\alpha_3)\in \mathcal{T}_{p_+,p_-}\setminus (\mathcal{T}^0_{p_+,p_-}\sqcup\mathcal{T}^{\rm Min}_{p_+,p_-})$. By Proposition \ref{sankakuthick}, we have
%2023P
\begin{align*}
{\rm Ext}^1_{\mathcal{L}}(K(\tau),L(h_{\alpha_1},h_{\alpha_2}))\simeq \C.
\end{align*}
Let $\widehat{P}^{\downarrow}(\tau)$ be an indecomposable module representing the nonzero class in the above ${\rm Ext}^1$-group.
\item Fix $1<r<p_+$, $1<s<p_-$and $n\geq 2$. By Proposition \ref{sankakuthick}, we have
\begin{align*}
&{\rm Ext}^1_{\mathcal{L}}(\widetilde{K}(h_{r,s;n}),L(h_{r^\vee,s;n-1},h_{r,s;n}))\simeq \mathbb{C},\\
&{\rm Ext}^1_{\mathcal{L}}(\widetilde{K}(h_{r,s;n}),L(h_{r,s^\vee;n-1},h_{r,s;n}))\simeq \mathbb{C}.
\end{align*}
Thus indecomposable modules in
\begin{align*}
{\rm Ext}^1_{\mathcal{L}}(\widetilde{K}(h_{r,s;n}),L(h_{r^\vee,s;n-1},h_{r,s;n})\oplus L(h_{r,s^\vee;n-1},h_{r,s;n}) )
\end{align*}
are unique up to isomorphism. 
Let $\widetilde{P}^{\downarrow}(h_{r,s;n})$ be such a module.
\end{enumerate}
\label{dfn20231225}
\end{dfn}

\mbox{}\\
Note that $\widehat{P}^{\downarrow}(\tau)$ and $\widetilde{P}^{\downarrow}(h_{r,s;n})$ have quotients isomorphic to ${P}(\tau)$ and $\widetilde{P}(h_{r,s;n})$, respectively.
By definition and Proposition~\ref{VirC}, we see that
\begin{align}
&{\rm Hom}_{U(\mathcal{L})}(\widehat{P}^{\downarrow}(\tau),K(\tau))\simeq \C^2,\label{2bas7}\\
&{\rm Hom}_{U(\mathcal{L})}(\widetilde{P}^{\downarrow}(h_{r,s;n}),\widetilde{K}(h_{r,s;n}))\simeq \C^3.\label{2bas77}
\end{align}
The following proposition can be proved in a simlar way as Proposition \ref{sankaku0000}, by using the structure of the logarithmic modules $(F(\tau),J_{\tau})$ and $(\widetilde{F}(\tau),\widetilde{J}_{\tau})$ (see also the schematic diagrams in Figures~\ref{51figkkk} and \ref{51figkkkk}).
\begin{prop}
\begin{enumerate}
\item
Let $\tau=(\alpha_1,\alpha_2,\alpha_3)\in \mathcal{T}_{p_+,p_-}\setminus (\mathcal{T}^0_{p_+,p_-}\sqcup\mathcal{T}^{\rm Min}_{p_+,p_-})$.
Then, two surjections can be chosen as a basis of the Hom-space (\ref{2bas7}).
\item
Let $1<r<p_+$, $1<s<p_-$and $n\geq 2$.
Then, three surjections can be chosen as a basis of the Hom-space (\ref{2bas77}).
\end{enumerate}
\label{bt3}
\end{prop}

\begin{figure}[htbp]
  \centering
\begin{tikzpicture}[scale=0.7]

\node[inner sep=0.7pt] (p1) at (0, -1.5) {${\ket{h_{\alpha_2}}}$};
\node[inner sep=0.7pt] (p2) at (1.5, 0) {${\ket{h_{\alpha_1}}}$};
\node[inner sep=0.7pt] (q1) at (1.5, -3) {${\ket{h_{\alpha_3}}}$};
\node[inner sep=0.7pt] (t) at (3, -1.5) {${\ket{h_{\alpha_2}}}$};
\node[inner sep=0.7pt] (v) at (3, 0) {${\ket{h_{\alpha_1}}}$};
\draw[arrows = {-Stealth[scale=0.9]}] (p1) to (q1);
\draw[arrows = {-Stealth[scale=0.9]}] (p1) to (q1);
\draw[arrows = {-Stealth[scale=0.9]},very thick] (p1) to (p2);
\draw[arrows = {-Stealth[scale=0.9]},dotted,shift={(7,0)}] (q1) to (t);
\draw[arrows = {-Stealth[scale=0.9]}] (p2) to (t);
\draw[arrows = {-Stealth[scale=0.9]},dashed] (p1) to (t);
\draw[arrows = {-Stealth[scale=0.9]},very thick] (p1) to (v);
\draw[arrows = {-Stealth[scale=0.9]}] (v) to (t);
\end{tikzpicture}
\caption{The schematic diagram of the logarithmic module $\widehat{P}^{\downarrow}(\tau)$ defined in Definition \ref{dfn20231225}. 
Each ket vector $\ket{h_{\alpha_\bullet}}$ represents the lowest weight vector of the simple module $L(h_{\alpha_\bullet})$ in some quotient.
The two thick arrows represent that the above ket vectors $\ket{h_{\alpha_1}}$ cannot be obtained independently from the ket vector $\ket{h_{\alpha_2}}$ at the starting point of the arrows by the action of $U(\mathcal{L})$ (see (\ref{2bas7})). The meaning of the remain arrows is the same as in Figure \ref{51figk}.  \label{51figkkk}}
\mbox{}\\
\end{figure}

\begin{figure}[htbp]
  \centering
\begin{tikzpicture}[scale=0.7]

\node[inner sep=0.7pt,shift={(8.5,0)}] (dd) at (-4.5, 0) {${\ket{h_{r^\vee,s;n-1}}}$};
\node[inner sep=0.7pt,shift={(8.5,0)}] (ddd) at (4.5, 0) {${\ket{h_{r^\vee,s;-n+1}}}$};
\node[inner sep=0.7pt,shift={(8.5,0)}] (d2) at (-1.5, 0) {${\ket{h_{r^\vee,s;n-1}}}$};
\node[inner sep=0.7pt,shift={(8.5,0)}] (f1) at (-1.5, -3) {${\ket{h_{r^\vee,s;n+1}}}$};
\node[inner sep=0.7pt,shift={(8.5,0)}] (s) at (-3, -1.5) {${\ket{h_{r,s;n}}}$};

\node[inner sep=0.7pt,shift={(8.5,0)}] (pp1) at (0, -1.5) {${\ket{h_{r,s;n}}}$};
\node[inner sep=0.7pt,shift={(8.5,0)}] (pp2) at (1.5, 0) {${\ket{h_{r^\vee,s;-n+1}}}$};
\node[inner sep=0.7pt,shift={(8.5,0)}] (qq1) at (1.5, -3) {${\ket{h_{r^\vee,s-n-1}}}$};
\node[inner sep=0.7pt,shift={(8.5,0)}] (tt) at (3, -1.5) {${\ket{h_{r,s;n}}}$};

\draw[arrows = {-Stealth[scale=0.9]},very thick,shift={(7,0)}] (pp1) to (dd);
\draw[arrows = {-Stealth[scale=0.9]},very thick,shift={(7,0)},dotted,] (pp1) to (ddd);
\draw[arrows = {-Stealth[scale=0.9]},shift={(7,0)}] (dd) to (s);
\draw[arrows = {-Stealth[scale=0.9]},shift={(7,0)}] (ddd) to (tt);

\draw[arrows = {-Stealth[scale=0.9]},shift={(7,0)}] (pp1) to (f1);
\draw[arrows = {-Stealth[scale=0.9]},very thick,shift={(7,0)}] (pp1) to (d2);
\draw[arrows = {-Stealth[scale=0.9]},dotted,shift={(7,0)}] (f1) to (s);
\draw[arrows = {-Stealth[scale=0.9]},shift={(7,0)}] (d2) to (s);
\draw[arrows = {-Stealth[scale=0.9]},dashed,shift={(7,0)}] (pp1) to (s);

\draw[arrows = {-Stealth[scale=0.9]},shift={(7,0)}] (pp1) to (qq1);
\draw[arrows = {-Stealth[scale=0.9]},shift={(7,0)}] (pp1) to (qq1);
\draw[arrows = {-Stealth[scale=0.9]},very thick,shift={(7,0)},dotted,] (pp1) to (pp2);
\draw[arrows = {-Stealth[scale=0.9]},dotted,shift={(7,0)}] (qq1) to (tt);
\draw[arrows = {-Stealth[scale=0.9]},shift={(7,0)}] (pp2) to (tt);
\draw[arrows = {-Stealth[scale=0.9]},dashed,shift={(7,0)}] (pp1) to (tt);

\end{tikzpicture}
\caption{The schematic diagram of the logarithmic module $\widetilde{P}^{\downarrow}(h_{r,s;n})$ defined in Definition \ref{dfn20231225}. 
Each ket vector $\ket{h_{\bullet,\bullet;\bullet}}$ represents the lowest weight vector of the simple module $L(h_{\bullet,\bullet;\bullet})$ in some quotient.
The two thick arrows represent that the above ket vectors $\ket{h_{r^\vee,s;n-1}}$ cannot be obtained independently from the ket vector $\ket{h_{r,s;n}}$ at the starting point of the arrows by the action of $U(\mathcal{L})$ (see (\ref{2bas77})). 
The two thick dotted arrows are interpreted in the same way.
The meaning of the remain arrows is the same as in Figure \ref{51figk}.  \label{51figkkkk}}
\mbox{}\\
\end{figure}

\section{Projective covers of simple $\W_{p_+,p_-}$-modules}
\label{section proj}
Since $\W_{p_+,p_-}$ is $C_2$-cofinite, it follows from \cite{H} that every simple $\W_{p_+,p_-}$-module admits a projective cover. In this section, we determine some ${\rm Ext}^1$-groups between certain indecomposable modules and simple modules. Based on these ${\rm Ext}^1$ groups, we determine the structure of projective covers of the simple modules in each thick block and thin block.
From this section, we denote the ${\rm Ext}^1$-groups in $\mathcal{C}_{p_+,p_-}$ as ${\rm Ext}^1(\bullet,\bullet)$ simply and identify any $\W_{p_+,p_-}$-modules that are isomorphic among each other.

\subsection{Structure of the logarithmic modules $\mathcal{Q}(\X^\pm_{\bullet,\bullet})_{\bullet,\bullet}$}
\label{structureQ}
Fix integers $r,s\in \Z^2_{\geq 1}$ satisfying $1\leq r< p_+$ and $1\leq s< p_-$. In this subsection we consider the structure of the indecomposable modules $\mathcal{Q}(\X^\pm_{\bullet,\bullet})_{\bullet,\bullet}$ in the blocks $C^{thin}_{r,p_-}$, $C^{thin}_{p_+,s}$ and $C^{thick}_{r,s}$. 
First let us consider the structure of the logarithmic module $\mathcal{Q}(\X^+_{r,p_-})_{r^\vee,p_-}\in C^{thin}_{r,p_-}$. Recall that as the vector space $\mathcal{Q}(\X^+_{r,p_-})_{r^\vee,p_-}=\V^-_{r^\vee,p_-}\oplus \V^+_{r,p_-}$. 
Fix an arbitrary $v^+_{1}\in F_{r,p_-}[rp_-]$ satisfying $Q^{[r]}_+(v^+_{1})\in\mathbb{C}^\times \ket{\alpha_{r^\vee,p_-;1}}$, and set
\begin{align}
\label{relsym1}
\begin{aligned}
u_0&=\ket{\alpha_{r^\vee,p_-;-1}},&u_1&=\ket{\alpha_{r,p_-}},\\
v^+_2&=\ket{\alpha_{r^\vee,p_-;1}},&v^-_1&=\ket{\alpha_{r,p_-;-2}},&v^-_2&=Q^{[r]}_+(\ket{\alpha_{r,p_-;-2}}).
\end{aligned}
\end{align}
Note that $\mathcal{W}_{p_+,p_-}.u_1\simeq \mathcal{X}^+_{r,p_-}$ (see Definition \ref{dfn0830}). 
Consider the following Virasoro submodule of $\mathcal{Q}(\X^+_{r,p_-})_{r^\vee,p_-}$
\begin{align*}
I=U(\mathcal{L}).u_0+U(\mathcal{L}).v^+_1\ \in\mathcal{L}_{c_{p_+,p_-}}\mathchar`-{\rm mod},
\end{align*}
where we have omitted the symbol $J^+_{r,p_-}$ associated to the module action defined in Proposition~\ref{action0303}.
Note that 
\begin{align}
\label{renketsu2023}
&(L_0-\Delta^+_{r,p_-;0})u_0\in\C^\times u_1,
&\sigma(S_{r,p_-})v^+_1\in\C^\times u_1
\end{align}
in $\mathcal{Q}(\X^+_{r,p_-})_{r^\vee,p_-}$.
Then, by (\ref{renketsu2023}), we see that $I$ has a quotient isomorphic to $\widehat{P}(\tau)$, where $\tau=(\alpha_{r^\vee,p_-;-1},\alpha_{r,p_-},\alpha_{r^\vee,p_-;1})$.
Thus, by Proposition \ref{sankaku0000}, we have
\begin{align}
S_{r,p_-}{u_0}&\in \mathbb{C}^\times {v^+_1}+\C^\times {v^-_2}+U(\mathcal{L}).u_1,\label{20221021KK-0}\\
\sigma(S_{r,p_-}) v^-_2&\in \C^\times u_1\label{20221021KK},
\end{align}
in $\mathcal{Q}(\mathcal{X}^+_{r,p_-})_{r^\vee,p_-}$. 
By (\ref{20221021KK-0}), we have
\begin{align}
\label{20221021KK-01}
\begin{aligned}
&W^+[0]W^-[0]S_{r,p_-}{u_0}&\in \mathbb{C}^\times {v^+_1}+U(\mathcal{L}).u_1,\\
&W^-[0]W^+[0]S_{r,p_-}{u_0}&\in \mathbb{C}^\times {v^-_2}+U(\mathcal{L}).u_1
\end{aligned}
\end{align}
(see Figure~\ref{fig61} for a schematic illustration of (\ref{20221021KK-01})).
Then, by (\ref{20221021KK-01}) and (\ref{20221021KK}), we see that $\mathcal{Q}(\X^+_{r,p_-})_{r^\vee,p_-}$ has two submodules $\W_{p_+,p_-}.v^+_1$ and $\W_{p_+,p_-}.v^-_2$ satisfying
\begin{equation}
\label{generate2023}
\begin{split}
&[\W_{p_+,p_-}.v^+_1]\in {\rm Ext}^1(\X^-_{r^\vee,p_-},\X^+_{r,p_-})\setminus\{0\},\\
&[\W_{p_+,p_-}.v^-_2]\in {\rm Ext}^1(\X^-_{r^\vee,p_-},\X^+_{r,p_-})\setminus\{0\}.
\end{split}
\end{equation}
By (\ref{generate2023}), we see that $\mathcal{Q}(\X^+_{r,p_-})_{r^\vee,p_-}$ is generated from the top composition factor $\X^+_{r,p_-}$ and ${\rm S}(\mathcal{Q}(\X^+_{r,p_-})_{r^\vee,p_-})=\X^+_{r,p_-}$.

We have similar results for the other indecomposable module $\mathcal{Q}(\X^+_{p_+,s})_{p_+,s^\vee}$ in $C^{thin}_{p_+,s}$. Thus we obtain the following proposition.
\begin{prop}
The logarithmic modules $\mathcal{Q}(\X^+_{r,p_-})_{r^\vee,p_-}$ and $\mathcal{Q}(\X^+_{p_+,s})_{p_+,s^\vee}$ are generated from the top composition factors and have the following length three socle series:
\begin{enumerate}
\item The socle series of $\mathcal{Q}(\X^+_{r,p_-})_{r^\vee,p_-}$ is given by
\begin{align*}
\hspace{-6mm}{\rm S}_1=\X^+_{r,p_-},\quad {\rm S}_2/{\rm S}_1=2\X^{-}_{r^\vee,p_-},\quad \mathcal{Q}(\X^+_{r,p_-})_{r^\vee,p_-}/{\rm S}_2=\X^+_{r,p_-}.
\end{align*}
\item The socle series of $\mathcal{Q}(\X^+_{p_+,s})_{p_+,s^\vee}$ is given by
\begin{align*}
\hspace{-6mm}{\rm S}_1=\X^+_{p_+,s},\quad {\rm S}_2/{\rm S}_1=2\X^{-}_{p_+,s^\vee},\quad \mathcal{Q}(\X^+_{p_+,s})_{p_+,s^\vee}/{\rm S}_2=\X^+_{p_+,s}.
\end{align*}
\end{enumerate}
\label{thinst}
\end{prop}

\begin{figure}[htbp]
  \centering
\begin{tikzpicture}[scale=0.85]
\node[inner sep=0.7pt] (a1) at (0, 0) {${u_0}$};
\node[inner sep=0.7pt] (b1) at (-2, -2) {${v^+_1}$};
\node[inner sep=0.7pt] (b2) at (2, -2) {${v^-_2}$};
\node[inner sep=0.7pt] (c1) at (0,-4) {${u_1}$};
\node[inner sep=0.7pt] (d1) at (-4.5,-2) {${v^-_1}$};
\node[inner sep=0.7pt] (d2) at (4.5,-2) {${v^+_2}$};
\node[inner sep=0.7pt,scale=0.9] (k1) at (-3.25,-2+0.6) {${W^+_{0}}$};
\node[inner sep=0.7pt,scale=0.9] (k2) at (-3.25,-2-0.6) {${W^-_{0}}$};
\node[inner sep=0.7pt,scale=0.9] (l1) at (3.25,-2+0.6) {${W^+_{0}}$};
\node[inner sep=0.7pt,scale=0.9] (l2) at (3.25,-2-0.6) {${W^-_{0}}$};
\node[inner sep=0.7pt,scale=0.9] (j1) at (-1.75,-1+0.5) {${W^+_{0}W^-_{0}S_{r,p_-}}$};
\node[inner sep=0.7pt,scale=0.9] (j2) at (-1.75,-2.5-0.6) {${\sigma(S_{r,p_-})}$};
\node[inner sep=0.7pt,scale=0.9] (i1) at (1.75,-1+0.5) {${W^-_{0}W^+_{0}S_{r,p_-}}$};
\node[inner sep=0.7pt,scale=0.9] (i2) at (1.75,-2.5-0.6) {${\sigma(S_{r,p_-})}$};
\draw[arrows = {-Stealth[scale=0.9]}] (a1) to (b1);
\draw[arrows = {-Stealth[scale=0.9]}] (a1) to (b2);
\draw[arrows = {-Stealth[scale=0.9]}] (b1) to (c1);
\draw[arrows = {-Stealth[scale=0.9]}] (b2) to (c1);
\draw[->] (d1) to[bend left=20]  (b1);
\draw[->] (b1) to[bend left=20]  (d1);
\draw[->] (d2) to[bend left=20]  (b2);
\draw[->] (b2) to[bend left=20]  (d2);
\end{tikzpicture}
\caption{A schematic diagram of the $\mathcal{W}_{p_+,p_-}$-action between generators of the composition factors of $\mathcal{Q}(\X^+_{r,p_-})_{r^\vee,p_-}$, where we set $W^\pm_0=W^\pm[0]$ and the symobls correspond to those introduced in (\ref{relsym1}). \label{fig61}}
\end{figure}
\mbox{}\\
Next let us consider the logarithmic module $\mathcal{Q}(\X^-_{r^\vee,p_-})_{r,p_-}$. 
Recall that as the vector space $\mathcal{Q}(\X^-_{r^\vee,p_-})_{r,p_-}=\V^-_{r^\vee,p_-}\oplus \V^+_{r,p_-}$.
Fix an arbitrary $\tilde{v}_+\in F_{r,p_-}[rp_-]$ satisfying $Q^{[r]}_+(\tilde{v}_+)\in\mathbb{C}^\times \ket{\alpha_{r^\vee,p_-;1}}$, and set
\begin{align}
\label{relsym2}
\begin{aligned}
u_-&=\ket{\alpha_{r^\vee,p_-;-1}},&u_+&=\ket{\alpha_{r,p_-}},\\
{v}_+&=\ket{\alpha_{r^\vee,p_-;1}},&v_-&=Q^{[r]}_+(\ket{\alpha_{r,p_-;-2}}),&\tilde{v}_-&=\ket{\alpha_{r,p_-;-2}}.
\end{aligned}
\end{align}
Note that $\mathcal{W}_{p_+,p_-}.v_{+}+\mathcal{W}_{p_+,p_-}.v_{-}\simeq \mathcal{X}^-_{r^\vee,p_-}$ (see Definition \ref{dfn0830}). 
Consider the indecomposable Virasoro modules
$
I^\pm=U(\mathcal{L}).\tilde{v}_\pm.
$
Note that 
\begin{align}
\label{0402eq0}
(L_0-\Delta^-_{r^\vee,p_-;0})\tilde{v}_\pm\in \C^\times v_{\pm},
\end{align}
in $\mathcal{Q}(\X^-_{r^\vee,p_-})_{r,p_-}$ (we omit the symbol $J^-_{r^\vee,p_-}$ associated the module action defined in Proposition~\ref{action0303}).
Then, noting Theorem \ref{KR0}, Proposition (\ref{sankakuthick}) and (\ref{0402eq0}), we see that $I^{\epsilon}(\epsilon=\pm)$ has a quotient isomorphic to $P(\tau)$, where $\tau=(\alpha_{r,p_-},\alpha_{r^\vee,p_-;1},\alpha_{r,p_-;2})$ (for the definition of $P(\tau)$, see Definition \ref{sankakuthickdf}).
Then, from the structure of $P(\tau)$, we obtain
\begin{align}
\label{2022KKK}
S_{r,p_-}\sigma(S_{r,p_-})\tilde{v}_\pm\in\C^\times v_\pm
\end{align}
(see Figure~\ref{fig62} for a schematic illustration of (\ref{2022KKK})).
By (\ref{2022KKK}) we see that $\mathcal{Q}(\X^-_{r^\vee,p_-})_{r,p_-}$ has the submodules $\W_{p_+,p_-}.\sigma(S_{r,p_-})\tilde{v}_\pm$ satisfying
\begin{align*}
[\W_{p_+,p_-}.\sigma(S_{r,p_-})\tilde{v}_\pm]\in {\rm Ext}^1(\X^+_{r,p_-},\X^-_{r^\vee,p_-})\setminus\{0\}.
\end{align*}
In particular we see that $\mathcal{Q}(\X^-_{r^\vee,p_-})_{r,p_-}$ is generated from the top composition factor $\X^-_{r^\vee,p_-}$. We have similar results for the other indecomposable module $\mathcal{Q}(\X^-_{p_+,s^\vee})_{p_+,s}$ in $C^{thin}_{p_+,s}$. Thus we obtain the following theorem.
\begin{prop}
\label{thinst2}
The logarithmic modules $\mathcal{Q}(\X^-_{r^\vee,p_-})_{r,p_-}$ and $\mathcal{Q}(\X^-_{p_+,s^\vee})_{p_+,s}$ are generated from the top composition factors and have the following length three socle series:
\begin{enumerate}
\item The socle series of $\mathcal{Q}(\X^-_{r^\vee,p_-})_{r,p_-}$ is given by
\begin{align*}
\hspace{-6mm}{\rm S}_1=\X^-_{r^\vee,p_-},\quad {\rm S}_2/{\rm S}_1=2\X^{+}_{r,p_-},\quad \mathcal{Q}(\X^-_{r^\vee,p_-})_{r,p_-}/{\rm S}_2=\X^-_{r^\vee,p_-}.
\end{align*}
\item The socle series of $\mathcal{Q}(\X^-_{p_+,s^\vee})_{p_+,s}$ is given by
\begin{align*}
\hspace{-6mm}{\rm S}_1=\X^-_{p_+,s^\vee},\quad {\rm S}_2/{\rm S}_1=2\X^{+}_{p_+,s},\quad \mathcal{Q}(\X^-_{p_+,s^\vee})_{p_+,s}/{\rm S}_2=\X^-_{p_+,s^\vee}.
\end{align*}
\end{enumerate}
\end{prop}

\vspace{3mm}
\begin{figure}[htbp]
  \centering
\begin{tikzpicture}[scale=0.85]
\node[inner sep=0.7pt] (a1) at (0, 0) {${\tilde{v}}_-$};
\node[inner sep=0.7pt] (a2) at (2.5, 0) {${\tilde{v}}_+$};
\node[inner sep=0.7pt] (b1) at (-2.5, -2) {${u_-}$};
\node[inner sep=0.7pt] (b2) at (5, -2) {${u_+}$};
\node[inner sep=0.7pt] (c1) at (0,-4) {${v_-}$};
\node[inner sep=0.7pt] (c2) at (2.5,-4) {${v_+}$};
\node[inner sep=0.7pt,scale=0.9] (k1) at (1.25,0.6) {${W^+[0]}$};
\node[inner sep=0.7pt,scale=0.9] (k2) at (1.25,-0.6) {${W^-[0]}$};
\node[inner sep=0.7pt,scale=0.9] (l1) at (1.25,-4+0.6) {${W^+[0]}$};
\node[inner sep=0.7pt,scale=0.9] (l2) at (1.25,-4-0.6) {${W^-[0]}$};
\node[inner sep=0.7pt,scale=0.9] (j2) at (-1.95,-1+0.2) {${\sigma(S_{r,p_-})}$};
\node[inner sep=0.7pt,scale=0.9] (i2) at (-1.64,-2.5-0.6) {$S_{r,p_-}$};
\node[inner sep=0.7pt,scale=0.9] (j3) at (1.75+2.7,-1+0.2) {${\sigma(S_{r,p_-})}$};
\node[inner sep=0.7pt,scale=0.9] (i3) at (1.75+2.7,-2.5-0.6) {${S_{r,p_-}}$};

\draw[arrows = {-Stealth[scale=0.9]}] (a1) to (b1);
\draw[arrows = {-Stealth[scale=0.9]}] (a2) to (b2);
\draw[arrows = {-Stealth[scale=0.9]}] (b1) to (c1);
\draw[arrows = {-Stealth[scale=0.9]}] (b2) to (c2);
\draw[->] (a1) to[bend left=20]  (a2);
\draw[->] (a2) to[bend left=20]  (a1);
\draw[->] (c1) to[bend left=20]  (c2);
\draw[->] (c2) to[bend left=20]  (c1);
\end{tikzpicture}
\caption{A schematic diagram of the $\mathcal{W}_{p_+,p_-}$-action between generators of the composition factors of $\mathcal{Q}(\X^-_{r^\vee,p_-})_{r,p_-}$, where the symobls correspond to those introduced in (\ref{relsym2}).  \label{fig62}}
\end{figure}
\mbox{}\\

As in Theorems \ref{thinst} and \ref{thinst2}, by using Proposition \ref{sankakuthick} and Proposition \ref{sankaku0000}, we obtain the following propositions for the indecomposable modules $\mathcal{Q}(\X^\pm_{\bullet,\bullet})_{\bullet,\bullet}$ in $C^{thick}_{r,s}$. We omit the proofs.
\begin{prop}
\label{Q+thick}
Let $(a,b,c,d)$ be any element in
\begin{align*}
\{(a,b,c,d)\}=&\bigl\{(r,s,r^\vee,s),(r,s,r,s^\vee),(r^\vee,s^\vee,r^\vee,s),(r^\vee,s^\vee,r,s^\vee)\}.
\end{align*}
Then the socle series of $\mathcal{Q}(\X^+_{a,b})_{c,d}$ is given by
\begin{align*}
\hspace{-3mm}{\rm S}_1=\X^+_{a,b},\quad 
{\rm S}_2/{\rm S}_1=2\X^{-}_{c,d}\oplus L(h_{a,b}),\quad 
\mathcal{Q}(\X^+_{a,b})_{c,d}/{\rm S}_2=\X^+_{a,b}.
\end{align*}
Moreover, $\mathcal{Q}(\X^+_{a,b})_{c,d}$ is generated from the top composition factor $\X^+_{a,b}$.
\end{prop}
\begin{prop}
\label{Q-thick}
Let $(a,b,c,d)$ be an arbitrary element in
\begin{align*}
\{(a,b,c,d)\}=&\bigl\{(r^\vee,s,r,s),(r^\vee,s,r^\vee,s^\vee),(r,s^\vee,r,s),(r,s^\vee,r^\vee,s^\vee)\}.
\end{align*}
Then the socle series of $\mathcal{Q}(\X^-_{a,b})_{c,d}$ is given by
\begin{align*}
\hspace{-3mm}{\rm S}_1=\X^-_{a,b},\quad 
{\rm S}_2/{\rm S}_1=2\X^+_{c,d},\quad 
\mathcal{Q}(\X^-_{a,b})_{c,d}/{\rm S}_2=\X^-_{a,b}.
\end{align*}
Moreover, $\mathcal{Q}(\X^-_{a,b})_{c,d}$ is generated from the top composition factor $\X^-_{a,b}$.
\end{prop}

\subsection{${\rm Ext}^1$-groups between simple modules}
\label{W-mod}
Fix integers $r,s\in \Z^2_{\geq 1}$ such that $1\leq r< p_+$ and $1\leq s< p_-$. In this subsection, we determine the ${\rm Ext}^1$-groups between all simple modules.

\begin{dfn}
\label{kaisei}
Fix an arbitrary pair $(a,b,c,d)$ in
\begin{align*}
\{(a,b,c,d)\}=&\bigl\{(r^\vee,s,r,s),(r^\vee,s,r^\vee,s^\vee),(r,s^\vee,r,s),(r,s^\vee,r^\vee,s^\vee),\\
&(r^\vee,p_-,r,p_-),(p_+,s^\vee,p_+,s)\}.
\end{align*}
\begin{enumerate}
\item For $\mathcal{Q}(\X^-_{a,b})_{c,d}$, let $\{v_+,v_-\}$ be a basis of the lowest weight space of the submodule $\X^-_{a,b}\subset \mathcal{Q}(\X^-_{a,b})_{c,d}$ such that
\begin{align*}
&W^\pm[0]v_\pm=0,
&W^\pm[0]v_\mp\in \C^\times v_\pm,
\end{align*}
and let $u_{\pm}$ be the vectors in $\mathcal{Q}(\X^-_{a,b})_{c,d}[\Delta^+_{c,d;0}]$ satisfying $v_{\pm}\in U(\mathcal{L}).u_\pm$. 
Then we define
\begin{align*}
&\mathcal{E}^+(\X^+_{c,d})_{a,b}:=\W_{p_+,p_-}.u_{+},
&\mathcal{E}^-(\X^+_{c,d})_{a,b}:=\W_{p_+,p_-}.u_{-},
\end{align*}
which give different extension classes in ${\rm Ext}^1(\X^+_{c,d},\X^-_{a,b})\setminus\{0\}$. 
\item We define
\begin{align*}
&\mathcal{E}^\pm(\X^-_{a,b})_{c,d}:=\mathcal{Q}(\X^-_{a,b})_{c,d}/\mathcal{E}^\mp(\X^+_{c,d})_{a,b}.
\end{align*}
\item As the quotient of $\mathcal{Q}(\X^+_{c,d})_{a,b}$, we define the indecomposable module $\widetilde{\mathcal{E}}(\X^+_{c,d})_{a,b}$ satisfying the exact sequence
\begin{align*}
0\rightarrow \X^-_{a,b}\oplus \X^-_{a,b}\rightarrow \widetilde{\mathcal{E}}(\X^+_{c,d})_{a,b}\rightarrow \X^+_{c,d}\rightarrow 0.
\end{align*}
\item We define $\widetilde{\mathcal{E}}(\X^-_{a,b})_{c,d}:=\mathcal{Q}(\X^-_{a,b})_{c,d}/\X^-_{a,b}$, which satisfies the exact sequence
\begin{align*}
0\rightarrow \X^+_{c,d}\oplus \X^+_{c,d}\rightarrow \widetilde{\mathcal{E}}(\X^-_{a,b})_{c,d}\rightarrow \X^-_{a,b}\rightarrow 0.
\end{align*}
\end{enumerate}
\end{dfn}

\begin{lem}
\label{lemE3}
Let $n\geq 1$. Let $(a,b,c,d)$ be an arbitrary pair in
\begin{align*}
\{(a,b,c,d)\}=&\bigl\{(r,s,r^\vee,s),(r^\vee,s^\vee,r^\vee,s),(r,s,r,s^\vee),(r^\vee,s^\vee,r,s^\vee),\\
&(r,p_-,r^\vee,p_-),(p_+,s,p_+,s^\vee)\}.
\end{align*}
Any extension 
$
[E]\in {\rm Ext}^1(\X^+_{a,b},n\X^-_{c,d})
$
is trivial if $E$ satisfies
\begin{align}
\label{airline0428}
{\rm Hom}_{U(\mathcal{L})}(L(\Delta^+_{a,b;0}),E)\neq 0.
\end{align}
\begin{proof}
We only prove in the case of $(a,b,c,d)=(r,s,r^\vee,s)$ and $n=1$. The other cases can be proved in the same way.

Fix a non-trivial extension $[E]\in {\rm Ext}(\X^+_{r,p_-},\X^-_{r^\vee,p_-})\setminus\{0\}$.
Let $u$ be the lowest weight vector in $E[\Delta^+_{r,s;0}]$.
Assume (\ref{airline0428}). Then we have
\begin{align}
\label{assumptE}
S_{r,s^\vee+p_-}u=0.
\end{align}
Let $\{v_+,v_-\}$ be a basis of the lowest weight space of the submodule $\X^-_{r^\vee,s}\subset E$ such that
\begin{align*}
&W^\pm[0]v_\pm=0,
&W^\pm[0]v_\mp\in \C^\times v_\pm.
\end{align*}
Let $v^*_{+}$ and $v^*_-$ be $L_0$-homogeneous vectors of $E^*$ such that $\langle v^*_{\pm},v_{\pm}\rangle\neq 0$, and $L_kv^*_{\pm}=0$ for $k\geq 1$. 
Assume that for any $W=W^\pm,W^0$
\begin{align*}
W[k]v^*_{\pm}=0,\ \ k\geq 1.
\end{align*}
Then the vector space $\C v^*_++\C v^*_-$ becomes a $A(\W_{p_+,p_-})$-module and this vector space is isomorphic to the lowest weight space of $\X^-_{r^\vee,s}$ as a $A(\W_{p_+,p_-})$-module. Thus $E^*$ has the submodule $\W_{p_+,p_-}.(\C v^*_++\C v^*_-)\simeq \X^-_{r^\vee,s}$ and thus $E^*=\X^+_{r,s}\oplus \X^-_{r^\vee,s}$. But this contradicts the assumption that $E$ is non-trivial.
Therefore we have
\begin{align*}
\langle v^*_{\pm},{Y}_{E}(W^{\bullet};z)u\rangle\neq 0,
\end{align*}
where $W^{\bullet}$ is one of $W^+$, $W^0$ or $W^-$. On the other hand, using Proposition \ref{FusionInt}, we have
\begin{align*}
&\langle v^*_{\pm},{Y}_{E}(W^{\bullet};z)S_{r,s^\vee+p_-}u\rangle\\
&=z^{-r(s^\vee+p_-)}\prod_{i=1}^{r}\prod_{j=1}^{s^\vee+p_-}(h_{4p_+-1,1}-h_{r+r+2p_+-2i+1,s^\vee+p_-+s-2j+1})\langle v^*_{\pm},{Y}_{E}(W^{\bullet};z)u\rangle\\
&\neq 0.
\end{align*}
But this contradicts (\ref{assumptE}).
\end{proof}
\end{lem}

As in the triplet $W$-algebras $\W_p$ \cite{AML,FF4,McRae,NT}, we have the following proposition.
\begin{prop}
\label{thinsim}
Let $(a,b,c,d)$ be $(r,p_-,r^\vee,p_-)$ or $(p_+,s,p_+,s^\vee)$.
In the thin block $C^{thin}_{a,b}$, we have 
\begin{align*}
{\rm Ext}^1(\X^+_{a,b},\X^-_{c,d})={\rm Ext}^1(\X^-_{c,d},\X^+_{a,b})=\C^2.
\end{align*}
The other extensions between the simple modules in this block are trivial.
\begin{proof}
We prove only 
\begin{align*}
&{\rm Ext}^1(\X^-_{r^\vee,p_-},\X^-_{r^\vee,p_-})=0,&{\rm Ext}^1(\X^+_{r,p_-},\X^-_{r^\vee,p_-})&=\C^2,\\
&{\rm Ext}^1(\X^+_{r,p_-},\X^+_{r,p_-})=0.
\end{align*}
The other ${\rm Ext}^1$-groups can be proved in the same way.

First let us prove ${\rm Ext}^1(\X^-_{r^\vee,p_-},\X^-_{r^\vee,p_-})=0$. 
Fix an arbitrary extension $[E_1]\in{\rm Ext}^1(\X^-_{r^\vee,p_-},\X^-_{r^\vee,p_-})$. By Proposition \ref{Vir0}, we see that $L_0$ acts semisimply on $E_1$. Let $\overline{E}_1$ be the lowest weight space of $E_1$. Note that $E_1$ is generated from $\overline{E}_1$. Let $\widetilde{E}_1$ be the $\W_{p_+,p_-}$-module induced from $\overline{E}_1$. Then from Proposition~\ref{simpleclass}, $\widetilde{E}_1=E_1$. By Proposition \ref{sl2}, we see that as the $A(\W_{p_+,p_-})$-module
\begin{align*}
\overline{E}_1\simeq\overline{\X^-_{r^\vee,p_-}}\oplus \overline{\X^-_{r^\vee,p_-}},
\end{align*}
where $\overline{\X^-_{r^\vee,p_-}}$ is the lowest weight space of $\X^-_{r^\vee,p_-}$.
Note that the $\W_{p_+,p_-}$-module induced from $\overline{\X^-_{r^\vee,p_-}}$ is isomorphic to $\X^-_{r^\vee,p_-}$.
Thus we have $\widetilde{E}_1\simeq {\X^-_{r^\vee,p_-}}\oplus {\X^-_{r^\vee,p_-}}$.

Next, we prove ${\rm Ext}^1(\X^+_{r,p_-},\X^-_{r^\vee,p_-})=\C^2$. 
Since ${\rm Ext}^1(\X^-_{r^\vee,p_-},\X^-_{r^\vee,p_-})=0$, it is sufficient to show that
\begin{align}
\label{eq:ext-etil}
{\rm Ext}^1(\widetilde{\mathcal{E}}(\X^+_{r,p_-})_{r^\vee,p_-},\X^-_{r^\vee,p_-})=0,
\end{align}
where $\widetilde{\mathcal{E}}(\X^+_{r,p_-})_{r^\vee,p_-}$ is defined in Definition \ref{kaisei}.
Assume that
\begin{align*}
{\rm Ext}^1(\widetilde{\mathcal{E}}(\X^+_{r,p_-})_{r^\vee,p_-},\X^-_{r^\vee,p_-})\neq 0
\end{align*}
and fix a non-trivial extension $[E_2]$ in this ${\rm Ext}^1$-group.
By ${\rm Ext}^1(\X^-_{r^\vee,s},\X^-_{r^\vee,s})=0$, we see that
\begin{align}
\label{socle3tu}
{\rm S}(E_2)=\X^-_{r^\vee,p_-}\oplus \X^-_{r^\vee,p_-}\oplus \X^-_{r^\vee,p_-}.
\end{align}
Let $u$ be the lowest weight vector of $E_2$. Consider the submodule 
\begin{align*}
I=\W_{p_+,p_-}.S_{r,p_-}u 
\end{align*}
of $E_2$. Let $\overline{I}$ be the lowest weight space of $I$. By Proposition \ref{sl2}, 
\begin{align*}
\overline{I}_1\simeq\overline{\X^-_{r^\vee,p_-}}\oplus \overline{\X^-_{r^\vee,p_-}},
\end{align*}
as $A(\W_{p_+,p_-})$-modules.
Then by Proposition~\ref{simpleclass}, we have
\begin{align}
\label{socle2tu}
\W_{p_+,p_-}.S_{r,p_-}u\simeq \X^-_{r^\vee,p_-}\oplus \X^-_{r^\vee,p_-}.
\end{align}
Thus by (\ref{socle3tu}) and (\ref{socle2tu}), we have 
\begin{align*}
[E_2/\W_{p_+,p_-}.S_{r,p_-}u]\in {\rm Ext}^1(\X^+_{r,p_-},\X^-_{r^\vee,p_-})\setminus\{0\},
\end{align*}
But this contradicts Lemma \ref{lemE3}. Therefore we obtain (\ref{eq:ext-etil}).

Finally let us prove ${\rm Ext}^1(\X^+_{r,p_-},\X^+_{r,p_-})=0$. 
Let $\overline{\X^+_{r,p_-}}$ be the lowest weight space of $\X^+_{r,p_-}$. 
Since 
\begin{align*}
&{\rm Ext}^1(\X^-_{r^\vee,p_-},\X^-_{r^\vee,p_-})=0,
&{\rm Ext}^1(\X^+_{r,p_-},\X^-_{r^\vee,p_-})=\C^2,
\end{align*}
we see that the induced $\W_{p_+,p_-}$-module from $\overline{\X^+_{r,s}}$ is isomorphic to $\widetilde{\mathcal{E}}(\X^+_{r,p_-})_{r^\vee,p_-}$. 
Let $[E_3]$ be an extension in ${\rm Ext}^1(\X^+_{r,p_-},\X^+_{r,p_-})$. By Proposition \ref{Vir0}, we see that $L_0$ acts semisimply on $E_3$. 
Let $\overline{E}_3$ be the lowest weight space of $E_3$. 
Note that $E_3$ is generated from $\overline{E}_3$.
Assume that $\overline{E}_3\ncong \overline{\X^+_{r,p_-}}\oplus \overline{\X^+_{r,p_-}}$ as a $A(\W_{p_+,p_-})$-module. Then, from the $\W_{p_+,p_-}$-action on $E_3$, we have a non-trivial non-logarithmic Virasoro intertwining operator of type
\begin{align*}
\begin{pmatrix}
   \ L(\Delta^+_{r,p_-;0})  \\
   L(h_{4p_+-1,1})\ \ L(\Delta^+_{r,p_-;0})
\end{pmatrix}
.
\end{align*}
But we have a contradiction by using Proposition \ref{VirasoroFusion}. 
Thus, as the $A(\W_{p_+,p_-})$-module, $\overline{E}_3\simeq \overline{\X^+_{r,p_-}}\oplus \overline{\X^+_{r,p_-}}$. 
Let $\widetilde{E}_3$ be the induced module from $\overline{E}_3$. Then we have
\begin{align*}
\widetilde{E}_3\simeq\widetilde{\mathcal{E}}(\X^+_{r,p_-})_{r^\vee,p_-}\oplus \widetilde{\mathcal{E}}(\X^+_{r,p_-})_{r^\vee,p_-}.
\end{align*}
Therefore we obtain
$
E_3\simeq {\X^+_{r,p_-}}\oplus {\X^+_{r,p_-}}.
$
\end{proof}
\end{prop}

\begin{prop}
\label{Ext}
In the thick block $C^{thick}_{r,s}$, we have 
\begin{align*}
{\rm Ext}^1(\mathcal{X}^{\pm},\mathcal{X}^{\mp})=\C^2,\ \ \ \ \ {\rm Ext}^1(L(h_{r,s}),\X^+)={\rm Ext}^1(\X^+,L(h_{r,s}))=\C,
\end{align*}
where $\mathcal{X}^+=\mathcal{X}^+_{r,s}$ or $\mathcal{X}^+_{r^\vee,s^\vee}$ and $\mathcal{X}^-=\mathcal{X}^-_{r^\vee,s}$ or $\mathcal{X}^-_{r,s^\vee}$. The other extensions between the simple modules in $C^{thick}_{r,s}$ are trivial.
\begin{proof}
We prove only ${\rm Ext}^1(L(h_{r,s}),\mathcal{X}^+_{r,s})=\C$.
The other ${\rm Ext}^1$-groups can be determined in a similar way as Proposition \ref{thinsim}.

Note that 
\begin{align*}
[\W_{p_+,p_-}.\ket{\alpha_{r,s}}]\in {\rm Ext}^1(L(h_{r,s}),\X^+_{r,s})\setminus\{0\},
\end{align*}
and, as the Virasoro module
\begin{align*}
\W_{p_+,p_-}.\ket{\alpha_{r,s}}=L(h_{r,s},\Delta^+_{r,s;0})\oplus \bigoplus_{n\geq 1}(2n+1)L(\Delta^+_{r,s;n}).
\end{align*} 
Fix a nonzero extension 
$
[E]\in {\rm Ext}^1(L(h_{r,s}),\X^+_{r,s}).
$
Let $t$ be the lowest weight vector of $E$. Assume that $S_{r^\vee,s^\vee}t\neq 0$. 
Then, as the Virasoro module
\begin{align*}
E=L(h_{r,s},\Delta^+_{r,s;0})\oplus \bigoplus_{n\geq 1}(2n+1)L(\Delta^+_{r,s;n}).
\end{align*} 
Since 
$
{\rm Ext}^1_{\mathcal{L}}(L(h_{r,s}),L(\Delta^+_{r,s;0}))=\C,
$
as the Baer sum of extensions obtained from the indecomposable modules $E$ and $\W_{p_+,p_-}.\ket{\alpha_{r,s}}$, we have an extension $[E']\in {\rm Ext}^1(L(h_{r,s}),\X^+_{r,s})$ such that $S_{r^\vee,s^\vee}t'=0$, where $t'$ is the lowest weight vector of $E'$. Thus, by Theorem $\ref{VirC}$, we have the following decomposition as the Virasoro module
\begin{align*}
E'=L(h_{r,s})\oplus \bigoplus_{n\geq 0}(2n+1)L(\Delta^+_{r,s;n}).
\end{align*} 
Assume $[E']\neq 0$.
Then, from the $\W_{p_+,p_-}$-module action on $E'$, we must have a non-trivial Virasoro intertwining operator of type
\begin{align*}
\begin{pmatrix}
   \ L(\Delta^+_{r,s;n})  \\
   L(h_{4p_+-1,1})\ \ L(h_{r,s})
\end{pmatrix}
\end{align*}
for some $n\geq 0$. But, by using Proposition \ref{VirasoroFusion}, we have a contradiction. In case $S_{r^\vee,s^\vee}t=0$, we see that $[E]=0$ as shown above.
\end{proof}
\end{prop}

\subsection{Projective modules in the thin blocks}
\label{54}
Fix arbitrary two thin blocks $C^{thin}_{r,p_-}$, $C^{thin}_{p_+,s}$ $(1\leq r< p_+,\ 1\leq s< p_-)$. In this subsection, we will show the logarithmic modules $\mathcal{Q}(\X^+_{r,p_-})_{r^\vee,p_-},\mathcal{Q}(\X^-_{r^\vee,p_-})_{r,p_-}\in C^{thin}_{r,p_-}$ and $\mathcal{Q}(\X^+_{p_+,s})_{p_+,s^\vee},\mathcal{Q}(\X^-_{p_+,s^\vee})_{p_+,s}\in C^{thin}_{p_+,s}$ are projective.

\begin{lem}
\label{bt3lem}
Let $(a,b,c,d)$ be $(r^\vee,p_-,r,p_-)$ or $(p_+,s^\vee,p_+,s)$. Fix an arbitrary indecomposable module $E$ corresponding to a class in ${\rm Ext}^1(\X^-_{a,b},\X^+_{c,d})$. For the surjection $\pi:E\twoheadrightarrow \X^-_{a,b}$, let $(v^-,v^+)\in E[\Delta^-_{a,b;0}]^2$ be $L_0$-homogeneous vectors of $E$ satisfying $\pi (v^\pm)=v^{(0)}_{\pm \frac{1}{2}}$, where $v^{(0)}_{\pm \frac{1}{2}}$ are the Virasoro singular vectors defined in Definition \ref{dfn0830}.
Then, for $\epsilon=\pm$, the Virasoro module $U(\mathcal{L}).v^\epsilon$ has a quotient isomorphic to $
L(\Delta^-_{a,b;0},\Delta^+_{c,d;1}).
$
\end{lem}
\begin{proof}
We only prove the case of $(a,b,c,d)=(r^\vee,p_-,r,p_-)$. The second case can be proved in the same way.

Recall that $\widetilde{\mathcal{E}}(\X^-_{r^\vee,p_-})_{r,p_-}$ is defined by the quotient $\mathcal{Q}(\mathcal{X}^-_{r^\vee,p_-})_{r,p_-}/\X^-_{r^\vee,p_-}$. Then 
\begin{align}
\label{eq:eex}
[\widetilde{\mathcal{E}}(\X^-_{r^\vee,p_-})_{r,p_-}]\in {\rm Ext}^1(\mathcal{E}^\pm(\X^-_{r^\vee,p_-})_{r,p_-},\X^-_{r^\vee,p_-}).
\end{align}
Viewing $\widetilde{\mathcal{E}}(\X^-_{r^\vee,p_-})_{r,p_-}$ as the Virasoro module, it follows from the Virasoro module structure of $\mathcal{Q}(\mathcal{X}^-_{r^\vee,p_-})_{r,p_-}$ that $\widetilde{\mathcal{E}}(\X^-_{r^\vee,p_-})_{r,p_-}$ has a quotient isomorphic to $2L(\Delta^-_{r^\vee,p_-;0},\Delta^+_{r,p_-;1})$. Then, noting (\ref{eq:eex}), we see that the assertion of the lemma holds for $\mathcal{E}^\pm(\X^-_{r^\vee,p_-})_{r,p_-}$.
Since from Proposition~\ref{thinsim}, $\{[\mathcal{E}^\epsilon(\X^-_{r^\vee,p_-})_{r,p_-}]|\ \epsilon=\pm\}$ gives a basis of ${\rm Ext}^1(\X^-_{r^\vee,p_-},\X^+_{r,p_-})$, it follows from the structure of $\mathcal{E}^\pm(\X^-_{r^\vee,p_-})_{r,p_-}$ that the assertion of the lemma holds.
\end{proof}
From this subsection, we will frequently use the following notation.
\begin{dfn}
Let $M\in \mathcal{L}_{c_{p_+,p_-}}\mathchar`-{\rm mod}$ and fix an arbitrary $u\in M$. When two vectors $v_1,v_2\in M$ coincide except for an element in $U(\mathcal{L}).u$, denote it by 
\begin{align*}
v_1\equiv v_2\ \ {\rm mod}\ U(\mathcal{L}).u.
\end{align*}
\end{dfn}
\begin{prop}
\label{Ext2+0}
We have
\begin{align*}
&{\rm Ext}^1(\mathcal{Q}(\X^+_{r,p_-})_{r^\vee,p_-},\X^+_{r,p_-})=0,
&{\rm Ext}^1(\mathcal{Q}(\X^+_{p_+,s})_{p_+,s^\vee},\X^+_{p_+,s})=0.
\end{align*}
\end{prop}
\begin{proof}
We only prove the first equality. The second equality can be proved in the same way.

Assume that ${\rm Ext}^1(\mathcal{Q}(\X^+_{r,p_-})_{r^\vee,p_-},\X^+_{r,p_-})\neq 0$ and fix a non-trivial extension $[E]$ in this ${\rm Ext}^1$-group. 
Let $\iota_1$ and $\iota_2$ be injections from $\X^+_{r,p_-}$ to $E$ such that
\begin{align*}
\iota_1(\X^+_{r,p_-})\oplus \iota_2(\X^+_{r,p_-})={\rm S}(E).
\end{align*}
For the surjection $\pi :E\rightarrow \X^+_{r,p_-}$, we fix $L_0$-homogeneous vectors $\widetilde{u}^{(0)}_0,\widetilde{u}^{(1)}_{\pm 1}, \widetilde{u}^{(1)}_{0}\in E$ such that
\begin{align*}
&\pi(\widetilde{u}^{(0)}_0)=u^{(0)}_0,
&\pi(\tilde{u}^{(1)}_0)=u^{(1)}_0,&& &\pi(\tilde{u}^{(1)}_{\pm 1})=u^{(1)}_{\pm 1},
\end{align*}
where ${u}^{(0)}_0,{u}^{(1)}_{\pm 1}$ and ${u}^{(1)}_{0}$ are the Virasoro singular vectors defined in Definition~\ref{dfn0830}.
Set $v^\pm=W^\pm[0]W^\mp[0]S_{r,p_-}\widetilde{u}^{(0)}_0$. Note that $\{v^-,W^+[0]v^-\}$ and $\{W^-[0]v^+,v^+\}$ correspond to bases of the lowest weight spaces of $\X^-_{r^\vee,p_-}\oplus \X^-_{r^\vee,p_-}$ (see also Figure~\ref{fig61}).

Consider the Virasoro module
\begin{align*}
M=U(\mathcal{L}).\widetilde{u}^{(1)}_0+U(\mathcal{L}).v^++U(\mathcal{L}).v^-.
\end{align*}
Let $\pi_Q$ be the surjection from $E$ to $\mathcal{Q}(\X^+_{r,p_-})_{r^\vee,p_-}$.
Note that by definition, we can choose the vectors $\widetilde{u}^{(1)}_0$ and $v^\pm$ to satisfy $\pi_Q(\widetilde{u}^{(1)}_0),\pi_Q(v^+)\in F_{r^\vee,p_-;-1}$ and $\pi_Q(v^-)\in F_{r,p_-}$ (see Figure \ref{figQQ}). 
Then, from the Virasoro structure of $\mathcal{Q}(\X^+_{r,p_-})_{r^\vee,p_-}$, we have
\begin{align*}
(L_0-\Delta^+_{r,p_-;1})\pi_Q(\widetilde{u}^{(1)}_0)\neq 0.
\end{align*}
Then, we see that the Virasoro module $M$ has an indecomposable quotient isomorphic to 
\begin{align*}
\widehat{P}^{\downarrow }(\tau),\ \ \  \tau=(\alpha_{r^\vee,p_-;1},\alpha_{r,p_-;2},\alpha_{r^\vee,p_-;3}), 
\end{align*}
and that the quotient module
\begin{equation*}
M/(U(\mathcal{L}).v^++U(\mathcal{L}).v^-)
\end{equation*}
does not contain $L(\Delta^-_{r^\vee,p_-;0})$ as a composition factor.
In particular, noting the Hom-structure of $\widehat{P}^{\downarrow }(\tau)$ described in Proposition~\ref{bt3}, we have
\begin{align}
\label{bt1}
 \sigma(S_{r^\vee,2p_-})\widetilde{u}^{(1)}_0\in \C^\times v^++\C^\times v^-+U(\mathcal{L}).\iota_1({u}^{(0)}_0)+U(\mathcal{L}).\iota_2({u}^{(0)}_0).
\end{align} 
By the structure of $\widehat{P}(\tau)$, we have 
\begin{equation}
\label{ntt0004230}
S_{r^\vee,2p_-}\sigma(S_{r^\vee,2p_-})\widetilde{u}^{(1)}_0\not\equiv 0\ {\rm mod}\ U(\mathcal{L}).\iota_1({u}^{(0)}_0)+U(\mathcal{L}).\iota_2({u}^{(0)}_0).
\end{equation}
and 
\begin{equation}
\label{ntt0004231}
\begin{split}
S_{r^\vee,2p_-}\sigma(S_{r^\vee,2p_-})\widetilde{u}^{(1)}_0 \in &\C\iota_1({u}^{(1)}_0)+\C\iota_2({u}^{(1)}_0)+\sum_{i=1,2}\sum_{\epsilon=\pm}\C\iota_i({u}^{(1)}_{\epsilon 1})\\
&+U(\mathcal{L}).\iota_1({u}^{(0)}_0)+U(\mathcal{L}).\iota_2({u}^{(0)}_0).
\end{split}
\end{equation}
Assume that the coefficient of $\iota_i({u}^{(1)}_{\pm 1})$ in $S_{r^\vee,2p_-}\sigma(S_{r^\vee,2p_-})\widetilde{u}^{(1)}_0$ is nonzero. Then, from the structure of $\widehat{P}(\tau)$, we see that the coefficient of $\iota_i({u}^{(1)}_{\pm 1})$ in $(L_0-\Delta^+_{r,p_-;1})\widetilde{u}^{(1)}_0$ is also nonzero. Thus, noting Proposition \ref{sl2action2} and multiplying $(L_0-\Delta^+_{r,p_-;1})\widetilde{u}^{(1)}_0$ by $(W^\mp[0])^2$, we have $\iota_i({u}^{(1)}_{\pm 1})\in U(\mathcal{L}).\iota_1({u}^{(0)}_0)+U(\mathcal{L}).\iota_2({u}^{(0)}_0)$. But this is a contradiction. 
Therefore, from (\ref{ntt0004231}), we have
\begin{equation}
\label{ntt0004232}
S_{r^\vee,2p_-}\sigma(S_{r^\vee,2p_-})\widetilde{u}^{(1)}_0\equiv k_1\iota_1({u}^{(1)}_0)+k_2\iota_2({u}^{(1)}_0)\ {\rm mod}\ \sum_{i=1}^2U(\mathcal{L}).\iota_i({u}^{(0)}_0),
\end{equation}
where $(k_1,k_2)(\neq (0,0))$ are some constants. Since ${\rm Ext}^1(\X^+_{r,p_-},\X^+_{r,p_-})=0$, we have
\begin{align}
\label{al0516}
\W_{p_+,p_-}.(k_1\iota_1({u}^{(1)}_0)+k_2\iota_2({u}^{(1)}_0))\simeq \X^+_{r,p_-}. 
\end{align}
Let $E'$ be the quotient module defined by
\begin{align*}
E'=\frac{E}{\W_{p_+,p_-}.(k_1\iota_1({u}^{(1)}_0)+k_2\iota_2({u}^{(1)}_0))}.
\end{align*}
Then by the assumption that $[E]\neq 0$ and ${\rm Ext}^1(\mathcal{X}^+_{r,p_-},\mathcal{X}^+_{r,p_-})=0$, we have 
\begin{align*}
[E']\in {\rm Ext}^1(\widetilde{\mathcal{E}}(\mathcal{X}^+_{r,p_-})_{r^\vee,p_-},\X^+_{r,p_-})\setminus\{0\},
\end{align*}
where $\widetilde{\mathcal{E}}(\mathcal{X}^+_{r,p_-})_{r^\vee,p_-}$ is defined in Definition~\ref{kaisei}.
Let $\phi$ be the surjection from $E$ to $E'$. 
From (\ref{ntt0004230}), (\ref{ntt0004232}) and (\ref{al0516}), we have
\begin{align}
\label{bt2}
S_{r^\vee,2p_-}\sigma(S_{r^\vee,2p_-})\phi(\widetilde{u}^{(1)}_0)\in \sum_{i=1,2}U(\mathcal{L}).\phi\circ\iota_i({u}^{(0)}_0).
\end{align}
By (\ref{bt2}), we see that $\phi(M)$ does not have $\widehat{P}(\tau)$ as any subquotient.
Then, from Proposition~\ref{bt3}, (\ref{bt1}) and (\ref{ntt0004232}), we have
\begin{align}
\label{bt4}
S_{r^\vee,2p_-}\phi(v^\pm)\in \sum_{i=1,2}U(\mathcal{L}).\phi\circ\iota_i({u}^{(0)}_0).
\end{align}
By Lemma \ref{bt3lem} and (\ref{bt4}), $E'$ has $\X^-_{r^\vee,p_-}\oplus \X^-_{r^\vee,p_-}$ as submodules. Thus $E'$ has a indecomposable quotient representing a class in ${\rm Ext}^1(\X^+_{r,p_-},\X^+_{r,p_-})$. But this contradicts Proposition \ref{thinsim}.
\end{proof}

\begin{figure}[htbp]
  \centering
\begin{tikzpicture}[scale=0.85]
\node[inner sep=0.7pt] (b1) at (0, -3) {${\pi_Q(v^-)}$};
\node[inner sep=0.7pt] (b2) at (3, -3) {${\pi_Q(v^+)}$};
\node[inner sep=0.7pt] (aa1) at (0, -6) {${\pi_Q(\widetilde{u}^{(1)}_0)}$};
\node[inner sep=0.9pt] (aa2) at (3, -6) {${\bullet}$};
\node[inner sep=0.7pt] (dd1) at (0, -1.5) {$F_{r^\vee,p_-;-1}$};
\node[inner sep=0.7pt] (dd2) at (3, -1.5) {$F_{r,p_-}$};
\node[inner sep=0.7pt] (cc1) at (0, -7) {};
\node[inner sep=0.7pt] (cc2) at (3, -7) {};

\draw[arrows = {-Stealth[scale=0.9]},very thick] (aa1) to (b1);
\draw[arrows = {-Stealth[scale=0.9]},very thick] (aa1) to (b2);
\draw[arrows = {-Stealth[scale=0.9]}] (b1) to (aa2);
\draw[arrows = {-Stealth[scale=0.9]}] (b2) to (aa2);

\draw[arrows = {-Stealth[scale=0.9]},dashed] (aa1) to (aa2);
\draw[dotted,very thick] (cc1) to (aa1);
\draw[dotted,very thick] (cc2) to (aa2);
\draw[dotted,very thick] (dd1) to (b1);
\draw[dotted,very thick] (dd2) to (b2);
\end{tikzpicture}
\caption{\label{figQQ}The schematic diagram of the logarithmic Virasoro modules $\pi_Q(M)=U(\mathcal{L}).\pi_Q(\widetilde{u}^{(1)}_0)$ defined in the proof of Proposition \ref{Ext2+0}. The meaning of the arrows is the same as in Figures~\ref{51figk},~\ref{51figkk}~and~\ref{51figkkk}. 
The black circle represents the lowest weight vectors of the simple module $L(\Delta_{r,p_-;1})\subset {\rm S}(F_{r,p_-})$.
The dotted vertical lines represent that the vectors are contained in the Fock modules $F_{r^\vee,p_-;-1}\oplus F_{r,p_-}$.}
\mbox{}\\
\end{figure}

Recall the indecomposable modules in Definition~\ref{kaisei}.
\begin{prop}
\label{Ext2-+0}
We have
\begin{align*}
&{\rm Ext}^1(\mathcal{Q}(\X^-_{r^\vee,p_-})_{r,p_-},\X^+_{r,p_-})=0,
&{\rm Ext}^1(\mathcal{Q}(\X^-_{p_+,s^\vee})_{p_+,s},\X^+_{p_+,s})=0.
\end{align*}
\begin{proof}
We only prove the first equality. The second equality can be proved in the same way.

By Proposition \ref{Ext2+0}, we have
\begin{align}
{\rm Ext}^1(\mathcal{E}^+(\X^+_{r,p_-})_{r^\vee,p_-},\X^+_{r,p_-})=0.
\label{E+}
\end{align}
From Definition~\ref{kaisei}, we have the following exact sequence
\begin{align*}
0\rightarrow \mathcal{E}^+(\X^+_{r,p_-})_{r^\vee,p_-}\xrightarrow{} \mathcal{Q}(\X^-_{r^\vee,p_-})_{r,p_-}\rightarrow \mathcal{E}^{-}(\X^-_{r^\vee,p_-})_{r,p_-}\rightarrow 0.
\end{align*}
By this exact sequence and (\ref{E+}), we have the following exact sequence
\begin{align*}
0\rightarrow \C\rightarrow {\rm Ext}^1(\mathcal{E}^{-}(\X^-_{r^\vee,p_-})_{r,p_-},\X^+_{r,p_-})\xrightarrow{} {\rm Ext}^1(\mathcal{Q}(\X^-_{r^\vee,p_-})_{r,p_-},\X^+_{r,p_-})\rightarrow 0.
\end{align*}
By Proposition \ref{thinsim}, we have ${\rm Ext}^1(\mathcal{E}^{-}(\X^-_{r^\vee,p_-})_{r,p_-},\X^+_{r,p_-})\simeq\C$. Therefore we obtain 
$
{\rm Ext}^1(\mathcal{Q}(\X^-_{r^\vee,p_-})_{r,p_-},\X^+_{r,p_-})=0.
$
\end{proof}
\end{prop}

\begin{lem}
\label{lemL0}
We have
\begin{align*}
{\rm Ext}^1(\mathcal{E}^\pm(\X^-_{r^\vee,p_-})_{r,p_-},\X^-_{r^\vee,p_-})={\rm Ext}^1(\mathcal{E}^\pm(\X^-_{p_+,s^\vee})_{p_+,s},\X^-_{p_+,s^\vee})=0.
\end{align*}
\begin{proof}
We only prove ${\rm Ext}^1(\mathcal{E}^+(\X^-_{r^\vee,p_-})_{r,p_-},\X^-_{r^\vee,p_-})=0$. The other equalities can be proved in the same way.

Assume that 
$
{\rm Ext}^1(\mathcal{E}^+(\X^-_{r^\vee,p_-})_{r,p_-},\X^-_{r^\vee,p_-})\neq 0.
$
Fix a non-trivial extension
\begin{align*}
0\rightarrow \X^-_{r^\vee,p_-}\xrightarrow{\iota} E\xrightarrow{} \mathcal{E}^+(\X^-_{r^\vee,p_-})_{r,p_-}\rightarrow 0.
\end{align*}
For the surjection $\pi: E\rightarrow \X^-_{r^\vee,p_-}$, let $\tilde{v}_\pm$ be any $L_0$ homogeneous vectors of $E$ such that $\pi(\tilde{v}_\pm)=v^{(0)}_{\pm\frac{1}{2}}$, where $v^{(0)}_{\frac{1}{2}}$ and $v^{(0)}_{-\frac{1}{2}}$ are the Virasoro singular vectors defined in Definition \ref{dfn0830}. 
Since ${\rm Ext}^1(\X^-_{r^\vee,p_-},\X^-_{r^\vee,p_-})=0$, we see that $E$ has an indecomposable submodule representing a class in ${\rm Ext}^1(\X^+_{r,p_-},\X^-_{r^\vee,p_-})$. Then, noting Lemma \ref{lemE3} and the definition of $\mathcal{E}^+(\X^-_{r^\vee,p_-})_{r,p_-}$ (see Definition~\ref{kaisei}), we see that the Virasoro module $U(\mathcal{L}). \tilde{v}_+$ a quotient isomorphic to $P(\tau)$, where $\tau=(\alpha_{r,p_-},\alpha_{r^\vee,p_-;1},\alpha_{r,p_-;2})$ (for the definition of the logarithmic modules $P(\tau)$, see Definition \ref{sankakuthickdf}). Note that, from the structure of $\mathcal{E}^+(\X^-_{r^\vee,p_-})_{r,p_-}$, $L_n\tilde{v}_-=0$ for $n\geq 1$. Thus, by Proposition \ref{sankakuthick}, we have
\begin{align}
\label{L02023}
&(L_0-\Delta^-_{r^\vee,p_-;0})\tilde{v}_+ =k_+\iota(v^{(0)}_{\frac{1}{2}})+k_-\iota(v^{(0)}_{-\frac{1}{2}}),\\
&(L_0-\Delta^-_{r^\vee,p_-;0})\tilde{v}_-=0,
\label{L020232}
\end{align}
where $(k_+,k_-)\neq (0,0)$. Asuume $k_+\neq 0$. Then, multiplying both sides of (\ref{L02023}) by $W^-[0]$, we have
\begin{align*}
(L_0-\Delta^-_{r^\vee,p_-;0})W^-[0]\tilde{v}_+ \in \C^\times \iota(v^{(0)}_{-\frac{1}{2}}).
\end{align*}
But this contradicts (\ref{L020232}). Next assume $k_-\neq 0$. 
Then, multiplying both sides of (\ref{L02023}) by $W^+[0]$, we have
\begin{align*}
(L_0-\Delta^-_{r^\vee,p_-;0})W^+[0]\tilde{v}_+ \in \C^\times \iota(v^{(0)}_{\frac{1}{2}}).
\end{align*}
On the other hand, by the definition of $\tilde{v}_+$, we have $(L_0-\Delta^-_{r^\vee,p_-;0})W^+[0]\tilde{v}_+=0$. Thus we have a contradiction.
\end{proof}
\end{lem}

\begin{prop}
\label{Ext2-0} 
We have
\begin{align*}
&{\rm Ext}^1(\mathcal{Q}(\X^-_{r^\vee,p_-})_{r,p_-},\X^-_{r^\vee,p_-})=0,
&{\rm Ext}^1(\mathcal{Q}(\X^-_{p_+,s^\vee})_{p_+,s},\X^-_{p_+,s^\vee})=0
\end{align*}
\begin{proof}
We only prove the first equality. The second equality can be proved in the same way.

Assume that 
$
{\rm Ext}^1(\mathcal{Q}(\X^-_{r^\vee,p_-})_{r,p_-},\X^-_{r^\vee,p_-})\neq 0.
$
Fix a non-trivial extension
\begin{align*}
0\rightarrow \X^-_{r^\vee,p_-}\xrightarrow{\iota} E\xrightarrow{\phi} \mathcal{Q}(\X^-_{r^\vee,p_-})_{r,p_-}\rightarrow 0.
\end{align*}
Since ${\rm Ext}^1(\X^-_{r^\vee,p_-},\X^-_{r^\vee,p_-})=0$, $E$ has a submodule representing a class in
\begin{align*}
{\rm Ext}^1(\widetilde{\mathcal{E}}(\X^-_{r^\vee,p_-})^*_{r,p_-},\X^-_{r^\vee,p_-})\setminus\{0\}.
\end{align*}
Then, from Propositions \ref{thinsim} and Lemma \ref{lemL0}, we see that the following sequence of submodules holds
\begin{align}
\label{2023Hold2}
\iota(\X^-_{r^\vee,p_-})\subset \widetilde{\mathcal{E}}(\X^+_{r,p_-})_{r^\vee,p_-}\subset E.
\end{align}
Let $\{v_+,v_-\}$ be a basis of the lowest weight space of $\X^-_{r^\vee,p_-}$ such that
\begin{align*}
&W^\pm[0]v_\pm=0,
&W^\pm[0]v_\mp\in\C^\times v_{\pm}.
\end{align*}
Let $\{v^0_+,v^0_-\}$ be a basis of the lowest weight space of the submodule $\X^-_{r^\vee,p_-}\subset E$ such that $\phi(v^0_{\pm})\neq 0$ and 
\begin{align*}
&W^\pm[0]v^0_\pm=0,
&W^\pm[0]v^0_\mp\in\C^\times v^0_{\pm}.
\end{align*}
For the surjection $\pi : E\rightarrow \X^-_{r^\vee,p_-}$, we fix $L_0$-homogeneous vectors $\tilde{v}_-,\tilde{v}_+\in E$ such that $\pi(\tilde{v}_{\pm})=v_{\pm}$.
Note that the Virasoro module $U(\mathcal{L}).\tilde{v}_\pm$ has a quotient isomorphic to $P(\tau)$, where $\tau=(\alpha_{r,p_-},\alpha_{r^\vee,p_-;1},\alpha_{r,p_-;2})$. 
Thus, by Proposition \ref{thinsim} and (\ref{2023Hold2}), we have
\begin{equation}
\label{2022KKK0518}
\begin{split}
&S_{r,p_-}\sigma(S_{r,p_-})\tilde{v}_- \in \C^\times\iota(v_+)+\C\iota(v_-)+\C^\times v^0_-,\\
&S_{r,p_-}\sigma(S_{r,p_-})\tilde{v}_+ \in \C\iota(v_+)+\C^\times\iota(v_-)+\C^\times v^0_-.
\end{split}
\end{equation}
Then, by Proposition \ref{sankakuthick} and (\ref{2022KKK0518}), we see that one of the followings holds
\begin{align*}
&(L_0-\Delta^-_{r^\vee,p_-;0})\tilde{v}_-=k_+\iota(v_+)+k_-\iota(v_-)+\C^\times v^0_-,\ \ \ k_+\neq 0,\\
&(L_0-\Delta^-_{r^\vee,p_-;0})\tilde{v}_+=l_-\iota(v_-)+l_+\iota(v_+)+\C^\times v^0_+,\ \ \ l_-\neq 0.
\end{align*}
Assume that the first equality is true. Multiplying the first equation by $W^-[0]$, we have
\begin{align*}
(L_0-\Delta^-_{r^\vee,p_-;0})W^-[0]\tilde{v}_-=k_+\iota(v_-).
\end{align*}
By the definition of $\tilde{v}_-$, the left hand side becomes zero. But this is a contradiction. In a similar way, assuming the second equality, we can show the contradiction.
\end{proof}
\end{prop}

\begin{prop}
\label{Ext2+-0}
We have
\begin{align*}
{\rm Ext}^1(\mathcal{Q}(\X^+_{r,p_-})_{r^\vee,p_-},\X^-_{r^\vee,p_-})={\rm Ext}^1(\mathcal{Q}(\X^+_{p_+,s})_{p_+,s^\vee},\X^-_{p_+,s^\vee})=0.
\end{align*}
\begin{proof}
We only prove the first equality. The second equality can be proved in the same way.

Note that $\mathcal{Q}(\X^+_{r,p_-})_{r^\vee,p_-}$ satisfies the following exact sequence
\begin{align*}
0\rightarrow \mathcal{E}^{-}(\X^-_{r^\vee,p_-})_{r,p_-}\xrightarrow{} \mathcal{Q}(\X^+_{r,p_-})_{r^\vee,p_-}\rightarrow \mathcal{E}^{+}(\X^+_{r,p_-})_{r^\vee,p_-}\rightarrow 0.
\end{align*}
By this exact sequence and Proposition \ref{Ext2-0}, we have the following exact sequence
\begin{align*}
0\rightarrow \C\rightarrow {\rm Ext}^1(\mathcal{E}^{+}(\X^+_{r,p_-})_{r^\vee,p_-},\X^-_{r^\vee,p_-})\xrightarrow{} {\rm Ext}^1(\mathcal{Q}(\X^+_{r,p_-})_{r^\vee,p_-},\X^-_{r^\vee,p_-})\rightarrow 0.
\end{align*}
By Proposition \ref{thinsim} we have ${\rm Ext}^1(\mathcal{E}^{+}(\X^+_{r,p_-})_{r^\vee,p_-},\X^-_{r^\vee,p_-})\simeq\C$. Therefore we obtain ${\rm Ext}^1(\mathcal{Q}(\X^+_{r,p_-})_{r^\vee,p_-},\X^-_{r^\vee,p_-})=0$.

\end{proof}
\end{prop}

By Propositions \ref{Ext2+0}, \ref{Ext2-+0}, \ref{Ext2-0} and \ref{Ext2+-0}, we obtain the following theorem.
\begin{thm}
\mbox{}
\begin{enumerate}
\item $\mathcal{Q}(\X^+_{r,p_-})_{r^\vee,p_-}$ and $\mathcal{Q}(\X^-_{r^\vee,p_-})_{r,p_-}$ are the projective covers of $\X^+_{r,p_-}$ and $\X^-_{r^\vee,p_-}$, respectively.
\item $\mathcal{Q}(\X^+_{p_+,s})_{p_+,s^\vee}$ and $\mathcal{Q}(\X^-_{p_+,s^\vee})_{p_+,s}$ are the projective covers of $\X^+_{p_+,s}$ and $\X^-_{p_+,s^\vee}$, respectively.
\end{enumerate}
\end{thm}

\subsection{Projective modules in the thick blocks}
\label{53}
In this subsection, we fix an arbitrary thick block $C^{thick}_{r,s}$ and compute ${\rm Ext}^1$ groups between certain indecomposable $\W_{p_+,p_-}$-modules and the simple modules in this block. Based on these ${\rm Ext}^1$-groups, we prove that the logarithmic modules $\mathcal{P}^\pm_{\bullet,\bullet}$ are projective $\W_{p_+,p_-}$-modules.

\begin{dfn}
\label{Kzenhan}
Let $(a,b)$ be $(r,s)$ or $(r^\vee,s^\vee)$. We denote by $\mathcal{K}_{a,b}$ an indecomposable module representing a class in the ${\rm Ext}^1$-group 
$
{\rm Ext}^1(L(h_{r,s}),\X^+_{a,b})(\simeq \mathbb{C})
$
(see Proposition \ref{Ext}).
\end{dfn}

The following proposition can be proved in the same way as Proposition \ref{Ext2+0}, but we give a different proof.
\begin{prop}
\label{Ext2+}
Let $(a,b)$ be $(r,s)$ or $(r^\vee,s^\vee)$. Then we have
\begin{align*}
{\rm Ext}^1(\mathcal{Q}(\X^+_{a,b})_{a^\vee,b},\X^+_{a,b})={\rm Ext}^1(\mathcal{Q}(\X^+_{a,b})_{a,b^\vee},\X^+_{a,b})=0.
\end{align*}
\begin{proof}
We prove only ${\rm Ext}^1(\mathcal{Q}(\X^+_{a,b})_{a^\vee,b},\X^+_{a,b})=0$ in the case $(a,b)=(r,s)$. The other cases can be proved in the same way. 

Note that $\mathcal{Q}(\X^+_{r,s})_{r^\vee,s}$ has a submodule isomorphic to $\K_{r,s}$. Assume that ${\rm Ext}^1(\mathcal{Q}(\X^+_{r,s})_{r^\vee,s},\X^+_{r,s})\neq 0$. Then, by Theorem \ref{Ext}, we have 
\begin{align*}
{\rm Ext}^1(\mathcal{Q}(\X^+_{r,s})_{r^\vee,s}/\K_{r,s},\X^+_{r,s})\neq 0.
\end{align*}
Fix a non-trivial extension $E$ in this ${\rm Ext}^1$-group. By Proposition~\ref{sankakuthick}, we see that $L_0$ acts semisimply on the lowest weight space of $E$. Thus, by Proposition \ref{sl2action2}, we see that $L_0$ acts semisimply on $E$. 
Let $u_0\in E[\Delta^+_{r,s;0}]$ be an $L_0$-homogeneous vector whose image in $E/\X^+_{r,s}\simeq \mathcal{Q}(\X^+_{r,s})_{r^\vee,s}/\K_{r,s}$ is a generator. Let $u_1$ be the lowest weight vector of the submodule $\X^+_{r,s}\subset E$. Fix an arbitrary $L_0$-homogeneous vector $u^*_1\in E^*$ such that $\langle u^*_1,u_1\rangle \neq 0$. 
Noting that $[E]\neq 0$ and Proposition~\ref{Ext}, we see that $E$ has one of $\mathcal{E}^{+}(\X^-_{r^\vee,s})_{r,s}$ or $\mathcal{E}^{-}(\X^-_{r^\vee,s})_{r,s}$, as a submodule. Thus, by the structure of $\mathcal{Q}(\X^+_{r,s})_{r^\vee,s}/\K_{r,s}$ and $\mathcal{E}^{\pm}(\X^-_{r^\vee,s})_{r,s}$, we have
\begin{align*}
\langle u^*_1,\sigma(S_{r,s^\vee+p_-})Y_{E}(W^\bullet;z)S_{r,s^\vee+p_-}u_0\rangle \neq 0,
\end{align*}
where $W^\bullet$ is one of $W^\pm$ or $W^0$. In particular, we have
\begin{align}
\label{20220629}
\langle u^*_1,Y_{E}(W^\bullet;z)u_0\rangle \neq 0.
\end{align}
Note that $S_{r^\vee+p_+,s}u_0=0$. 
Thus, by Proposition \ref{FusionInt}, we have
\begin{align*}
0&=\langle u^*_1,{Y}_{E}(W^\bullet;z)S_{r^\vee+p_+,s}u_0\rangle\\
&=z^{-(r^\vee+p_+)s}\prod_{i=1}^{r^\vee+p_+}\prod_{j=1}^{s}(h_{4p_+-1,1}-h_{2r^\vee+2p_+-2i+1,2s-2j+1})\langle u^*_1,{Y}_{E}(W^\bullet;z)u_0\rangle.
\end{align*}
The coefficient in the above equation is nonzero. Then $\langle u^*_1,{Y}_{E}(W^\bullet;z)u_0\rangle{=}0$. But this contradicts (\ref{20220629}). 
\end{proof}
\end{prop}

The following proposition can be proved in the same way as Propositions \ref{Ext2-+0}, \ref{Ext2-0} and \ref{Ext2+-0}, so we omit the proofs.
\begin{prop}
\begin{enumerate}
\item Let $(a,b)$ be $(r^\vee,s)$ or $(r,s^\vee)$. Then we have
\begin{align*}
{\rm Ext}^1(\mathcal{Q}(\X^-_{a,b})_{a^\vee,b},\X^+_{a^\vee,b})&={\rm Ext}^1(\mathcal{Q}(\X^-_{a,b})_{a,b^\vee},\X^+_{a,b^\vee})=0,\\
{\rm Ext}^1(\mathcal{Q}(\X^-_{a,b})_{a^\vee,b},\X^-_{a,b})&={\rm Ext}^1(\mathcal{Q}(\X^-_{a,b})_{a,b^\vee},\X^-_{a,b})=0.
\end{align*}
\item Let $(a,b)$ be $(r,s)$ or $(r^\vee,s^\vee)$. Then we have 
\begin{align*}
{\rm Ext}^1(\mathcal{Q}(\X^+_{a,b})_{a^\vee,b},\X^-_{a^\vee,b})={\rm Ext}^1(\mathcal{Q}(\X^+_{a,b})_{a,b^\vee},\X^-_{a,b^\vee})=0.
\end{align*}
\end{enumerate}
\label{Ext2-+}
\end{prop}

We set
\begin{align}
\label{eq:s-1}
\mathbb{S}_1:=\{(r,s,+),(r^\vee,s^\vee,+),(r^\vee,s,-),(r,s^\vee,-)\}.
\end{align}
\begin{prop}
\label{Ext26+-}
Let $(a,b,\epsilon)\in \mathbb{S}_1$.
Then we have
\begin{align*}
{\rm Ext}^1(\mathcal{Q}(\X^\epsilon_{a,b})_{a^\vee,b},\X^\epsilon_{a^\vee,b^\vee})={\rm Ext}^1(\mathcal{Q}(\X^\epsilon_{a,b})_{a,b^\vee},\X^\epsilon_{a^\vee,b^\vee})=0.
\end{align*}
\begin{proof}
We prove only
\begin{align*}
&{\rm Ext}^1(\mathcal{Q}(\X^+_{r,s})_{r^\vee,s},\X^+_{r^\vee,s^\vee})=0,
&{\rm Ext}^1(\mathcal{Q}(\X^-_{r^\vee,s})_{r,s},\X^-_{r,s^\vee})=0.
\end{align*}
The other equalities can be proved in the same way.

First we prove ${\rm Ext}^1(\mathcal{Q}(\X^+_{r,s})_{r^\vee,s},\X^+_{r^\vee,s^\vee})=0$. 
Recall that $\widetilde{K}(\Delta^+_{r,s;0})$ is the indecomposable Virasoro module satisfying the following exact sequence
\begin{align*}
0\rightarrow L(\Delta^-_{r,s^\vee;0})\rightarrow \widetilde{K}(\Delta^+_{r,s;0})\rightarrow K(\tau)\rightarrow 0,
\end{align*}
where $\tau=(\alpha_{r^\vee,s^\vee},\alpha_{r,s^\vee;1},\alpha_{r^\vee,s^\vee;2})$.
By the structure of Virasoro Verma modules (cf. \cite{FF,IK}), we see that 
\begin{align*}
{\rm Ext}^1_{{\mathcal{L}}}(\widetilde{K}(\Delta^+_{r,s;0})/L(h_{r,s}),L(\Delta^+_{r^\vee,s^\vee;0}))=0.
\end{align*}
Thus, by the structure of the Fock module $F_{r^\vee,s^\vee}$, we have
\begin{align}
{\rm Ext}^1_{{\mathcal{L}}}(\widetilde{K}(\Delta^+_{r,s;0}),L(\Delta^+_{r^\vee,s^\vee;0}))\simeq \C.
\label{BGG00}
\end{align}
Let $V(\Delta^+_{r,s;0})$ be an indecomposable module representing a class in (\ref{BGG00}).
By (\ref{BGG00}), we have 
\begin{align}
\label{202304040}
{\rm Ext}^1_{\mathcal{L}}(V(\Delta^+_{r,s;0}),L(\Delta^+_{r^\vee,s^\vee;0}))=0.
\end{align}
Note that $V(\Delta^+_{r,s;0})$ satisfies the exact sequence
\begin{align*}
0\rightarrow L(\Delta^+_{r^\vee,s^\vee;0},\Delta^-_{r,s^\vee;0})^*\xrightarrow{} V(\Delta^+_{r,s;0})\rightarrow K(\tau)\rightarrow 0.
\end{align*}
By this exact sequence and (\ref{202304040}), we obtain
\begin{align}
{\rm Ext}^1_{\mathcal{L}}(K(\tau),L(\Delta^+_{r^\vee,s^\vee;0}))=0.
\label{2022090301}
\end{align}
Assume that ${\rm Ext}^1(\mathcal{Q}(\X^+_{r,s})_{r^\vee,s},\X^+_{r^\vee,s^\vee})\neq 0$. Then, by Proposition \ref{Ext}, we have 
$
{\rm Ext}^1(\mathcal{Q}(\X^+_{r,s})_{r^\vee,s}/\X^+_{r,s},\X^+_{r^\vee,s^\vee})\neq 0.
$
Fix a non-trivial extension $[E]\in {\rm Ext}^1(\mathcal{Q}(\X^+_{r,s})_{r^\vee,s}/\X^+_{r,s},\X^+_{r^\vee,s^\vee})$.
Note that $L_0$ acts semisimply on $E$. 
By (\ref{2022090301}), we can take a $L_0$-homogeneous generator $u_0$ of $E$ satisfying
\begin{align}
\label{20230416gen}
U(\mathcal{L}).u_0\simeq K(\tau).
\end{align}
Let $u_1$ be the lowest weight vector of the submodule of $E$ isomorphic to $\X^+_{r^\vee,s^\vee}$. 
By (\ref{2022090301}) and (\ref{20230416gen}), we see that there exists a homogeneous vector $u^*_1\in E^*$ satisfying $\langle u^*_1,u_1\rangle \neq 0$ and 
\begin{align}
L_n u^*_1=0,\ \ \ {\rm for}\ n\geq 1.
\label{2022090302}
\end{align}
By Proposition~\ref{Ext}, (\ref{2022090301}) and the nontriviality of $[E]$, we see that $E$ has a submodule isomorphic to $\mathcal{E}^+(\X^+_{r^\vee,s^\vee})^*_{r^\vee,s}$ or $\mathcal{E}^-(\X^+_{r^\vee,s^\vee})^*_{r^\vee,s}$ (for $\mathcal{E}^\pm(\X^+_{r^\vee,s^\vee})_{r^\vee,s}$, see Definition~\ref{kaisei}).
In particular, we have
\begin{align}
\label{Q120}
\langle u^*_1,{Y}_E(W^\bullet;z)u_0\rangle\neq 0,
\end{align}
where $W^\bullet$ is one of $W^\pm$ or $W^0$.
Note that by (\ref{20230416gen}),
\begin{align*}
S_{r^\vee+p_+,s}u_0\in L(h_{r,s}).
\end{align*}
Then by Proposition \ref{FusionInt} and (\ref{2022090302}), we have
\begin{align*}
0&=\langle u^*_1,{Y}_E(W^\bullet;z)S_{r^\vee+p_+,s}u_0\rangle\\
&=z^{-(r^\vee +p_+)s}\prod_{i=1}^{r^\vee+p_+}\prod_{j=1}^{s}(h_{4p_+-1,1}-h_{r^\vee+r+2p_+-2i+1,s+s^\vee-2j+1})\langle u^*_1,{Y}_E(W^\bullet;z)u_0\rangle.
\end{align*}
The coefficient in the above equation is nonzero, so we have $\langle u^*_{1},{Y}_E(W^\bullet;z)u_0\rangle=0$. But this contradicts (\ref{Q120}).

Next we prove ${\rm Ext}^1(\mathcal{Q}(\X^-_{r^\vee,s})_{r,s},\X^-_{r,s^\vee})=0$. Note that, by the structure of Virasoro Verma modules and by the structure of the Fock module $F_{r,s^\vee;1}$,
\begin{align}
{\rm Ext}^1_{\mathcal{L}}(\widetilde{K}(\Delta^-_{r^\vee,s;0})/L(\Delta^+_{r^\vee,s^\vee;0}),L(\Delta^-_{r,s^\vee;0}))\simeq \C
\label{BGG0}
\end{align} 
(see Definitions \ref{dfnC2} for the definitions of Virasoro module $\widetilde{K}(\Delta^-_{r^\vee,s;0})$). 
Let $V(\Delta^-_{r^\vee,s;0})$ be an indecomposable module representing a class in this ${\rm Ext}^1$-group. 
By (\ref{BGG0}), we have 
\begin{align}
\label{202304041}
{\rm Ext}^1_{\mathcal{L}}(V(\Delta^-_{r^\vee,s;0}),L(\Delta^-_{r,s^\vee;0}))=0.
\end{align}
Note that $V(\Delta^-_{r^\vee,s;0})$ satisfies the exact sequence
\begin{align*}
0\rightarrow L(\Delta^-_{r,s^\vee;0},\Delta^+_{r^\vee,s^\vee;1})^*\xrightarrow{} V(\Delta^-_{r^\vee,s;0})\rightarrow K(\tau')\rightarrow 0,
\end{align*}
where $\tau'=(\alpha_{r,s^\vee;1},\alpha_{r^\vee,s^\vee;2},\alpha_{r,s^\vee;3})$.
By this exact sequence and (\ref{202304041}), we obtain
\begin{align}
{\rm Ext}^1_{\mathcal{L}}(K(\tau'),L(\Delta^-_{r,s^\vee;0}))=0.
\label{sp}
\end{align}
Assume that
$
{\rm Ext}^1(\mathcal{Q}(\X^-_{r^\vee,s})_{r,s},\X^-_{r,s^\vee})\neq 0.
$
Then, since ${\rm Ext}^1(\X^-_{r^\vee,s},\X^-_{r,s^\vee})=0$, we have
\begin{align*}
{\rm Ext}^1(\mathcal{Q}(\X^-_{r^\vee,s})_{r,s}/\X^-_{r^\vee,s},\X^-_{r,s^\vee})\neq 0.
\end{align*}
Note that $L_0$ acts semisimply on any extensions of this ${\rm Ext}^1$-group. 
Fix a non-trivial extension $[E]\in {\rm Ext}^1(\mathcal{Q}(\X^-_{r^\vee,s})_{r,s}/\X^-_{r^\vee,s},\X^-_{r,s^\vee})$. 
Then, from Lemma \ref{lemE3} and Proposition \ref{Ext}, $E$ has a submodule isomorphic to $\mathcal{E}^+(\X^+_{r,s})_{r^\vee,s}$ or $\mathcal{E}^-(\X^+_{r,s})_{r^\vee,s}$.
Thus, considering $E$ as an object in $\mathcal{L}_{c_{p_+,p_-}}$-Mod, as its quotient, we have a Virasoro module representing a nontrivial class in 
\begin{align*}
{\rm Ext}^1_{\mathcal{L}}(K(\tau'),L(\Delta^-_{r,s^\vee;0})).
\end{align*}
But this contradicts (\ref{sp}).
\end{proof}
\end{prop}

\begin{prop}
\label{Extcoro2}
Let $(a,b)$ be $(r,s)$ or $(r^\vee,s^\vee)$. Then we have
\begin{align*}
{\rm Ext}^1(\mathcal{Q}(\X^+_{a,b})_{a^\vee,b},L(h_{a,b}))={\rm Ext}^1(\mathcal{Q}(\X^+_{a,b})_{a,b^\vee},L(h_{a,b}))=0.
\end{align*}
\begin{proof}
By the exact sequence
\begin{align*}
0\rightarrow \K_{r,s}\xrightarrow{} \mathcal{Q}(\X^+_{r,s})_{r^\vee,s}\rightarrow \mathcal{Q}(\X^+_{r,s})_{r^\vee,s}/\K_{r,s}\rightarrow 0,
\end{align*}
we have the exact sequence
\begin{align*}
0\rightarrow \C\xrightarrow{} \C\rightarrow {\rm Ext}(\mathcal{Q}(\X^+_{r,s})_{r^\vee,s},L(h_{r,s}))\rightarrow {\rm Ext}^1(\K_{r,s},L(h_{r,s}))\rightarrow 0.
\end{align*}
Thus we have ${\rm Ext}(\mathcal{Q}(\X^+_{r,s})_{r^\vee,s},L(h_{r,s}))\simeq {\rm Ext}^1(\K_{r,s},L(h_{r,s}))$. Assume 
\begin{align*}
{\rm Ext}^1(\K_{r,s},L(h_{r,s}))\neq 0,
\end{align*}
and fix a non-trivial extension $[E]\in {\rm Ext}^1(\K_{r,s},L(h_{r,s}))$. Since 
\begin{align*}
{\rm Ext}^1(L(h_{r,s}),L(h_{r,s}))=0,
\end{align*}
$E$ has a submodule isomorphic to $\K_{r,s}^*$. Thus, by Proposition \ref{Y}, we see that $E$ has $L_0$-nilpotent rank two. Let $\{u_0,u_1\}$ be a basis of the lowest weight space of $E$ such that
\begin{align}
\label{QL}
(L_0-h_{r,s})u_0\in\C^\times u_1.
\end{align}
Then, by Proposition \ref{Y} and (\ref{QL}), we have
$
\sigma(S_{r,s})S_{r,s}u_0\in \C^\times u_1.
$
In particular we have $S_{r,s}u_0\neq 0$. Thus $E$ has $\X^+_{r^\vee,s^\vee}$ as a composition factor. But this is a contradiction. Thus we obtain ${\rm Ext}^1(\mathcal{Q}(\X^+_{r,s})_{r^\vee,s},L(h_{r,s}))=0$. The other equalities can be proved in the same way, so we omit the proofs.
\end{proof}

\begin{comment}
\begin{proof}
Assume ${\rm Ext}^1(\mathcal{Q}(h_{r,s}),L(h_{r,s}))\neq 0$. Fix a non-trivial exact sequence
\begin{align*}
0\rightarrow L(h_{r,s})\xrightarrow{} E\rightarrow \mathcal{Y}(h_{r,s})\rightarrow 0.
\end{align*}
Let $E_0$ be the submodule of $E$ generated from $\X^+_{r,s}$ and $\X^+_{r^\vee,s^\vee}$ at level 1. We have
\begin{align*}
E_0\simeq \mathcal{I}^*_{r,s}\oplus \mathcal{I}^*_{r^\vee,s^\vee}.
\end{align*}
Let $u$ be the lowest weight vector of $E$ which corresponds to the lowest weight vector of $L(h_{r,s})$ at level 0. Then $S^*_{r,s}S_{r,s}u$ and $S^*_{r^\vee,s^\vee}S_{r^\vee,s^\vee}u$ generate ${\rm S}(E)=L(h_{r,s})\oplus L(h_{r,s})$. On the other hand, we have
\begin{align*}
&(L_0-h_{r,s})u=c_1S^*_{r,s}S_{r,s}u
&(L_0-h_{r,s})u=c_2S^*_{r^\vee,s^\vee}S_{r^\vee,s^\vee}u
\end{align*}
where $c_1$ and $c_2$ are non-zero constants. In particular $c_1S^*_{r,s}S_{r,s}u=c_2S^*_{r^\vee,s^\vee}S_{r^\vee,s^\vee}u$. But this is a contradiction.
\end{proof}
\end{comment}
\end{prop}

\begin{prop}
\label{Extcoro1}
Let $(a,b)$ be $(r^\vee,s)$ or $(r,s^\vee)$. Then we have
\begin{align*}
{\rm Ext}^1(\mathcal{Q}(\X^-_{a,b})_{a^\vee,b},L(h_{a^\vee,b}))={\rm Ext}^1(\mathcal{Q}(\X^-_{a,b})_{a,b^\vee},L(h_{a,b^\vee}))=0.
\end{align*}
\begin{proof}
By Proposition \ref{Ext}, we have 
\begin{align*}
{\rm Ext}^1(\mathcal{Q}(\X^-_{r^\vee,s})_{r,s},L(h_{r,s}))\simeq {\rm Ext}^1(\mathcal{Q}(\X^-_{r^\vee,s})_{r,s}/\X^-_{r^\vee,s},L(h_{r,s})).
\end{align*}
Assume ${\rm Ext}^1(\mathcal{Q}(\X^-_{r^\vee,s})_{r,s}/\X^-_{r^\vee,s},L(h_{r,s}))\neq 0$. Then any non-trivial extension of ${\rm Ext}^1(\mathcal{Q}(\X^-_{r^\vee,s})_{r,s}/\X^-_{r^\vee,s},L(h_{r,s}))$ has an indecomposable submodule corresponding to a class in ${\rm Ext}^1(\X^-_{r^\vee,s},\K_{r,s}^*)$. Thus we have
$
{\rm Ext}^1(\K_{r,s},\X^-_{r^\vee,s})\neq 0.
$
Since ${\rm Ext}^1(L(h_{r,s}),\X^-_{r^\vee,s})=0$, any non-trivial extension of ${\rm Ext}^1(\K_{r,s},\X^-_{r^\vee,s})$ has a submodule corresponding to a nontrivial class in ${\rm Ext}^1(\X^+_{r,s},\X^-_{r^\vee,s})$.
In particular, considering $E$ as an object in $\mathcal{L}_{c_{p_+,p_-}}$-Mod, we have
\begin{align*}
{\rm Ext}^1_{\mathcal{L}}(L(h_{r,s},\Delta^+_{r,s;0}),L(\Delta^-_{r^\vee,s;0}))\neq 0.
\end{align*}
But this contradicts the quotient structure of the Verma module $M(h_{r,s},c_{p_+,p_-})$ (cf.~\cite{FF,IK}).
The other equalities can be proved in the same way. 
\end{proof}
\end{prop}
We set
\begin{align*}
\mathbb{S}_2=&\bigl\{(+,r,s,r^\vee,s),(+,r,s,r,s^\vee),(+,r^\vee,s^\vee,r^\vee,s),(+,r^\vee,s^\vee,r,s^\vee),\\
&(-,r^\vee,s,r,s),(-,r^\vee,s,r^\vee,s^\vee),(-,r,s^\vee,r,s),(-,r,s^\vee,r^\vee,s^\vee)\bigr\}.
\end{align*}
The following is a summary of Propositions~\ref{Ext2+}-\ref{Extcoro1}.
\begin{prop}
For each $(\epsilon,a,b,c,d)\in \mathbb{S}_2$,
we have
\begin{align*}
&{\rm Ext}^1(\mathcal{Q}(\X^\epsilon_{a,b})_{c,d},L(h_{r,s}))=0,\ \ \ {\rm Ext}^1(\mathcal{Q}(\X^\epsilon_{a,b})_{c,d},\X^{\epsilon}_{a,b})=0,\\
&{\rm Ext}^1(\mathcal{Q}(\X^\epsilon_{a,b})_{c,d},\X^{\epsilon}_{a^\vee,b^\vee})=0,\ \ \ \ {\rm Ext}^1(\mathcal{Q}(\X^\epsilon_{a,b})_{c,d},\X^{-\epsilon}_{c,d})=0.
\end{align*}
\label{ExtQ}
\end{prop}

In the following, let us introduce the socle structure of the logarithmic modules $\mathcal{P}^\pm_{\bullet,\bullet}$.
By Propositions \ref{Ext}, \ref{Ext2+} and \ref{Ext2-+}, we have the following lemma.
\begin{lem} 
\label{koreda}
Let $(\epsilon,a,b,c,d)\in \mathbb{S}_2$.
Then, any indecomposable module whose composition factors are the same as those of $\mathcal{Q}(\X^\epsilon_{a,b})_{c,d}$ is isomorphic to $\mathcal{Q}(\X^\epsilon_{a,b})_{c,d}$.
\end{lem}
Recall the notation (\ref{eq:s-1}).
\begin{prop}
\label{quotient}
Let $(a,b,\epsilon)\in \mathbb{S}_1$.
Then the logarithmic module $\mathcal{P}^\epsilon_{a,b}$ has the following sequences of submodules
\begin{align*}
&U_1(\mathcal{P}^\epsilon_{a,b})\subset U_2(\mathcal{P}^\epsilon_{a,b})\subset U_3(\mathcal{P}^\epsilon_{a,b})\subset U_4(\mathcal{P}^\epsilon_{a,b})=\mathcal{P}^\epsilon_{a,b},\\
&V_1(\mathcal{P}^\epsilon_{a,b})\subset V_2(\mathcal{P}^\epsilon_{a,b})\subset V_3(\mathcal{P}^\epsilon_{a,b})\subset V_4(\mathcal{P}^\epsilon_{a,b})=\mathcal{P}^\epsilon_{a,b},
\end{align*}
such that
\begin{align*}
&U_1= \mathcal{Q}(\X^{\epsilon}_{a,b})_{a^\vee,b},\ U_2/U_1=U_3/U_2= \mathcal{Q}(\X^{-\epsilon}_{a,b^\vee})_{a^\vee,b^\vee},\ U_4/U_3=\mathcal{Q}(\X^{\epsilon}_{a,b})_{a^\vee,b},\\
&V_1= \mathcal{Q}(\X^{\epsilon}_{a,b})_{a,b^\vee},\ V_2/V_1=V_3/V_2= \mathcal{Q}(\X^{-\epsilon}_{a^\vee,b})_{a^\vee,b^\vee},\ V_4/V_3=\mathcal{Q}(\X^{\epsilon}_{a,b})_{a,b^\vee}.
\end{align*}
\begin{proof}
Let us show that $\mathcal{P}^+_{r,s}$ has a sequence
\begin{align*}
U_1(\mathcal{P}^+_{r,s})\subset U_2(\mathcal{P}^+_{r,s})\subset U_3(\mathcal{P}^+_{r,s})\subset U_4(\mathcal{P}^+_{r,s})=\mathcal{P}^+_{r,s}
\end{align*}
such that
\begin{align*}
U_1= \mathcal{Q}(\X^{+}_{r,s})_{r^\vee,s},\ U_2/U_1=U_3/U_2= \mathcal{Q}(\X^{-}_{r,s^\vee})_{r^\vee,s^\vee},\ U_4/U_3=\mathcal{Q}(\X^{+}_{r,s})_{r^\vee,s}.
\end{align*}
The other cases can be proved in the same way.

By Theorem~\ref{logarithmic1}, we see that $\mathcal{P}^+_{r,s}$ has a logarithmic submodule $M$ satisfying
\begin{align*}
&\mathcal{V}^+_{r,s}\subset M,
&M/\mathcal{V}^+_{r,s}\simeq \mathcal{V}^-_{r^\vee,s}.
\end{align*}
From the structure of $\mathcal{V}^+_{r,s}$ and $\mathcal{V}^-_{r^\vee,s}$, we see that $M$ has two indecomposable subquotients whose composition factors are the same as those of $\mathcal{Q}(\X^+_{r,s})_{r^\vee,s}$ and $\mathcal{Q}(\X^-_{r,s^\vee})_{r^\vee,s^\vee}$, respectively.
Then, from Lemma \ref{koreda}, we see that $M$ satisfies
\begin{align}
\label{al04280}
&\mathcal{Q}(\X^+_{r,s})_{r^\vee,s}\subset M,
&M/\mathcal{Q}(\X^+_{r,s})_{r^\vee,s}\simeq \mathcal{Q}(\X^-_{r,s^\vee})_{r^\vee,s^\vee}.
\end{align}
Set $N=\mathcal{P}^+_{r,s}/M$. 
From Theorem~\ref{logarithmic1}, we see that $N$ is logarithmic and satisfies
\begin{align*}
&\mathcal{V}^-_{r,s^\vee}\subset N,
&N/\mathcal{V}^-_{r,s^\vee}\simeq \mathcal{V}^+_{r^\vee,s^\vee}.
\end{align*}
From the structure of $\mathcal{V}^-_{r,s^\vee}$ and $\mathcal{V}^+_{r^\vee,s^\vee}$, we see that $N$ has two indecomposable subquotients whose composition factors are the same as those of $\mathcal{Q}(\X^-_{r,s^\vee})_{r^\vee,s^\vee}$ and $\mathcal{Q}(\X^+_{r,s})_{r^\vee,s}$, respectively.
Then, from Lemma \ref{koreda}, we see that $N$ satisfies
\begin{align}
\label{al04281}
&\mathcal{Q}(\X^-_{r,s^\vee})_{r^\vee,s^\vee}\subset N,
&N/\mathcal{Q}(\X^-_{r,s^\vee})_{r^\vee,s^\vee}\simeq \mathcal{Q}(\X^+_{r,s})_{r^\vee,s}.
\end{align}
By (\ref{al04280}) and (\ref{al04281}), we obtain the desired assertion.
\end{proof}
\end{prop}

By Propositions \ref{ExtQ} and \ref{quotient}, we obtain the following proposition.
\begin{prop}
For each $(a,b,\epsilon)\in \mathbb{S}_1$,
we have
\begin{align*}
&{\rm Ext}^1(\mathcal{P}^\epsilon_{a,b},L(h_{r,s}))={\rm Ext}^1(\mathcal{P}^\epsilon_{a,b},\X^\epsilon_{a^\vee,b^\vee})=0,\\
&{\rm Ext}^1(\mathcal{P}^\epsilon_{a,b},\X^{-\epsilon}_{a^\vee,b})={\rm Ext}^1(\mathcal{P}^\epsilon_{a,b},\X^{-\epsilon}_{a,b^\vee})=0.
\end{align*}
\label{P4}
\end{prop}

By Proposition \ref{quotient}, we obtain the following proposition.
\begin{prop}
\label{logarithmic2}
Let $(a,b)=(r,s)$ or $(r^\vee,s^\vee)$. Each of the logarithmic modules $\mathcal{P}^+_{a,b},\mathcal{P}^-_{a,b^\vee}\in C^{thick}_{r,s}$ is generated from the top composition factor and has the following length five socle series:
\mbox{}
\begin{enumerate}
\item
For the socle series of $\mathcal{P}^+_{a,b}$, we have
\begin{align*}
\hspace{-4mm}&{\rm S}_1=\mathcal{X}^+_{a,b},\quad 
{\rm S}_2/{\rm S}_1=2\mathcal{X}^{-}_{a,b^\vee}\oplus L(h_{r,s})\oplus 2\mathcal{X}^-_{a^\vee,b},\quad
{\rm S}_3/{\rm S}_2=2\mathcal{X}^+_{a,b}\oplus4\mathcal{X}^+_{a^{\vee},b^{\vee}},\\
\hspace{-4mm}&{\rm S}_4/{\rm S}_3=2\mathcal{X}^-_{a^{\vee},b}\oplus L(h_{r,s})\oplus 2\mathcal{X}^{-}_{a,b^{\vee}},\quad
\mathcal{P}^+_{a,b}/{\rm S}_4=\mathcal{X}^+_{a,b}.
\end{align*}
\item
For the socle series of $\mathcal{P}^-_{a,b^{\vee}}$, we have
\begin{align*}
\hspace{-6mm}&{\rm S}_1=\mathcal{X}^-_{a,b^\vee},\quad
{\rm S}_2/{\rm S}_1=2\mathcal{X}^{+}_{a,b}\oplus 2\mathcal{X}^+_{a^\vee,b^\vee},\quad
{\rm S}_3/{\rm S}_2=2\mathcal{X}^-_{a,b^\vee}\oplus 2 L(h_{r,s}) \oplus 4\mathcal{X}^-_{a^\vee,b},\\
\hspace{-6mm}&{\rm S}_4/{\rm S}_3=2\mathcal{X}^+_{a^\vee,b^\vee}\oplus 2\mathcal{X}^{+}_{a,b},\quad
\mathcal{P}^-_{a,b^{\vee}}/{\rm S}_4=\mathcal{X}^-_{a,b^\vee}.
\end{align*}
\end{enumerate}
\end{prop}

\begin{figure}[htbp]
  \centering
\begin{tikzpicture}[scale=1.05]
\node[inner sep=0.5pt,scale=0.8] (a0) at (0, 0.43) {$\mathcal{P}^+_{r,s}$};
\node[inner sep=0.5pt,scale=0.8] (a1) at (0, 0) {${\heartsuit}$};
\node[inner sep=0.5pt,scale=0.8] (b1) at (-1.5,-0.8) {$\clubsuit$};
\node[inner sep=0.5pt,scale=0.8] (b2) at (-0.75,-0.8) {$\clubsuit$};
\node[inner sep=0.5pt,scale=0.8] (b3) at (0,-0.8) {$\triangle$};
\node[inner sep=0.5pt,scale=0.8] (b4) at (0.75,-0.8) {$\spadesuit$};
\node[inner sep=0.5pt,scale=0.8] (b5) at (1.5,-0.8) {$\spadesuit$};
\node[inner sep=0.5pt,scale=0.8] (c1) at (-1.5-0.375,-1.6) {$\heartsuit$};
\node[inner sep=0.5pt,scale=0.8] (c2) at (-0.75-0.375,-1.6) {$\diamondsuit$};
\node[inner sep=0.5pt,scale=0.8] (c3) at (-0.375,-1.6) {$\diamondsuit$};
\node[inner sep=0.5pt,scale=0.8] (c4) at (0.75-0.375,-1.6) {$\diamondsuit$};
\node[inner sep=0.5pt,scale=0.8] (c5) at (1.5-0.375,-1.6) {$\diamondsuit$};
\node[inner sep=0.5pt,scale=0.8] (c6) at (1.5+0.375,-1.6) {$\heartsuit$};
\node[inner sep=0.5pt,scale=0.8] (d1) at (-1.5,-2.4) {$\spadesuit$};
\node[inner sep=0.5pt,scale=0.8] (d2) at (-0.75,-2.4) {$\spadesuit$};
\node[inner sep=0.5pt,scale=0.8] (d3) at (0,-2.4) {$\triangle$};
\node[inner sep=0.5pt,scale=0.8] (d4) at (0.75,-2.4) {$\clubsuit$};
\node[inner sep=0.5pt,scale=0.8] (d5) at (1.5,-2.4) {$\clubsuit$};
\node[inner sep=0.5pt,scale=0.8] (e1) at (0, -3.2) {$\heartsuit$};

\node[inner sep=0.5pt,scale=0.8] (0a) at (5.5, 0.43) {$\mathcal{P}^-_{r^\vee,s}$};
\node[inner sep=0.5pt,scale=0.8] (1a) at (5.5, 0) {${\clubsuit}$};
\node[inner sep=0.5pt,scale=0.8] (1b) at (5.5-1.75,-0.8) {$\heartsuit$};
\node[inner sep=0.5pt,scale=0.8] (2b) at (5.5-1.1,-0.8) {$\heartsuit$};
\node[inner sep=0.5pt,scale=0.8] (3b) at (5.5+1.1,-0.8) {$\diamondsuit$};
\node[inner sep=0.5pt,scale=0.8] (4b) at (5.5+1.75,-0.8) {$\diamondsuit$};
\node[inner sep=0.5pt,scale=0.8] (1c) at (5.5-1.75+5.5-1.75-5.5+1.1,-1.6) {$\clubsuit$};
\node[inner sep=0.5pt,scale=0.8] (2c) at (5.5-1.75,-1.6) {$\spadesuit$};
\node[inner sep=0.5pt,scale=0.8] (3c) at (5.5-1.1,-1.6) {$\spadesuit$};
\node[inner sep=0.5pt,scale=0.8] (4c) at (5.5-0.4,-1.6) {$\triangle$};
\node[inner sep=0.5pt,scale=0.8] (5c) at (5.5+0.4,-1.6) {$\triangle$};
\node[inner sep=0.5pt,scale=0.8] (6c) at (5.5+1.1,-1.6) {$\spadesuit$};
\node[inner sep=0.5pt,scale=0.8] (7c) at (5.5+1.75,-1.6) {$\spadesuit$};
\node[inner sep=0.5pt,scale=0.8] (8c) at (5.5+1.75+1.75-1.1,-1.6) {$\clubsuit$};
\node[inner sep=0.5pt,scale=0.8] (1d) at (5.5-1.75,-2.4) {$\diamondsuit$};
\node[inner sep=0.5pt,scale=0.8] (2d) at (5.5-1.1,-2.4) {$\diamondsuit$};
\node[inner sep=0.5pt,scale=0.8] (3d) at (5.5+1.1,-2.4) {$\heartsuit$};
\node[inner sep=0.5pt,scale=0.8] (4d) at (5.5+1.75,-2.4) {$\heartsuit$};
\node[inner sep=0.5pt,scale=0.8] (1e) at (5.5, -3.2) {$\clubsuit$};

\draw[arrows = {-Stealth[scale=0.6]}] (a1) to (b1);
\draw[arrows = {-Stealth[scale=0.6]}] (a1) to (b2);
\draw[arrows = {-Stealth[scale=0.6]}] (a1) to (b3);
\draw[arrows = {-Stealth[scale=0.6]}] (a1) to (b4);
\draw[arrows = {-Stealth[scale=0.6]}] (a1) to (b5);
\draw[arrows = {-Stealth[scale=0.6]}] (b1) to (c1);
\draw[arrows = {-Stealth[scale=0.6]}] (b1) to (c2);
\draw[arrows = {-Stealth[scale=0.6]}] (b1) to (c4);
\draw[arrows = {-Stealth[scale=0.6]}] (b2) to (c1);
\draw[arrows = {-Stealth[scale=0.6]}] (b2) to (c3);
\draw[arrows = {-Stealth[scale=0.6]}] (b2) to (c5);
\draw[arrows = {-Stealth[scale=0.6]}] (b3) to (c1);
\draw[arrows = {-Stealth[scale=0.6]}](b3) to (c3);
\draw[arrows = {-Stealth[scale=0.6]}] (b3) to (c4);
\draw[arrows = {-Stealth[scale=0.6]}] (b3) to (c6);
\draw[arrows = {-Stealth[scale=0.6]}](b4) to (c2);
\draw[arrows = {-Stealth[scale=0.6]}] (b4) to (c3);
\draw[arrows = {-Stealth[scale=0.6]}] (b4) to (c6);
\draw[arrows = {-Stealth[scale=0.6]}] (b5) to (c4);
\draw[arrows = {-Stealth[scale=0.6]}] (b5) to (c5);
\draw[arrows = {-Stealth[scale=0.6]}] (b5) to (c6);
\draw[arrows = {-Stealth[scale=0.6]}] (c1) to (d1);
\draw[arrows = {-Stealth[scale=0.6]}] (c1) to (d2);
\draw[arrows = {-Stealth[scale=0.6]}] (c1) to (d3);
\draw[arrows = {-Stealth[scale=0.6]}] (c2) to (d1);
\draw[arrows = {-Stealth[scale=0.6]}] (c2) to (d4);
\draw[arrows = {-Stealth[scale=0.6]}] (c3) to (d1);
\draw[arrows = {-Stealth[scale=0.6]}] (c3) to (d3);
\draw[arrows = {-Stealth[scale=0.6]}] (c3) to (d5);
\draw[arrows = {-Stealth[scale=0.6]}] (c4) to (d2);
\draw[arrows = {-Stealth[scale=0.6]}] (c4) to (d3);
\draw[arrows = {-Stealth[scale=0.6]}] (c4) to (d4);
\draw[arrows = {-Stealth[scale=0.6]}] (c5) to (d2);
\draw[arrows = {-Stealth[scale=0.6]}] (c5) to (d5);
\draw[arrows = {-Stealth[scale=0.6]}] (c6) to (d3);
\draw[arrows = {-Stealth[scale=0.6]}] (c6) to (d4);
\draw[arrows = {-Stealth[scale=0.6]}] (c6) to (d5);
\draw[arrows = {-Stealth[scale=0.6]}] (d1) to (e1);
\draw[arrows = {-Stealth[scale=0.6]}] (d2) to (e1);
\draw[arrows = {-Stealth[scale=0.6]}] (d3) to (e1);
\draw[arrows = {-Stealth[scale=0.6]}] (d4) to (e1);
\draw[arrows = {-Stealth[scale=0.6]}] (d5) to (e1);

\draw[arrows = {-Stealth[scale=0.6]}] (1a) to (1b);
\draw[arrows = {-Stealth[scale=0.6]}] (1a) to (2b);
\draw[arrows = {-Stealth[scale=0.6]}] (1a) to (3b);
\draw[arrows = {-Stealth[scale=0.6]}] (1a) to (4b);
\draw[arrows = {-Stealth[scale=0.6]}] (1b) to (1c);
\draw[arrows = {-Stealth[scale=0.6]}] (1b) to (2c);
\draw[arrows = {-Stealth[scale=0.6]}] (1b) to (3c);
\draw[arrows = {-Stealth[scale=0.6]}] (1b) to (4c);
\draw[arrows = {-Stealth[scale=0.6]}] (2b) to (1c);
\draw[arrows = {-Stealth[scale=0.6]}] (2b) to (5c);
\draw[arrows = {-Stealth[scale=0.6]}] (2b) to (6c);
\draw[arrows = {-Stealth[scale=0.6]}] (2b) to (7c);
\draw[arrows = {-Stealth[scale=0.6]}] (3b) to (2c);
\draw[arrows = {-Stealth[scale=0.6]}] (3b) to (4c);
\draw[arrows = {-Stealth[scale=0.6]}] (3b) to (6c);
\draw[arrows = {-Stealth[scale=0.6]}] (3b) to (8c);
\draw[arrows = {-Stealth[scale=0.6]}] (4b) to (3c);
\draw[arrows = {-Stealth[scale=0.6]}] (4b) to (5c);
\draw[arrows = {-Stealth[scale=0.6]}] (4b) to (7c);
\draw[arrows = {-Stealth[scale=0.6]}] (4b) to (8c);
\draw[arrows = {-Stealth[scale=0.6]}] (1c) to (1d);
\draw[arrows = {-Stealth[scale=0.6]}] (1c) to (2d);
\draw[arrows = {-Stealth[scale=0.6]}] (2c) to (1d);
\draw[arrows = {-Stealth[scale=0.6]}] (2c) to (3d);
\draw[arrows = {-Stealth[scale=0.6]}] (3c) to (2d);
\draw[arrows = {-Stealth[scale=0.6]}] (3c) to (3d);
\draw[arrows = {-Stealth[scale=0.6]}] (4c) to (1d);
\draw[arrows = {-Stealth[scale=0.6]}] (4c) to (3d);
\draw[arrows = {-Stealth[scale=0.6]}] (5c) to (2d);
\draw[arrows = {-Stealth[scale=0.6]}] (5c) to (4d);
\draw[arrows = {-Stealth[scale=0.6]}] (6c) to (1d);
\draw[arrows = {-Stealth[scale=0.6]}] (6c) to (4d);
\draw[arrows = {-Stealth[scale=0.6]}] (7c) to (2d);
\draw[arrows = {-Stealth[scale=0.6]}] (7c) to (4d);
\draw[arrows = {-Stealth[scale=0.6]}] (8c) to (3d);
\draw[arrows = {-Stealth[scale=0.6]}] (8c) to (4d);
\draw[arrows = {-Stealth[scale=0.6]}] (1d) to (1e);
\draw[arrows = {-Stealth[scale=0.6]}] (2d) to (1e);
\draw[arrows = {-Stealth[scale=0.6]}] (3d) to (1e);
\draw[arrows = {-Stealth[scale=0.6]}] (4d) to (1e);

%%第二列
\node[inner sep=0.5pt,scale=0.8] (a02) at (0, 0.43-4.5) {$\mathcal{P}^+_{r^\vee,s^\vee}$};
\node[inner sep=0.5pt,scale=0.8] (a12) at (0, -4.5) {${\diamondsuit}$};
\node[inner sep=0.5pt,scale=0.8] (b12) at (-1.5,-0.8-4.5) {$\spadesuit$};
\node[inner sep=0.5pt,scale=0.8] (b22) at (-0.75,-0.8-4.5) {$\spadesuit$};
\node[inner sep=0.5pt,scale=0.8] (b32) at (0,-0.8-4.5) {$\triangle$};
\node[inner sep=0.5pt,scale=0.8] (b42) at (0.75,-0.8-4.5) {$\clubsuit$};
\node[inner sep=0.5pt,scale=0.8] (b52) at (1.5,-0.8-4.5) {$\clubsuit$};
\node[inner sep=0.5pt,scale=0.8] (c12) at (-1.5-0.375,-1.6-4.5) {$\diamondsuit$};
\node[inner sep=0.5pt,scale=0.8] (c22) at (-0.75-0.375,-1.6-4.5) {$\heartsuit$};
\node[inner sep=0.5pt,scale=0.8] (c32) at (-0.375,-1.6-4.5) {$\heartsuit$};
\node[inner sep=0.5pt,scale=0.8] (c42) at (0.75-0.375,-1.6-4.5) {$\heartsuit$};
\node[inner sep=0.5pt,scale=0.8] (c52) at (1.5-0.375,-1.6-4.5) {$\heartsuit$};
\node[inner sep=0.5pt,scale=0.8] (c62) at (1.5+0.375,-1.6-4.5) {$\diamondsuit$};
\node[inner sep=0.5pt,scale=0.8] (d12) at (-1.5,-2.4-4.5) {$\clubsuit$};
\node[inner sep=0.5pt,scale=0.8] (d22) at (-0.75,-2.4-4.5) {$\clubsuit$};
\node[inner sep=0.5pt,scale=0.8] (d32) at (0,-2.4-4.5) {$\triangle$};
\node[inner sep=0.5pt,scale=0.8] (d42) at (0.75,-2.4-4.5) {$\spadesuit$};
\node[inner sep=0.5pt,scale=0.8] (d52) at (1.5,-2.4-4.5) {$\spadesuit$};
\node[inner sep=0.5pt,scale=0.8] (e12) at (0, -3.2-4.5) {$\diamondsuit$};

\node[inner sep=0.5pt,scale=0.8] (0a2) at (5.5, 0.43-4.5) {$\mathcal{P}^-_{r,s^\vee}$};
\node[inner sep=0.5pt,scale=0.8] (1a2) at (5.5, -4.5) {${\spadesuit}$};
\node[inner sep=0.5pt,scale=0.8] (1b2) at (5.5-1.75,-0.8-4.5) {$\diamondsuit$};
\node[inner sep=0.5pt,scale=0.8] (2b2) at (5.5-1.1,-0.8-4.5) {$\diamondsuit$};
\node[inner sep=0.5pt,scale=0.8] (3b2) at (5.5+1.1,-0.8-4.5) {$\heartsuit$};
\node[inner sep=0.5pt,scale=0.8] (4b2) at (5.5+1.75,-0.8-4.5) {$\heartsuit$};
\node[inner sep=0.5pt,scale=0.8] (1c2) at (5.5-1.75+5.5-1.75-5.5+1.1,-1.6-4.5) {$\spadesuit$};
\node[inner sep=0.5pt,scale=0.8] (2c2) at (5.5-1.75,-1.6-4.5) {$\clubsuit$};
\node[inner sep=0.5pt,scale=0.8] (3c2) at (5.5-1.1,-1.6-4.5) {$\clubsuit$};
\node[inner sep=0.5pt,scale=0.8] (4c2) at (5.5-0.4,-1.6-4.5) {$\triangle$};
\node[inner sep=0.5pt,scale=0.8] (5c2) at (5.5+0.4,-1.6-4.5) {$\triangle$};
\node[inner sep=0.5pt,scale=0.8] (6c2) at (5.5+1.1,-1.6-4.5) {$\clubsuit$};
\node[inner sep=0.5pt,scale=0.8] (7c2) at (5.5+1.75,-1.6-4.5) {$\clubsuit$};
\node[inner sep=0.5pt,scale=0.8] (8c2) at (5.5+1.75+1.75-1.1,-1.6-4.5) {$\spadesuit$};
\node[inner sep=0.5pt,scale=0.8] (1d2) at (5.5-1.75,-2.4-4.5) {$\heartsuit$};
\node[inner sep=0.5pt,scale=0.8] (2d2) at (5.5-1.1,-2.4-4.5) {$\heartsuit$};
\node[inner sep=0.5pt,scale=0.8] (3d2) at (5.5+1.1,-2.4-4.5) {$\diamondsuit$};
\node[inner sep=0.5pt,scale=0.8] (4d2) at (5.5+1.75,-2.4-4.5) {$\diamondsuit$};
\node[inner sep=0.5pt,scale=0.8] (1e2) at (5.5, -3.2-4.5) {$\spadesuit$};

\draw[arrows = {-Stealth[scale=0.6]}] (a12) to (b12);
\draw[arrows = {-Stealth[scale=0.6]}] (a12) to (b22);
\draw[arrows = {-Stealth[scale=0.6]}] (a12) to (b32);
\draw[arrows = {-Stealth[scale=0.6]}] (a12) to (b42);
\draw[arrows = {-Stealth[scale=0.6]}] (a12) to (b52);
\draw[arrows = {-Stealth[scale=0.6]}] (b12) to (c12);
\draw[arrows = {-Stealth[scale=0.6]}] (b12) to (c22);
\draw[arrows = {-Stealth[scale=0.6]}] (b12) to (c42);
\draw[arrows = {-Stealth[scale=0.6]}] (b22) to (c12);
\draw[arrows = {-Stealth[scale=0.6]}] (b22) to (c32);
\draw[arrows = {-Stealth[scale=0.6]}] (b22) to (c52);
\draw[arrows = {-Stealth[scale=0.6]}] (b32) to (c12);
\draw[arrows = {-Stealth[scale=0.6]}](b32) to (c32);
\draw[arrows = {-Stealth[scale=0.6]}] (b32) to (c42);
\draw[arrows = {-Stealth[scale=0.6]}] (b32) to (c62);
\draw[arrows = {-Stealth[scale=0.6]}](b42) to (c22);
\draw[arrows = {-Stealth[scale=0.6]}] (b42) to (c32);
\draw[arrows = {-Stealth[scale=0.6]}] (b42) to (c62);
\draw[arrows = {-Stealth[scale=0.6]}] (b52) to (c42);
\draw[arrows = {-Stealth[scale=0.6]}] (b52) to (c52);
\draw[arrows = {-Stealth[scale=0.6]}] (b52) to (c62);
\draw[arrows = {-Stealth[scale=0.6]}] (c12) to (d12);
\draw[arrows = {-Stealth[scale=0.6]}] (c12) to (d22);
\draw[arrows = {-Stealth[scale=0.6]}] (c12) to (d32);
\draw[arrows = {-Stealth[scale=0.6]}] (c22) to (d12);
\draw[arrows = {-Stealth[scale=0.6]}] (c22) to (d42);
\draw[arrows = {-Stealth[scale=0.6]}] (c32) to (d12);
\draw[arrows = {-Stealth[scale=0.6]}] (c32) to (d32);
\draw[arrows = {-Stealth[scale=0.6]}] (c32) to (d52);
\draw[arrows = {-Stealth[scale=0.6]}] (c42) to (d22);
\draw[arrows = {-Stealth[scale=0.6]}] (c42) to (d32);
\draw[arrows = {-Stealth[scale=0.6]}] (c42) to (d42);
\draw[arrows = {-Stealth[scale=0.6]}] (c52) to (d22);
\draw[arrows = {-Stealth[scale=0.6]}] (c52) to (d52);
\draw[arrows = {-Stealth[scale=0.6]}] (c62) to (d32);
\draw[arrows = {-Stealth[scale=0.6]}] (c62) to (d42);
\draw[arrows = {-Stealth[scale=0.6]}] (c62) to (d52);
\draw[arrows = {-Stealth[scale=0.6]}] (d12) to (e12);
\draw[arrows = {-Stealth[scale=0.6]}] (d22) to (e12);
\draw[arrows = {-Stealth[scale=0.6]}] (d32) to (e12);
\draw[arrows = {-Stealth[scale=0.6]}] (d42) to (e12);
\draw[arrows = {-Stealth[scale=0.6]}] (d52) to (e12);

\draw[arrows = {-Stealth[scale=0.6]}] (1a2) to (1b2);
\draw[arrows = {-Stealth[scale=0.6]}] (1a2) to (2b2);
\draw[arrows = {-Stealth[scale=0.6]}] (1a2) to (3b2);
\draw[arrows = {-Stealth[scale=0.6]}] (1a2) to (4b2);
\draw[arrows = {-Stealth[scale=0.6]}] (1b2) to (1c2);
\draw[arrows = {-Stealth[scale=0.6]}] (1b2) to (2c2);
\draw[arrows = {-Stealth[scale=0.6]}] (1b2) to (3c2);
\draw[arrows = {-Stealth[scale=0.6]}] (1b2) to (4c2);
\draw[arrows = {-Stealth[scale=0.6]}] (2b2) to (1c2);
\draw[arrows = {-Stealth[scale=0.6]}] (2b2) to (5c2);
\draw[arrows = {-Stealth[scale=0.6]}] (2b2) to (6c2);
\draw[arrows = {-Stealth[scale=0.6]}] (2b2) to (7c2);
\draw[arrows = {-Stealth[scale=0.6]}] (3b2) to (2c2);
\draw[arrows = {-Stealth[scale=0.6]}] (3b2) to (4c2);
\draw[arrows = {-Stealth[scale=0.6]}] (3b2) to (6c2);
\draw[arrows = {-Stealth[scale=0.6]}] (3b2) to (8c2);
\draw[arrows = {-Stealth[scale=0.6]}] (4b2) to (3c2);
\draw[arrows = {-Stealth[scale=0.6]}] (4b2) to (5c2);
\draw[arrows = {-Stealth[scale=0.6]}] (4b2) to (7c2);
\draw[arrows = {-Stealth[scale=0.6]}] (4b2) to (8c2);
\draw[arrows = {-Stealth[scale=0.6]}] (1c2) to (1d2);
\draw[arrows = {-Stealth[scale=0.6]}] (1c2) to (2d2);
\draw[arrows = {-Stealth[scale=0.6]}] (2c2) to (1d2);
\draw[arrows = {-Stealth[scale=0.6]}] (2c2) to (3d2);
\draw[arrows = {-Stealth[scale=0.6]}] (3c2) to (2d2);
\draw[arrows = {-Stealth[scale=0.6]}] (3c2) to (3d2);
\draw[arrows = {-Stealth[scale=0.6]}] (4c2) to (1d2);
\draw[arrows = {-Stealth[scale=0.6]}] (4c2) to (3d2);
\draw[arrows = {-Stealth[scale=0.6]}] (5c2) to (2d2);
\draw[arrows = {-Stealth[scale=0.6]}] (5c2) to (4d2);
\draw[arrows = {-Stealth[scale=0.6]}] (6c2) to (1d2);
\draw[arrows = {-Stealth[scale=0.6]}] (6c2) to (4d2);
\draw[arrows = {-Stealth[scale=0.6]}] (7c2) to (2d2);
\draw[arrows = {-Stealth[scale=0.6]}] (7c2) to (4d2);
\draw[arrows = {-Stealth[scale=0.6]}] (8c2) to (3d2);
\draw[arrows = {-Stealth[scale=0.6]}] (8c2) to (4d2);
\draw[arrows = {-Stealth[scale=0.6]}] (1d2) to (1e2);
\draw[arrows = {-Stealth[scale=0.6]}] (2d2) to (1e2);
\draw[arrows = {-Stealth[scale=0.6]}] (3d2) to (1e2);
\draw[arrows = {-Stealth[scale=0.6]}] (4d2) to (1e2);
\end{tikzpicture}
\caption{The Loewy diagrams of $\mathcal{P}^+_{r,s}$, $\mathcal{P}^+_{r^\vee,s^\vee}$, $\mathcal{P}^-_{r^\vee,s}$ and $\mathcal{P}^-_{r,s^\vee}$, where $\triangle=L(h_{r,s})$, $\heartsuit=\mathcal{X}^+_{r,s}$, $\diamondsuit=\mathcal{X}^+_{r^\vee,s^\vee}$, $\clubsuit=\mathcal{X}^-_{r^\vee,s}$.\label{fig4.1}}
\end{figure}

%[htbp]
\begin{remark}
Figure \ref{fig4.1} is the embedding structure of the logarithmic $\W_{p_+,p_-}$-modules defined in Theorem~\ref{logarithmic1}. The structure of these logarithmic modules are studied by {\rm \cite{GRW0}},{\rm \cite{GRW}} in the case $(p_+,p_-)=(2,3)$.  In one block $C^{thick}_{1,1}$, explicit realizations are given by {\rm \cite{AK}}. 
\end{remark}

We define the following notation.
\begin{dfn} 
Let $M$ be a $\mathcal{W}_{p_+,p_-}$-module of socle length $n$, and let $0\leq i\leq n$.
For a simple module $X\subset {\rm S}_{n-i}(M)/{\rm S}_{n-i-1}(M)$, let $[X:{\rm S}_{n-i}(M)/{\rm S}_{n-i-1}(M)]$ be the multiplicity of $X$ in ${\rm S}_{n-i}(M)/{\rm S}_{n-i-1}(M)$. Then we denote by $(M:X:i)$ the submodule of ${\rm S}_{n-i}(M)$ satisfying
\begin{align*}
(M:X:i)/{\rm S}_{n-i-1}(M)\simeq [X:{\rm S}_{n-i}(M)/{\rm S}_{n-i-1}(M)]X, 
\end{align*}
where in the right-hand side, we have used the notation defined in Definition \ref{not0519}.
\end{dfn}

\begin{dfn}
\label{P^+ud}
Let $(a,b)$ be $(r,s)$ or $(r^\vee,s^\vee)$. 
We define the indecomposable modules 
\begin{align*}
&\mathcal{P}^{+d}_{a,b}:=(\mathcal{P}^+_{a,b}: \X^+_{a,b}:2),
&\mathcal{P}^{+u}_{a,b}:=\mathcal{P}^+_{a,b}/(\mathcal{P}^+_{a,b}: \X^+_{a^\vee,b^\vee}:2).
\end{align*}
\end{dfn}

\begin{figure}[htbp]
  \centering
\begin{tikzpicture}[scale=1.05]
\node[inner sep=0.5pt,scale=0.8] (a0) at (0, 0.43) {$\mathcal{P}^{+u}_{r,s}$};
\node[inner sep=0.5pt,scale=0.8] (a1) at (0, 0) {${\heartsuit}$};
\node[inner sep=0.5pt,scale=0.8] (b1) at (-1.5,-0.8) {$\clubsuit$};
\node[inner sep=0.5pt,scale=0.8] (b2) at (-0.75,-0.8) {$\clubsuit$};
\node[inner sep=0.5pt,scale=0.8] (b3) at (0,-0.8) {$\triangle$};
\node[inner sep=0.5pt,scale=0.8] (b4) at (0.75,-0.8) {$\spadesuit$};
\node[inner sep=0.5pt,scale=0.8] (b5) at (1.5,-0.8) {$\spadesuit$};
\node[inner sep=0.5pt,scale=0.8] (c1) at (-1.5-0.375,-1.6) {$\heartsuit$};
\node[inner sep=0.5pt,scale=0.8] (c6) at (1.5+0.375,-1.6) {$\heartsuit$};

\draw[arrows = {-Stealth[scale=0.6]}] (a1) to (b1);
\draw[arrows = {-Stealth[scale=0.6]}] (a1) to (b2);
\draw[arrows = {-Stealth[scale=0.6]}] (a1) to (b3);
\draw[arrows = {-Stealth[scale=0.6]}] (a1) to (b4);
\draw[arrows = {-Stealth[scale=0.6]}] (a1) to (b5);
\draw[arrows = {-Stealth[scale=0.6]}] (b1) to (c1);
\draw[arrows = {-Stealth[scale=0.6]}] (b2) to (c1);
\draw[arrows = {-Stealth[scale=0.6]}] (b3) to (c1);
\draw[arrows = {-Stealth[scale=0.6]}] (b3) to (c6);
\draw[arrows = {-Stealth[scale=0.6]}] (b4) to (c6);
\draw[arrows = {-Stealth[scale=0.6]}] (b5) to (c6);

%%第二列
\node[inner sep=0.5pt,scale=0.8] (a02) at (0, 0.43-3.5) {$\mathcal{P}^{+u}_{r^\vee,s^\vee}$};
\node[inner sep=0.5pt,scale=0.8] (a12) at (0, -3.5) {${\diamondsuit}$};
\node[inner sep=0.5pt,scale=0.8] (b12) at (-1.5,-0.8-3.5) {$\spadesuit$};
\node[inner sep=0.5pt,scale=0.8] (b22) at (-0.75,-0.8-3.5) {$\spadesuit$};
\node[inner sep=0.5pt,scale=0.8] (b32) at (0,-0.8-3.5) {$\triangle$};
\node[inner sep=0.5pt,scale=0.8] (b42) at (0.75,-0.8-3.5) {$\clubsuit$};
\node[inner sep=0.5pt,scale=0.8] (b52) at (1.5,-0.8-3.5) {$\clubsuit$};
\node[inner sep=0.5pt,scale=0.8] (c12) at (-1.5-0.375,-1.6-3.5) {$\diamondsuit$};
\node[inner sep=0.5pt,scale=0.8] (c62) at (1.5+0.375,-1.6-3.5) {$\diamondsuit$};

\draw[arrows = {-Stealth[scale=0.6]}] (a12) to (b12);
\draw[arrows = {-Stealth[scale=0.6]}] (a12) to (b22);
\draw[arrows = {-Stealth[scale=0.6]}] (a12) to (b32);
\draw[arrows = {-Stealth[scale=0.6]}] (a12) to (b42);
\draw[arrows = {-Stealth[scale=0.6]}] (a12) to (b52);
\draw[arrows = {-Stealth[scale=0.6]}] (b12) to (c12);
\draw[arrows = {-Stealth[scale=0.6]}] (b22) to (c12);
\draw[arrows = {-Stealth[scale=0.6]}] (b32) to (c12);
\draw[arrows = {-Stealth[scale=0.6]}] (b32) to (c62);
\draw[arrows = {-Stealth[scale=0.6]}] (b42) to (c62);
\draw[arrows = {-Stealth[scale=0.6]}] (b52) to (c62);

\node[inner sep=0.5pt,scale=0.8] (a03) at (5, 0.43) {$\mathcal{P}^{+d}_{r,s}$};
\node[inner sep=0.5pt,scale=0.8] (c13) at (-1.5-0.375+5,-0.2) {$\heartsuit$};
\node[inner sep=0.5pt,scale=0.8] (c63) at (1.5+0.375+5,-0.2) {$\heartsuit$};
\node[inner sep=0.5pt,scale=0.8] (d13) at (-1.5+5,-2.4+1.6-0.2) {$\spadesuit$};
\node[inner sep=0.5pt,scale=0.8] (d23) at (-0.75+5,-2.4+1.6-0.2) {$\spadesuit$};
\node[inner sep=0.5pt,scale=0.8] (d33) at (0+5,-2.4+1.6-0.2) {$\triangle$};
\node[inner sep=0.5pt,scale=0.8] (d43) at (0.75+5,-2.4+1.6-0.2) {$\clubsuit$};
\node[inner sep=0.5pt,scale=0.8] (d53) at (1.5+5,-2.4+1.6-0.2) {$\clubsuit$};
\node[inner sep=0.5pt,scale=0.8] (e13) at (0+5, -3.2+1.6-0.2) {$\heartsuit$};

\draw[arrows = {-Stealth[scale=0.6]}] (c13) to (d13);
\draw[arrows = {-Stealth[scale=0.6]}] (c13) to (d23);
\draw[arrows = {-Stealth[scale=0.6]}] (c13) to (d33);
\draw[arrows = {-Stealth[scale=0.6]}] (c63) to (d33);
\draw[arrows = {-Stealth[scale=0.6]}] (c63) to (d43);
\draw[arrows = {-Stealth[scale=0.6]}] (c63) to (d53);
\draw[arrows = {-Stealth[scale=0.6]}] (d13) to (e13);
\draw[arrows = {-Stealth[scale=0.6]}] (d23) to (e13);
\draw[arrows = {-Stealth[scale=0.6]}] (d33) to (e13);
\draw[arrows = {-Stealth[scale=0.6]}] (d43) to (e13);
\draw[arrows = {-Stealth[scale=0.6]}] (d53) to (e13);

\node[inner sep=0.5pt,scale=0.8] (a023) at (5, 0.43-3.5) {$\mathcal{P}^{+d}_{r^\vee,s^\vee}$};
\node[inner sep=0.5pt,scale=0.8] (c123) at (-1.5-0.375+5,-1.6-3.5+1.4) {$\diamondsuit$};
\node[inner sep=0.5pt,scale=0.8] (c623) at (1.5+0.375+5,-1.6-3.5+1.4) {$\diamondsuit$};
\node[inner sep=0.5pt,scale=0.8] (d123) at (-1.5+5,-2.4-3.5+1.4) {$\clubsuit$};
\node[inner sep=0.5pt,scale=0.8] (d223) at (-0.75+5,-2.4-3.5+1.4) {$\clubsuit$};
\node[inner sep=0.5pt,scale=0.8] (d323) at (0+5,-2.4-3.5+1.4) {$\triangle$};
\node[inner sep=0.5pt,scale=0.8] (d423) at (0.75+5,-2.4-3.5+1.4) {$\spadesuit$};
\node[inner sep=0.5pt,scale=0.8] (d523) at (1.5+5,-2.4-3.5+1.4) {$\spadesuit$};
\node[inner sep=0.5pt,scale=0.8] (e123) at (0+5, -3.2-3.5+1.4) {$\diamondsuit$};

\draw[arrows = {-Stealth[scale=0.6]}] (c123) to (d123);
\draw[arrows = {-Stealth[scale=0.6]}] (c123) to (d223);
\draw[arrows = {-Stealth[scale=0.6]}] (c123) to (d323);
\draw[arrows = {-Stealth[scale=0.6]}] (c623) to (d323);
\draw[arrows = {-Stealth[scale=0.6]}] (c623) to (d423);
\draw[arrows = {-Stealth[scale=0.6]}] (c623) to (d523);
\draw[arrows = {-Stealth[scale=0.6]}] (d123) to (e123);
\draw[arrows = {-Stealth[scale=0.6]}] (d223) to (e123);
\draw[arrows = {-Stealth[scale=0.6]}] (d323) to (e123);
\draw[arrows = {-Stealth[scale=0.6]}] (d423) to (e123);
\draw[arrows = {-Stealth[scale=0.6]}] (d523) to (e123);
\end{tikzpicture}
\caption{The Loewy diagrams of $\mathcal{P}^{+u}_{\bullet,\bullet}$ and $\mathcal{P}^{+d}_{\bullet,\bullet}$. The triangle $\bigtriangleup$ corresponds to the simple module $L(h_{r,s})$, $\heartsuit$ to $\X^+_{r,s}$, $\diamondsuit$ to $\X^+_{r^\vee,s^\vee}$, $\spadesuit$ to $\X^-_{r,s^\vee}$ and $\clubsuit$ to $\X^-_{r^\vee,s}$. \label{fig65}}
\end{figure}

Figure \ref{fig65} represents the embedding structure of the logarithmic $\W_{p_+,p_-}$-modules defined in Definition \ref{P^+ud}.

\begin{prop}
\label{P^d+}
We have
$
{\rm Ext}^1(\PP^{+d}_{r,s},\X^+_{r,s})={\rm Ext}^1(\PP^{+d}_{r^\vee,s^\vee},\X^+_{r^\vee,s^\vee})=0.
$
\begin{proof}
By Proposition \ref{quotient}, we see that $\PP^{+d}_{r,s}$ has $\mathcal{Q}(\X^+_{r,s})_{r,s^\vee}$ as a submodule. Then by the exact sequence 
\begin{align*}
0\rightarrow \mathcal{Q}(\X^+_{r,s})_{r,s^\vee}\rightarrow \PP^{+d}_{r,s}\rightarrow \mathcal{Q}(\X^+_{r,s})_{r^\vee,s}/\K_{r,s}\rightarrow 0
\end{align*}
and by Proposition \ref{Ext2+}, we obtain
$
{\rm Ext}^1(\PP^{+d}_{r,s},\X^+_{r,s})=0.
$
The second case can be proved in the same way, so we omit the proofs.
\end{proof}
\end{prop}

As in Lemma~\ref{bt3lem}, from the structure of $\mathcal{E}^\pm(\X^-_{\bullet,\bullet})_{\bullet,\bullet}$, we obtain the following lemma.
\begin{lem}
\label{bt3lem2}
Let $(a,b,c,d)$ be any element in
\begin{align*}
\{(r^\vee,s,r,s),(r^\vee,s,r^\vee,s^\vee),(r,s^\vee,r,s),(r,s^\vee,r^\vee,s^\vee)\},
\end{align*}
and fix an indecomposable module $E$ in ${\rm Ext}^1(\X^-_{a,b},\X^+_{c,d})$. For the surjection $\pi:E\rightarrow \X^-_{a,b}$, let $(v^-,v^+)\in E[\Delta^-_{a,b;0}]^2$ be $L_0$-homogeneous vectors of $E$ satisfying
\begin{align*}
&W^\pm[0]\pi(v^\pm)=0,
&W^\pm[0]\pi(v^\mp)\in \C^\times \pi(v^\pm).
\end{align*}
Then, for $\epsilon=\pm$, the Virasoro module $U(\mathcal{L}).v^\epsilon$ has a quotient isomorphic to
$
L(\Delta^-_{a,b;0},\Delta^+_{c,d;1}).
$
\end{lem}
\begin{prop}
\label{P^u+2}
We have
\begin{align*}
{\rm Ext}^1(\PP^{+u}_{r,s}/(2\X^+_{r,s}),\X^+_{r,s})={\rm Ext}^1(\PP^{+u}_{r^\vee,s^\vee}/(2\X^+_{r^\vee,s^\vee}),\X^+_{r^\vee,s^\vee})=\C^2.
\end{align*}
\begin{proof}
Let us prove the first equality. The second equality can be proved in the same way.
Since ${\rm Ext}^1(\X^+_{r,s},\X^+_{r,s})=0$, it is sufficient to show that
\begin{align*}
{\rm Ext}^1(\PP^{+u}_{r,s},\X^+_{r,s})=0.
\end{align*}
Assume ${\rm Ext}^1(\PP^{+u}_{r,s},\X^+_{r,s})\neq 0$ and fix a non-trivial extension $[E]$ in this ${\rm Ext}^1$-group.

First let us introduce some symbols. 
Let $t$ be a nonzero vector of the one dimensional space $E[h_{r,s}]$.
Let ${u}^{(0)}_0$ and ${u}^{(1)}_{0}$ be the Virasoro singular vectors of $3L(\Delta^+_{r,s;1})\subset \X^+_{r,s}$ defined in Definition \ref{dfn0830}.
Let $\iota_1$, $\iota_2$ and $\iota_3$ be injections from $\X^+_{r,s}$ to $E$ such that
\begin{align*}
\iota_1(\X^+_{r,s})\oplus \iota_2(\X^+_{r,s})\oplus \iota_3(\X^+_{r,s})={\rm S}(E).
\end{align*}
For the surjection $\pi :E\rightarrow \X^+_{r,s}$, we fix $L_0$-homogeneous vectors $\tilde{u}_0,\tilde{w}_0\in E$ such that
\begin{align*}
&\pi(\tilde{u}^{(0)}_0)=u^{(0)}_0,
&\pi(\tilde{u}^{(1)}_0)=u^{(1)}_0.
\end{align*}
Set $v^\pm_{r^\vee,s}=W^\pm[0]W^\mp[0]S_{r,s^\vee+p_-}\tilde{u}^{(0)}_0$ and $v^\pm_{r^\vee,s}=W^\pm[0]W^\mp[0]S_{r^\vee+p_+,s}\tilde{u}^{(0)}_0$.
Note that the sets
\begin{align*}
&\{v^-_{r^\vee,s},W^+[0]v^-_{r^\vee,s}\},
&\{W^-[0]v^+_{r^\vee,s},v^+_{r^\vee,s}\}
\end{align*}
correspond to bases of the lowest weight spaces of the composition factors $\X^-_{r^\vee,s}\oplus \X^-_{r^\vee,s}$, and 
\begin{align*}
&\{v^-_{r,s^\vee},W^+[0]v^-_{r,s^\vee}\},
&\{W^-[0]v^+_{r,s^\vee},v^+_{r,s^\vee}\}
\end{align*}
to bases of the lowest weight spaces of the composition factors $\X^-_{r,s^\vee}\oplus \X^-_{r,s^\vee}$.

Consider the Virasoro module
\begin{align*}
M=U(\mathcal{L}).\tilde{u}^{(1)}_0+\sum_{\epsilon=\pm}U(\mathcal{L}).v^\epsilon_{r^\vee,s}+\sum_{\epsilon'=\pm}U(\mathcal{L}).v^{\epsilon'}_{r,s^\vee}.
\end{align*}
As in the argument in the proof of Proposition \ref{Ext2+0}, from the Virasoro structure of $\PP^{+u}_{r,s}$, we see that the Virasoro module $M$ has an indecomposable quotient isomorphic to $\widetilde{P}^{\downarrow}(\Delta^+_{r,s;1})$, and the quotient module
\begin{equation*}
M/\bigl(\sum_{\epsilon=\pm}U(\mathcal{L}).v^\epsilon_{r^\vee,s}+\sum_{\epsilon'=\pm}U(\mathcal{L}).v^{\epsilon'}_{r,s^\vee}\bigr)
\end{equation*}
does not contain $L(\Delta^-_{r^\vee,s;0})$ and $L(\Delta^-_{r,s^\vee;0})$ as the composition factors.
Thus by Proposition \ref{bt3}, we have
\begin{equation}
\label{bt17}
\begin{split}
&\sigma(S_{r^\vee,s^\vee+2p_-})\tilde{u}^{(1)}_0\in \C^\times v^+_{r^\vee,s}+\C^\times v^-_{r^\vee,s}+\sum _{i=1}^3U(\mathcal{L}).\iota_i(u^{(0)}_0)+U(\mathcal{L}).t,\\
&\sigma(S_{r^\vee+2p_+,s^\vee})\tilde{u}^{(1)}_0\in \C^\times v^+_{r,s^\vee}+\C^\times v^-_{r,s^\vee}+\sum _{i=1}^3U(\mathcal{L}).\iota_i(u^{(0)}_0)+U(\mathcal{L}).t.
\end{split}
\end{equation} 
As in the argument in the proof of Proposition \ref{Ext2+0}, from the structure of $\widetilde{P}(\Delta^+_{r,s;1})$, we have 
\begin{equation}
\label{ntt007}
\begin{split}
&S_{r^\vee,s^\vee+2p_-}\sigma(S_{r^\vee,s^\vee+2p_-})\tilde{u}^{(1)}_0\equiv \sum_{i=1}^3k_i\iota_i(u^{(1)}_0)\ {\rm mod}\ \sum _{j=1}^3U(\mathcal{L}).\iota_j(u^{(0)}_0)+U(\mathcal{L}).t,\\
&S_{r^\vee+2p_+,s^\vee}\sigma(S_{r^\vee+2p_+,s^\vee})\tilde{u}^{(1)}_0\equiv \sum_{i=1}^3l_i\iota_i(u^{(1)}_0)\ {\rm mod}\ \sum _{j=1}^3U(\mathcal{L}).\iota_j(u^{(0)}_0)+U(\mathcal{L}).t,
\end{split}
\end{equation}
where $(k_i,l_i)\ (i=1,2,3)$ are constants such that not all of the $k_i$ are zero and not all of the $l_i$ are zero.
Since ${\rm Ext}^1(\X^+_{r,s},\X^+_{r,s})=0$, we have
\begin{align}
\label{ntt0516}
\begin{aligned}
&\W_{p_+,p_-}.\bigl(\sum_{i=1}^3k_i\iota_i(u^{(1)}_0)\bigr)\simeq \X^+_{r,s},
&\W_{p_+,p_-}.\bigl(\sum_{i=1}^3l_i\iota_i(u^{(1)}_0)\bigr)\simeq \X^+_{r,s}.
\end{aligned}
\end{align}
We set
\begin{align*}
E'=\frac{E}{\W_{p_+,p_-}.\bigl(\sum_{i=1}^3k_i\iota_i(u^{(1)}_0)\bigr)\oplus \W_{p_+,p_-}.\bigl(\sum_{i=1}^3l_i\iota_i(u^{(1)}_0)\bigr)}.
\end{align*}
Since $[E]\neq 0$, from (\ref{ntt0516}), we have
\begin{align*}
[E']\in {\rm Ext}^1(\mathcal{P}^{+u}_{r,s}/(\X^+_{r,s}\oplus \X^+_{r,s}),\X^+_{r,s})\setminus\{0\}.
\end{align*}
Let $\phi$ be the surjection from $E$ to $E'$. 
Then, from (\ref{ntt007}), we see that
\begin{equation}
\label{bt27}
\begin{split}
&S_{r^\vee,s^\vee+2p_-}\sigma(S_{r^\vee,s^\vee+2p_-})\phi(\tilde{u}^{(1)}_0)\in \sum _{i=1}^3U(\mathcal{L}).\phi\circ\iota_i(u^{(0)}_0)+U(\mathcal{L}).\phi (t),\\
&S_{r^\vee+2p_+,s^\vee}\sigma(S_{r^\vee+2p_+,s^\vee})\phi(\tilde{u}^{(1)}_0)\in \sum _{i=1}^3U(\mathcal{L}).\phi\circ\iota_i(u^{(0)}_0)+U(\mathcal{L}).\phi (t).
\end{split}
\end{equation}
From Proposition~\ref{bt3}, (\ref{bt17}), (\ref{ntt007}) and (\ref{bt27}), we have
\begin{equation}
\label{bt377}
\begin{split}
&S_{r^\vee,s^\vee+2p_-}\phi(v^\pm_{r^\vee,s})\in \sum _{i=1}^3U(\mathcal{L}).\phi\circ\iota_i(u^{(0)}_0)+U(\mathcal{L}).\phi (t),\\
&S_{r^\vee+2p_+,s^\vee}\phi(v^\pm_{r,s^\vee})\in \sum _{i=1}^3U(\mathcal{L}).\phi\circ\iota_i(u^{(0)}_0)+U(\mathcal{L}).\phi (t).
\end{split}
\end{equation}
Then, by Lemma \ref{bt3lem2} and (\ref{bt377}), $E'$ has $2\X^-_{r^\vee,s}\oplus 2\X^-_{r,s^\vee}$ as submodules. Thus $E'$ has an indecomposable quotient in ${\rm Ext}^1(\mathcal{K}^*_{r,s},\X^+_{r,s})$. But this contradicts Proposition \ref{Ext2+}.
\end{proof}
\end{prop}

\begin{prop}
\label{P1}
For $\mathcal{P}^+_{r,s},\ \mathcal{P}^+_{r^\vee,s^\vee}\in C^{thick}_{r,s}$, we have
\begin{align*}
{\rm Ext}^1(\mathcal{P}^+_{r,s},\X^+_{r,s})={\rm Ext}^1(\mathcal{P}^+_{r^\vee,s^\vee},\X^+_{r^\vee,s^\vee})=0.
\end{align*}
\begin{proof}
By the exact sequence
$
0\rightarrow \PP^{+d}_{r,s}\rightarrow \PP^+_{r,s} \rightarrow \PP^+_{r,s}/\PP^{+d}_{r,s}\rightarrow 0
$
and Proposition \ref{P^d+}, we have the following exact sequence
\begin{align*}
0\rightarrow\C\rightarrow\C\rightarrow \C^2\rightarrow{\rm Ext}^1(\PP^+_{r,s}/\PP^{+d}_{r,s},\X^+_{r,s})\rightarrow{\rm Ext}^1(\PP^+_{r,s},\X^+_{r,s})\rightarrow 0.
\end{align*}
By Proposition \ref{P^u+2}, we have
$
{\rm Ext}^1(\PP^+_{r,s}/\PP^{+d}_{r,s},\X^+_{r,s})=\C^2.
$
Thus, from the above exact sequence, we have 
$
{\rm Ext}^1(\PP^+_{r,s},\X^+_{r,s})=0.
$
The second equality can be proved in the same way.
\end{proof}
\end{prop}

\begin{dfn}
Let $(a,b)$ be $(r^\vee,s)$ or $(r,s^\vee)$. We define the quotient module
\begin{align*}
\PP^{-u}_{a,b}:=\frac{\mathcal{P}^-_{a,b}}{(\mathcal{P}^-_{a,b}:\X^-_{a^\vee,b^\vee}:2)+(\mathcal{P}^-_{a,b}:L(h_{r,s}):2)}.
\end{align*}
\label{P^-ud}
\end{dfn}

\begin{figure}[htbp]
  \centering
\begin{tikzpicture}[scale=1.05]
\node[inner sep=0.5pt,scale=0.8] (0a) at (0, 0.43) {$\mathcal{P}^{-u}_{r,s^\vee}$};
\node[inner sep=0.5pt,scale=0.8] (1a) at (0, 0) {${\spadesuit}$};
\node[inner sep=0.5pt,scale=0.8] (1b) at (0-1.75,-0.8) {$\diamondsuit$};
\node[inner sep=0.5pt,scale=0.8] (2b) at (0-1.1,-0.8) {$\diamondsuit$};
\node[inner sep=0.5pt,scale=0.8] (3b) at (0+1.1,-0.8) {$\heartsuit$};
\node[inner sep=0.5pt,scale=0.8] (4b) at (0+1.75,-0.8) {$\heartsuit$};
\node[inner sep=0.5pt,scale=0.8] (1c) at (0-1.75+5.5-1.75-5.5+1.1,-1.6) {$\spadesuit$};
\node[inner sep=0.5pt,scale=0.8] (8c) at (0+1.75+1.75-1.1,-1.6) {$\spadesuit$};

\draw[arrows = {-Stealth[scale=0.6]}] (1a) to (1b);
\draw[arrows = {-Stealth[scale=0.6]}] (1a) to (2b);
\draw[arrows = {-Stealth[scale=0.6]}] (1a) to (3b);
\draw[arrows = {-Stealth[scale=0.6]}] (1a) to (4b);
\draw[arrows = {-Stealth[scale=0.6]}] (1b) to (1c);
\draw[arrows = {-Stealth[scale=0.6]}] (2b) to (1c);
\draw[arrows = {-Stealth[scale=0.6]}] (3b) to (8c);
\draw[arrows = {-Stealth[scale=0.6]}] (4b) to (8c);

\node[inner sep=0.5pt,scale=0.8] (0a2) at (5.7, 0.43) {$\mathcal{P}^{-u}_{r^\vee,s}$};
\node[inner sep=0.5pt,scale=0.8] (1a2) at (5.7, 0) {${\clubsuit}$};
\node[inner sep=0.5pt,scale=0.8] (1b2) at (5.7-1.75,-0.8) {$\heartsuit$};
\node[inner sep=0.5pt,scale=0.8] (2b2) at (5.7-1.1,-0.8) {$\heartsuit$};
\node[inner sep=0.5pt,scale=0.8] (3b2) at (5.7+1.1,-0.8) {$\diamondsuit$};
\node[inner sep=0.5pt,scale=0.8] (4b2) at (5.7+1.75,-0.8) {$\diamondsuit$};
\node[inner sep=0.5pt,scale=0.8] (1c2) at (5.7-1.75+5.5-1.75-5.5+1.1,-1.6) {$\clubsuit$};
\node[inner sep=0.5pt,scale=0.8] (8c2) at (5.7+1.75+1.75-1.1,-1.6) {$\clubsuit$};

\draw[arrows = {-Stealth[scale=0.6]}] (1a2) to (1b2);
\draw[arrows = {-Stealth[scale=0.6]}] (1a2) to (2b2);
\draw[arrows = {-Stealth[scale=0.6]}] (1a2) to (3b2);
\draw[arrows = {-Stealth[scale=0.6]}] (1a2) to (4b2);
\draw[arrows = {-Stealth[scale=0.6]}] (1b2) to (1c2);
\draw[arrows = {-Stealth[scale=0.6]}] (2b2) to (1c2);
\draw[arrows = {-Stealth[scale=0.6]}] (3b2) to (8c2);
\draw[arrows = {-Stealth[scale=0.6]}] (4b2) to (8c2);

\end{tikzpicture}
\caption{The Loewy diagrams of the $\W_{p_+,p_-}$-modules in Definition \ref{P^-ud}. The $\heartsuit$ corresponds to the simple module $\X^+_{r,s}$, $\diamondsuit$ to $\X^+_{r^\vee,s^\vee}$, $\spadesuit$ to $\X^-_{r,s^\vee}$ and $\clubsuit$ to $\X^-_{r^\vee,s}$.\label{fig20230302}}
\end{figure}

Figure \ref{fig20230302} represents the embedding structure of the quotient modules $\mathcal{P}^{-u}_{r,s^\vee}$ and $\mathcal{P}^{-u}_{r^\vee,s}$ defined in Definition \ref{P^-ud}.

\begin{prop}
\label{2-}
We have
\begin{align*}
{\rm Ext}^1(\mathcal{P}^{-u}_{r^\vee,s}/(2\X^-_{r^\vee,s}),\X^-_{r^\vee,s})={\rm Ext}^1(\mathcal{P}^{-u}_{r,s^\vee}/(2\X^-_{r,s^\vee}),\X^-_{r,s^\vee})=\C^2.
\end{align*}
\begin{proof}
Let us prove the first equality. The second equality can be proved in the same way.
Since ${\rm Ext}^1(\X^-_{r^\vee,s},\X^-_{r^\vee,s})=0$, it is sufficient to show that
\begin{align*}
{\rm Ext}^1(\mathcal{P}^{-u}_{r^\vee,s},\X^-_{r^\vee,s})=0.
\end{align*}
Assume ${\rm Ext}^1(\mathcal{P}^{-u}_{r^\vee,s},\X^-_{r^\vee,s})\neq 0$. Fix a non-trivial extension
\begin{align*}
0\rightarrow \X^-_{r^\vee,s}\xrightarrow{\iota}E\xrightarrow{\phi} \mathcal{P}^{-u}_{r^\vee,s}\rightarrow 0.
\end{align*}
By Propositions \ref{Ext} and \ref{Ext2-+}, we see that
$
{\rm S}(E)=\iota(\X^-_{r^\vee,s})\oplus 2\X^-_{r^\vee,s}.
$
Let $\iota_1$ and $\iota_2$ be injections from $\X^-_{r^\vee,s}$ to $E$ such that
\begin{align*}
\phi\circ\iota_1(\X^-_{r^\vee,s})\oplus \phi\circ\iota_2(\X^-_{r^\vee,s})=\X^-_{r^\vee,s}\oplus \X^-_{r^\vee,s}={\rm S}(\mathcal{P}^{-u}_{r^\vee,s}),
\end{align*}
and let $\{v_+,v_-\}$ be a basis of the lowest weight space of $\X^-_{r^\vee,s}$ such that
\begin{align*}
&W^\pm[0]v_\pm=0,
&W^\pm [0]v_\mp\in \C^\times v_\pm.
\end{align*}
For the canonical surjection $\pi : E\rightarrow \X^-_{r^\vee,s}$, we fix $L_0$-homogeneous vectors $\tilde{v}_-,\tilde{v}_+\in E$ such that $\pi(\tilde{v}_{\pm})=v_{\pm}$.
Note that, by Proposition \ref{sankakuthick2}, the Virasoro module $U(\mathcal{L}).\tilde{v}_\pm$ has $\widetilde{P}(\Delta^-_{r^\vee,s;0})$ as the quotient (for the definition of $\widetilde{P}(\Delta^-_{r^\vee,s;0})$, see the proof of \ref{sankakuthick2}). 
Let
\begin{align*}
&E_1=(E:\X^+_{r,s}:1),
&E_2=(E:\X^+_{r^\vee,s^\vee}:1).
\end{align*}
Then, by Proposition \ref{Ext2-+}, we see that $E_1$ and $E_2$ must contain $\iota(\X^-_{r^\vee,s})$ as the submodule:
\begin{align}
\label{P^-0424}
(\iota(\X^-_{r^\vee,s})\subset E_1)\land (\iota(\X^-_{r^\vee,s})\subset E_2).
\end{align}
Noting the structure of $\widetilde{P}(\Delta^-_{r^\vee,s;0})$, $\widetilde{\mathcal{E}}(\X^+_{r,s})_{r^\vee,s}$ and $\widetilde{\mathcal{E}}(\X^+_{r^\vee,s^\vee})_{r^\vee,s}$, from (\ref{P^-0424}),
we have
\begin{equation}
\label{P^-04241}
\begin{split}
&(L_0-\Delta^-_{r^\vee,s;0})\tilde{v}_+\in \C^\times \iota(v_-)+\C\iota(v_+)+\C^\times\iota_1(v_+)+\C^\times\iota_2(v_+),\\
&(L_0-\Delta^-_{r^\vee,s;0})\tilde{v}_-\in \C^\times \iota(v_+)+\C\iota(v_-)+\C^\times\iota_1(v_-)+\C^\times\iota_2(v_-).
\end{split}
\end{equation}
Multiplying $(L_0-\Delta^-_{r^\vee,s;0})\tilde{v}_\pm$ by $W^\pm[0]$, from (\ref{P^-04241}), we have $\iota(v_\pm)=0$. But this is a contradiction.
\end{proof}
\end{prop}
\begin{prop}
\label{P2}
We have
$
{\rm Ext}^1(\PP^-_{r^\vee,s},\X^-_{r^\vee,s})={\rm Ext}^1(\PP^-_{r,s^\vee},\X^-_{r,s^\vee})=0.
$
\begin{proof}
By the exact sequence
\begin{align*}
0\rightarrow \mathcal{Q}(\X^-_{r^\vee,s})_{r^\vee,s^\vee}\rightarrow \PP^-_{r^\vee,s} \rightarrow\PP^-_{r^\vee,s}/\mathcal{Q}(\X^-_{r^\vee,s})_{r^\vee,s^\vee}\rightarrow 0
\end{align*}
and Proposition \ref{Ext2-+}, we have the following exact sequence
\begin{align}
0\rightarrow\C\rightarrow\C\rightarrow\C\rightarrow{\rm Ext}^1(\PP^-_{r^\vee,s}/\mathcal{Q}(\X^-_{r^\vee,s})_{r^\vee,s^\vee},\X^-_{r^\vee,s})\rightarrow{\rm Ext}^1(\PP^-_{r^\vee,s},\X^-_{r^\vee,s})\rightarrow 0.
\label{kuso}
\end{align}
Note that $\PP^-_{r^\vee,s}/\mathcal{Q}(\X^-_{r^\vee,s})_{r^\vee,s^\vee} $ has a submodule isomorphic to $\K_{r,s}\oplus\K_{r,s}$. Let $\mathcal{M}$ be the quotient module defined by the following exact sequence
\begin{align*}
0\rightarrow \K_{r,s}\oplus\K_{r,s}\rightarrow \PP^-_{r^\vee,s}/\mathcal{Q}(\X^-_{r^\vee,s})_{r^\vee,s^\vee} \rightarrow \mathcal{M}\rightarrow 0.
\end{align*}
By this exact sequence and by ${\rm Ext}^1(\K_{r,s},\X^-_{r^\vee,s})=0$ (see the proof of Proposition \ref{Extcoro1}), we have
\begin{align}
\label{20221026M}
{\rm Ext}^1(\mathcal{M},\X^-_{r^\vee,s})\simeq{\rm Ext}^1(\PP^-_{r^\vee,s}/\mathcal{Q}(\X^-_{r^\vee,s})_{r^\vee,s^\vee},\X^-_{r^\vee,s}).
\end{align}
Note that $\mathcal{M}$ satisfies the exact sequence
\begin{align*}
0\rightarrow \X^-_{r^\vee,s}\oplus 4\X^-_{r,s^\vee}\rightarrow \mathcal{M} \rightarrow \PP^{-u}_{r^\vee,s}/(2\X^-_{r^\vee,s})\rightarrow 0.
\end{align*}
By this exact sequence and Proposition \ref{2-}, we have the following exact sequence
\begin{equation*}
0\rightarrow \C\rightarrow\C\rightarrow\C\rightarrow\C^2\rightarrow {\rm Ext}^1(\mathcal{M},\X^-_{r^\vee,s})\rightarrow 0.
\end{equation*}
Thus, by this exact sequence, we have ${\rm Ext}^1(\mathcal{M},\X^-_{r^\vee,s})\simeq \C$. Therefore, by (\ref{kuso}) and (\ref{20221026M}), we obtain
$
{\rm Ext}^1(\PP^-_{r^\vee,s},\X^-_{r^\vee,s})=0.
$
The second equality can be proved in the same way.
\end{proof}
\end{prop}

Since all logarithmic modules $\mathcal{P}^\pm_{\bullet,\bullet}$ in $C^{thick}_{r,s}$ are generated from the top composition factors, by Propositions~\ref{P4}, \ref{P1}, and \ref{P2}, we obtain the following theorem.
\begin{thm}
\label{77h}
For each $(a,b,\epsilon)\in \mathbb{S}_1$, the indecomposable module $\mathcal{P}^\epsilon_{a,b}$ is the projective cover of $\X^\epsilon_{a,b}$.
\end{thm}

\subsection{Projective covers of the minimal modules $L(h_{r,s})$}
\label{Projminimal}
Fix an arbitrary thick block $C^{thick}_{r,s}$. Let $\mathcal{P}(h_{r,s})$ be the projective cover of the minimal simple module $L(h_{r,s})$. 
By Corollary \ref{prop202302110}, we see that $\mathcal{P}(h_{r,s})$ has $L_0$ nilpotent rank three.
In this subsection, we determine the structure of $\mathcal{P}(h_{r,s})$.

\begin{dfn}
Let $(a,b)$ be $(r,s)$ or $(r^\vee,s^\vee)$. We define 
\begin{align*}
\mathcal{N}_{a,b}:=(\mathcal{P}^+_{a,b}:L(h_{r,s}):1).
\end{align*}
\end{dfn}
For $(a,b)=(r,s)$ or $(r^\vee,s^\vee)$, $\mathcal{N}_{a,b}$ has the following lenghth four socle series
\begin{align*}
{\rm S}_{1}&=\X^+_{a,b}, 
&{\rm S}_{2}/{\rm S}_{1}&=2\X^-_{a^\vee,b}\oplus L(h_{r,s})\oplus 2\X^-_{a,b^\vee},\\
{\rm S}_3/{\rm S}_2&=\X^+_{a,b}\oplus \X^+_{a^\vee,b^\vee},
&\mathcal{N}_{a,b}/{\rm S}_3&=L(h_{r,s}).
\end{align*}

\vspace{3mm}
\begin{figure}[htbp]
  \centering
\begin{tikzpicture}[scale=1.05]
\node[inner sep=0.5pt] (a1) at (0, 0) {${\triangle}$};
\node[inner sep=0.5pt] (b1) at (-0.75,-1) {$\heartsuit$};
\node[inner sep=0.5pt] (b2) at (0.75,-1) {$\diamondsuit$};
\node[inner sep=0.5pt] (c1) at (-2,-2) {$\spadesuit$};
\node[inner sep=0.5pt] (c2) at (-1.0,-2) {$\spadesuit$};
\node[inner sep=0.5pt] (c3) at (0,-2) {$\triangle$};
\node[inner sep=0.5pt] (c4) at (1.0,-2) {$\clubsuit$};
\node[inner sep=0.5pt] (c5) at (2,-2) {$\clubsuit$};
\node[inner sep=0.5pt] (d2) at (0.75,-3) {$\heartsuit$};

\draw[arrows = {-Stealth[scale=0.7]}] (a1) to (b1);
\draw[arrows = {-Stealth[scale=0.7]}] (a1) to (b2);
\draw[arrows = {-Stealth[scale=0.7]}] (b1) to (c1);
\draw[arrows = {-Stealth[scale=0.7]}] (b1) to (c2);
\draw[arrows = {-Stealth[scale=0.7]}] (b1) to (c3);
\draw[arrows = {-Stealth[scale=0.7]}] (b1) to (c4);
\draw[arrows = {-Stealth[scale=0.7]}] (b1) to (c5);
\draw[arrows = {-Stealth[scale=0.7]}] (b2) to (c1);
\draw[arrows = {-Stealth[scale=0.7]}] (b2) to (c2);
\draw[arrows = {-Stealth[scale=0.7]}] (b2) to (c3);
\draw[arrows = {-Stealth[scale=0.7]}] (b2) to (c4);
\draw[arrows = {-Stealth[scale=0.7]}] (b2) to (c5);
\draw[arrows = {-Stealth[scale=0.7]}] (c1) to (d2);
\draw[arrows = {-Stealth[scale=0.7]}] (c2) to (d2);
\draw[arrows = {-Stealth[scale=0.7]}] (c3) to (d2);
\draw[arrows = {-Stealth[scale=0.7]}] (c4) to (d2);
\draw[arrows = {-Stealth[scale=0.7]}] (c5) to (d2);

\node[inner sep=0.5pt] (a12) at (5.5, 0) {${\triangle}$};
\node[inner sep=0.5pt] (b12) at (5.5-0.75,-1) {$\heartsuit$};
\node[inner sep=0.5pt] (b22) at (5.5+0.75,-1) {$\diamondsuit$};
\node[inner sep=0.5pt] (c12) at (5.5-2,-2) {$\spadesuit$};
\node[inner sep=0.5pt] (c22) at (5.5-1.0,-2) {$\spadesuit$};
\node[inner sep=0.5pt] (c32) at (5.5+0,-2) {$\triangle$};
\node[inner sep=0.5pt] (c42) at (5.5+1.0,-2) {$\clubsuit$};
\node[inner sep=0.5pt] (c52) at (5.5+2,-2) {$\clubsuit$};
\node[inner sep=0.5pt] (d12) at (5.5-0.75,-3) {$\diamondsuit$};

\draw[arrows = {-Stealth[scale=0.7]}] (a12) to (b12);
\draw[arrows = {-Stealth[scale=0.7]}] (a12) to (b22);
\draw[arrows = {-Stealth[scale=0.7]}] (b12) to (c12);
\draw[arrows = {-Stealth[scale=0.7]}] (b12) to (c22);
\draw[arrows = {-Stealth[scale=0.7]}] (b12) to (c32);
\draw[arrows = {-Stealth[scale=0.7]}] (b12) to (c42);
\draw[arrows = {-Stealth[scale=0.7]}] (b12) to (c52);
\draw[arrows = {-Stealth[scale=0.7]}] (b22) to (c12);
\draw[arrows = {-Stealth[scale=0.7]}] (b22) to (c22);
\draw[arrows = {-Stealth[scale=0.7]}] (b22) to (c32);
\draw[arrows = {-Stealth[scale=0.7]}] (b22) to (c42);
\draw[arrows = {-Stealth[scale=0.7]}] (b22) to (c52);
\draw[arrows = {-Stealth[scale=0.7]}] (c12) to (d12);
\draw[arrows = {-Stealth[scale=0.7]}] (c22) to (d12);
\draw[arrows = {-Stealth[scale=0.7]}] (c32) to (d12);
\draw[arrows = {-Stealth[scale=0.7]}] (c42) to (d12);
\draw[arrows = {-Stealth[scale=0.7]}] (c52) to (d12);
\end{tikzpicture}
\caption{The Loewy diagram of $\mathcal{N}_{r,s}$ and $\mathcal{N}_{r^\vee,s^\vee}$. The triangle $\bigtriangleup$ corresponds to the simple module $L(h_{r,s})$, $\heartsuit$ to $\X^+_{r,s}$, $\diamondsuit$ to $\X^+_{r^\vee,s^\vee}$, $\spadesuit$ to $\X^-_{r,s^\vee}$ and $\clubsuit$ to $\X^-_{r^\vee,s}$. \label{fig20230125p}}
\end{figure}

Figure \ref{fig20230125p} represents the embedding structure of the logarithmic $\W_{p_+,p_-}$-modules $\mathcal{N}_{r,s}$ and $\mathcal{N}_{r^\vee,s^\vee}$.

We define 
\begin{align*}
\mathcal{Q}(h_{r,s}):=\frac{\mathcal{N}_{r,s}}{(\mathcal{N}_{r,s}:\X^-_{r^\vee,s}:2)+(\mathcal{N}_{r,s}:\X^-_{r,s^\vee}:2)}.
\end{align*}
By Propositions \ref{Ext} and by the proof of Proposition \ref{Extcoro2}, we obtain the following lemmas.
\begin{lem}
The socle series of the indecomposable module $\mathcal{Q}(h_{r,s})$ is given by
\begin{align*}
{\rm S}_1=L(h_{r,s}),\ \ \ \ \ \ {\rm S}_2/{\rm S}_1=\X^+_{r,s}\oplus \X^+_{r^\vee,s^\vee},\ \ \ \ \ \ {\rm S}_3/{\rm S}_2=L(h_{r,s}).
\end{align*}
\end{lem}

\begin{lem}
\label{lem202301240}
We have
$
{\rm Ext}^1(\mathcal{Q}(h_{r,s}),L(h_{r,s}))=0.
$
\end{lem}
Recall the indecomposable modules $\mathcal{P}^{+u}_{r,s}$ and $\mathcal{P}^{+u}_{r^\vee,s^\vee}$ defined in Definition \ref{P^+ud}.
Note that $\mathcal{P}^{+u}_{r,s}$ and $\mathcal{P}^{+u}_{r^\vee,s^\vee}$ have $\K_{r,s}$ and $\K_{r^\vee,s^\vee}$ as submodules respectively.
We define the indecomposable modules 
\begin{align*}
&\mathcal{R}_{r,s}:=\mathcal{P}^{+u}_{r,s}/\K_{r,s},
&\mathcal{R}_{r^\vee,s^\vee}:=\mathcal{P}^{+u}_{r^\vee,s^\vee}/\K_{r^\vee,s^\vee}.
\end{align*}
By Proposition \ref{Y} and by Proposition \ref{Ext}, we obtain the following lemma.
\begin{lem}
\label{R_{r,s}}
Let $(a,b)$ be $(r,s)$ or $(r^\vee,s^\vee)$. Then the indecomposable module $\mathcal{R}_{a,b}$ has $L_0$ nilpotent rank two and has the following length three socle series
\begin{align*}
 {\rm S}_{1}=\X^+_{a,b},\quad {\rm S}_{2}/{\rm S}_{1}=2\X^-_{r^\vee,s}\oplus 2\X^-_{r,s^\vee},\quad \mathcal{R}_{a,b}/{\rm S}_2=\X^+_{a,b}.
\end{align*} 
\end{lem}
By Propositions \ref{Ext}, \ref{P^u+2}, we obtain the following lemma. 
\begin{lem}
\label{lem20230124}
We have
$
{\rm Ext}^1(\mathcal{R}_{r,s},\X^+_{r,s})={\rm Ext}^1(\mathcal{R}_{r^\vee,s^\vee},\X^+_{r^\vee,s^\vee})=0.
$
\end{lem}

\begin{lem}
\label{lem202301142}
We have
\begin{align*}
{\rm Ext}^1(\mathcal{N}_{r,s}/\K_{r,s},\X^+_{r,s})={\rm Ext}^1(\mathcal{N}_{r^\vee,s^\vee}/\K_{r^\vee,s^\vee},\X^+_{r^\vee,s^\vee})=0.
\end{align*}
\begin{proof}
We only prove ${\rm Ext}^1(\mathcal{N}_{r,s}/\K_{r,s},\X^+_{r,s})=0$. The second case can be proved in the same way.

Assume that ${\rm Ext}^1(\mathcal{N}_{r,s}/\K_{r,s},\X^+_{r,s})\neq 0$ and fix a non-trivial extension 
$[E]$ in this ${\rm Ext}^1$-group. Let $t\in E$ be the lowest weight vector in the one dimensional space $E[h_{r,s}]$.
By Lemma \ref{lem20230124}, we see that $E$ has a submodule isomorphic to $\mathcal{R}_{r,s}$. Note that by Lemma \ref{R_{r,s}}, $\mathcal{R}_{r,s}$ has $L_0$ nilpotent rank two. Then we have a basis $\{u_0,u_1\}$ of the lowest weight space of the submodule $\mathcal{R}_{r,s}\subset E$ satisfying
\begin{align}
\label{202303210}
(L_0-\Delta^+_{r,s;0})u_0\in \C^\times u_1.
\end{align}
On the other hand, by
\begin{align*}
&(L_0-h_{r,s})t=0,
&S_{r^\vee,s^\vee}t\in \C^\times u_0+\C u_1,
\end{align*}
we have $(L_0-\Delta^+_{r,s;0})u_0=0$. But this contradicts (\ref{202303210}).
\end{proof}
\end{lem}

\begin{lem}
We have
\begin{align}
\label{20230322Ext}
{\rm Ext}^1(\mathcal{N}_{r,s},\X^+_{r^\vee,s^\vee})=\C.
\end{align}
\begin{proof}
By the exact sequence 
$
0\rightarrow \K_{r,s}\rightarrow \mathcal{N}_{r,s}\rightarrow \mathcal{N}_{r,s}/\K_{r,s}\rightarrow 0,
$
we have the exact sequence 
\begin{align}
\label{202303211}
0\rightarrow {\rm Ext}^1(\mathcal{N}_{r,s}/\K_{r,s},\X^+_{r^\vee,s^\vee})\rightarrow {\rm Ext}^1(\mathcal{N}_{r,s},\X^+_{r^\vee,s^\vee})\rightarrow \C.
\end{align}
By Proposition \ref{Ext}, we have $\mathcal{N}_{r,s}/\K_{r,s}\simeq \mathcal{N}_{r^\vee,s^\vee}/\K_{r^\vee,s^\vee}$. Then, by Lemma \ref{lem202301142} and (\ref{202303211}), we obtain 
\begin{align}
\label{202303212}
{\rm dim}_{\C}{\rm Ext}^1(\mathcal{N}_{r,s},\X^+_{r^\vee,s^\vee})\leq 1.
\end{align}
Note that by Proposition \ref{Ext} and Lemma \ref{lem202301240}, 
$
\mathcal{N}_{r,s}/\X^+_{r,s}\simeq \mathcal{N}_{r^\vee,s^\vee}/\X^+_{r^\vee,s^\vee}.
$
Then, we have ${\rm Ext}^1(\mathcal{N}_{r,s},\X^+_{r^\vee,s^\vee})\neq 0$. Thus, by (\ref{202303212}), we obtain (\ref{20230322Ext}).
\end{proof}
\end{lem}
Let $\widetilde{\mathcal{N}}(h_{r,s})$ be an indecomposable module representing a class in the ${\rm Ext}^1$-group (\ref{20230322Ext}).
By Lemma \ref{lem202301142}, we obtain the following lemma.
\begin{lem}
\label{null20230322+}
We have
$
{\rm Ext}^1(\widetilde{\mathcal{N}}(h_{r,s}),\X^+_{r,s})={\rm Ext}^1(\widetilde{\mathcal{N}}(h_{r,s}),\X^+_{r^\vee,s^\vee})=0.
$
\end{lem}
\begin{lem}
\label{null20230322-}
We have
$
{\rm Ext}^1(\widetilde{\mathcal{N}}(h_{r,s}),\X^-_{r^\vee,s})={\rm Ext}^1(\widetilde{\mathcal{N}}(h_{r,s}),\X^-_{r,s^\vee})=0.
$
\begin{proof}
We only prove the first case. The second case can be proved in the same way. Assume that
$
{\rm Ext}^1(\widetilde{\mathcal{N}}(h_{r,s}),\X^-_{r^\vee,s})\neq 0
$
and fix a non-trivial extension $[E]$ in this ${\rm Ext}^1$-group. Then we see that $E$ has an indecomposable submodule representing a class in
\begin{align*}
{\rm Ext}^1(\widetilde{\mathcal{E}}(\X^+_{r,s})^*_{r^\vee,s},\mathcal{E}^\pm(\X^+_{r^\vee,s^\vee})_{r^\vee,s})\ \ \ \ \ \ \ {\rm or}\ \ \ \ \ \ \ \ {\rm Ext}^1(\widetilde{\mathcal{E}}(\X^+_{r^\vee,s^\vee})^*_{r^\vee,s},\mathcal{E}^\pm(\X^+_{r,s})_{r^\vee,s})
\end{align*}
(for the definition of $\widetilde{\mathcal{E}}(\X^+_{r,s})_{r^\vee,s}$ and $\widetilde{\mathcal{E}}(\X^+_{r^\vee,s^\vee})_{r^\vee,s}$, see Definition \ref{kaisei}).
From Propositions \ref{Ext2+} and \ref{Ext2-+}, we see that these ${\rm Ext}^1$-groups are trivial, and thus we have a contradiction.
\end{proof}
\end{lem}
Assume that ${\rm Ext}^1(\widetilde{\mathcal{N}}(h_{r,s}),L(h_{r,s}))=0$. Then, by Lemmas~\ref{null20230322+} and \ref{null20230322-}, we see that $\widetilde{\mathcal{N}}(h_{r,s})$ is the projective cover of $L(h_{r,s})$. But this contradicts Corollary \ref{prop202302110}. Thus we obtain
$
{\rm Ext}^1(\widetilde{\mathcal{N}}(h_{r,s}),L(h_{r,s}))\neq 0.
$
Fix an arbitrary indecomposable module representing a class in the ${\rm Ext}^1$-group ${\rm Ext}^1(\widetilde{\mathcal{N}}(h_{r,s}),L(h_{r,s}))$ and denote it by $\mathcal{P}'(h_{r,s})$.
$\mathcal{P}'(h_{r,s})$ has the following length five socle series
\begin{align*}
\begin{aligned}
&{\rm S}_{1}=L(h_{r,s}),\quad {\rm S}_{2}/{\rm S}_{1}=\X^+_{r,s}\oplus \X^+_{r^\vee,s^\vee},\quad {\rm S}_3/{\rm S}_2=2\X^-_{r^\vee,s}\oplus L(h_{r,s})\oplus 2\X^-_{r,s^\vee},\\
&{\rm S}_4/{\rm S}_3=\X^+_{r,s}\oplus \X^+_{r^\vee,s^\vee},\qquad {\rm S}_5/{\rm S}_4=L(h_{r,s}).
\end{aligned}
\end{align*}
Consider the submodule $(\mathcal{P}'(h_{r,s}):L(h_{r,s}):2)$ of $\mathcal{P}'(h_{r,s})$. 
By Lemm \ref{lem202301240}, we see that this submodule is isomorphic to $\mathcal{Q}(h_{r,s})$. 
\begin{prop}
\label{null20230}
We have
$
{\rm Ext}^1(\mathcal{P}'(h_{r,s}),L(h_{r,s}))=0.
$
\begin{proof}
By Lemm \ref{lem202301240} and by the exact sequence
\begin{align*}
0\rightarrow \mathcal{Q}(h_{r,s})\rightarrow \mathcal{P}'(h_{r,s})\rightarrow \mathcal{P}'(h_{r,s})/\mathcal{Q}(h_{r,s})\rightarrow 0,
\end{align*}
we obtain the exact sequence
\begin{align*}
\begin{aligned}
0&\rightarrow \C\rightarrow \C\rightarrow \C\\
&\rightarrow {\rm Ext}^1(\mathcal{P}'(h_{r,s})/\mathcal{Q}(h_{r,s}),L(h_{r,s}))\rightarrow {\rm Ext}^1(\mathcal{P}'(h_{r,s}),L(h_{r,s}))\rightarrow 0.
\end{aligned}
\end{align*}
By Proposition~\ref{Ext} and Lemma~\ref{lem202301240}, we see that 
\begin{align*}
{\rm Ext}^1(\mathcal{P}'(h_{r,s})/\mathcal{Q}(h_{r,s}),L(h_{r,s}))=\mathbb{C}. 
\end{align*}
Thus from the above exact sequence, we obtain
$
{\rm Ext}^1(\mathcal{P}'(h_{r,s}),L(h_{r,s}))=0.
$
\end{proof}
\end{prop}
\begin{prop}
\label{null2023+}
We have
\begin{align*}
{\rm Ext}^1(\mathcal{P}'(h_{r,s}),\X^+_{r,s})={\rm Ext}^1(\mathcal{P}'(h_{r,s}),\X^+_{r^\vee,s^\vee})=0.
\end{align*}
\begin{proof}
We only prove the first equality. The second equality can be proved in the same way. By Proposition \ref{Ext2+}, we see that
\begin{align}
\label{kore20230}
{\rm Ext}^1(\mathcal{Q}(h_{r,s}),\X^+_{r,s})=0.
\end{align}
Note that $\mathcal{P}'(h_{r,s})/\mathcal{Q}(h_{r,s})\simeq \mathcal{N}_{r,s}/\mathcal{K}_{r,s}$. Thus, by Lemma \ref{lem202301142} and by (\ref{kore20230}), we obtain ${\rm Ext}^1(\mathcal{P}'(h_{r,s}),\X^+_{r,s})=0$.
\end{proof}
\end{prop}

By Proposition \ref{Ext} and Lemma \ref{null20230322-}, we obtain the following proposition.
\begin{prop}
\label{null2023-}
We have
\begin{align*}
{\rm Ext}^1(\mathcal{P}'(h_{r,s}),\X^-_{r^\vee,s})={\rm Ext}^1(\mathcal{P}'(h_{r,s}),\X^-_{r,s^\vee})=0.
\end{align*}
\end{prop}

By Propositions~\ref{null20230},~\ref{null2023+}, and~\ref{null2023-}, we have $\mathcal{P}(h_{r,s})\simeq \mathcal{P}'(h_{r,s})$. Therefore we obtain the following theorem.
\begin{thm}
The projective module $\mathcal{P}(h_{r,s})$ has the following length five socle series:
\begin{align*}
&{\rm S}_{1}=L(h_{r,s}),\quad 
{\rm S}_{2}/{\rm S}_{1}=\X^+_{r,s}\oplus \X^+_{r^\vee,s^\vee},\quad
{\rm S}_3/{\rm S}_2=2\X^-_{r^\vee,s}\oplus L(h_{r,s})\oplus 2\X^-_{r,s^\vee},\\
&{\rm S}_4/{\rm S}_3=\X^+_{r,s}\oplus \X^+_{r^\vee,s^\vee},\quad
\mathcal{P}(h_{r,s})/{\rm S}_4=L(h_{r,s}).
\end{align*}
\end{thm}

\vspace{3mm}
\begin{figure}[htbp]
  \centering
\begin{tikzpicture}[scale=1.05]
\node[inner sep=0.5pt] (a1) at (0, 0) {${\triangle}$};
\node[inner sep=0.5pt] (b1) at (-0.75,-1) {$\heartsuit$};
\node[inner sep=0.5pt] (b2) at (0.75,-1) {$\diamondsuit$};
\node[inner sep=0.5pt] (c1) at (-2,-2) {$\spadesuit$};
\node[inner sep=0.5pt] (c2) at (-1.0,-2) {$\spadesuit$};
\node[inner sep=0.5pt] (c3) at (0,-2) {$\triangle$};
\node[inner sep=0.5pt] (c4) at (1.0,-2) {$\clubsuit$};
\node[inner sep=0.5pt] (c5) at (2,-2) {$\clubsuit$};
\node[inner sep=0.5pt] (d1) at (-0.75,-3) {$\diamondsuit$};
\node[inner sep=0.5pt] (d2) at (0.75,-3) {$\heartsuit$};
\node[inner sep=0.5pt] (e1) at (0, -4) {$\triangle$};

\draw[arrows = {-Stealth[scale=0.7]}] (a1) to (b1);
\draw[arrows = {-Stealth[scale=0.7]}] (a1) to (b2);
\draw[arrows = {-Stealth[scale=0.7]}] (b1) to (c1);
\draw[arrows = {-Stealth[scale=0.7]}] (b1) to (c2);
\draw[arrows = {-Stealth[scale=0.7]}] (b1) to (c3);
\draw[arrows = {-Stealth[scale=0.7]}] (b1) to (c4);
\draw[arrows = {-Stealth[scale=0.7]}] (b1) to (c5);
\draw[arrows = {-Stealth[scale=0.7]}] (b2) to (c1);
\draw[arrows = {-Stealth[scale=0.7]}] (b2) to (c2);
\draw[arrows = {-Stealth[scale=0.7]}] (b2) to (c3);
\draw[arrows = {-Stealth[scale=0.7]}] (b2) to (c4);
\draw[arrows = {-Stealth[scale=0.7]}] (b2) to (c5);
\draw[arrows = {-Stealth[scale=0.7]}] (c1) to (d1);
\draw[arrows = {-Stealth[scale=0.7]}] (c1) to (d2);
\draw[arrows = {-Stealth[scale=0.7]}] (c2) to (d1);
\draw[arrows = {-Stealth[scale=0.7]}] (c2) to (d2);
\draw[arrows = {-Stealth[scale=0.7]}] (c3) to (d1);
\draw[arrows = {-Stealth[scale=0.7]}] (c3) to (d2);
\draw[arrows = {-Stealth[scale=0.7]}] (c4) to (d1);
\draw[arrows = {-Stealth[scale=0.7]}] (c4) to (d2);
\draw[arrows = {-Stealth[scale=0.7]}] (c5) to (d1);
\draw[arrows = {-Stealth[scale=0.7]}] (c5) to (d2);
\draw[arrows = {-Stealth[scale=0.7]}] (d1) to (e1);
\draw[arrows = {-Stealth[scale=0.7]}] (d2) to (e1);
\end{tikzpicture}
\caption{The Loewy diagram of $\mathcal{P}(h_{r,s})$. The triangle $\bigtriangleup$ corresponds to the simple module $L(h_{r,s})$, $\heartsuit$ to $\X^+_{r,s}$, $\diamondsuit$ to $\X^+_{r^\vee,s^\vee}$, $\spadesuit$ to $\X^-_{r,s^\vee}$ and $\clubsuit$ to $\X^-_{r^\vee,s}$.\label{fig202301250}}
\end{figure}

\begin{remark}
Figure \ref{fig202301250} represents the embedding structure of the projective module $\mathcal{P}(h_{r,s})$. This embedding structure is given by \cite{GRW}. 
\end{remark}

\section*{{Acknowledgement}}
We would like to thank Akihiro Tsuchiya, Koji Hasegawa, Gen Kuroki, and Masaru Sugawara for valuable discussions, and Robert McRae and Simon Wood for helpful comments.

\vspace{10mm}
\ \ \ \ H.~Nakano, \textsc{Advanced Mathematical Institute, Osaka Metropolitan University, Osaka 558-8585, Japan}\par\nopagebreak
  \textit{E-mail address} : \texttt{hiromutakati@gmail.com}
  

\begin{thebibliography}{99}
\bibitem[AM1]{AM}
     D. Adamovi\'{c} and A. Milas,
     On the triplet vertex algebra $W(p)$,
     {Adv. Math.} $\bold{217}$ (2008), 2664-2699.
\bibitem[AM2]{AML}D. Adamovi\'{c} and A. Milas,
     Lattice construction of logarithmic modules for certain vertex algebras,
     {Selecta Math}. (N.S.) \textbf{15} (2009), 535-561.
\bibitem[AM3]{AMW2p}
     D. Adamovi\'{c} and A. Milas,
     On $\mathcal{W}$-algebras associated to $(2,p)$ minimal models for certain vertex algebras,
     {Int. Math. Res. Not.} \textbf{10} (2010), 3896-3934.
\bibitem[AM4]{AMW3p}
     D. Adamovi\'{c} and A. Milas,
     On W-algebra extensions of (2, p) minimal models: p $>$3,
     {J. Algebra} \textbf{344} (2011), 313-332.
\bibitem[AM5]{AM Zhu}
     D. Adamovi\'{c} and A. Milas,
     The structure of Zhu’s algebras for certain W-algebras,
     {Adv. Math.} \textbf{227} (2011), 2425-2456
     %; arXiv:1006.5134.
\bibitem[AM6]{AK}
     D. Adamovi\'{c} and A. Milas,
     An explicit realization of logarithmic modules for the vertex operator algebra $W_{p_+,p_-}$,
     {J. Math. Phys.} \textbf{53}, (2012), 16pp.
\bibitem[AK]{AoK}
    K. Aomoto and M. Kita, 
     \textit{Theory of hypergeometric functions}. 
      Springer, Berlin, 2011.               
\bibitem[Ar]{Arike}
       Y. Arike,
       A matrix realization of the quantum group $\mathfrak{g}_{p,q}$, 
       {Int. J. Math.} \textbf{22}(03), (2011), 345-398.
\bibitem[BNW]{BNW}
     B. Boe, D. Nakano and E. Wiesner,
     Category $\mathcal{O}$ for the Virasoro algebra: cohomology and Koszulity,
     {Pacific J. Math}. \textbf{234} (2008), no. 1, 1-21.
\bibitem[CR]{CR}
     T. Creutzig and D. Ridout.
     Logarithmic conformal field theory: Beyond an introduction,
     {J. Phys. A} \textbf{46}(49) (2013), 494006. 
     %arXiv:1303.0847 [hep-th].
\bibitem[Cr]{C}     
M. Cromer, 
Free field realisations of staggered modules in 2D logarithmic CFTs,
arXiv:1612.02909 (2016).
%arXiv preprint arXiv:1612.02909 (2016).   

%追加下
\bibitem[DF1]{DF1}
     V. S. Dotsenko and V. A. Fateev, 
     Conformal algebra and multipoint correlation functions in 2D statistical models,
     {Nucl. Phys.} \textbf{B240} (1984), 312-348. 
\bibitem[DF2]{DF2}
     V. S. Dotsenko and V. A. Fateev, 
     Four-point correlation functions and the operator algebra in 2D conformal invariant theories with central charge $c\leq 1$,
     {Nucl. Phys.} \textbf{B251} (1985), 691-734.
     
\bibitem[FF]{FF}
     B. Feigin and D.B. Fuchs.
     Representations of the Virasoro algebra,
     {in Representations of infinite-dimensional Lie groups and Lie algebras, Gordon andd Breach, New York}(1989).
\bibitem[FGST1]{FF2}
     B. L. Feigin, A.M. Gainutdinov, A.M. Semikhatov, and I. Yu Tipunin,
     Logarithmic extensions of minimal models: characters and modular transformation,
     {Nucl. Phys. B} \textbf{757} (2006), 303-343.
\bibitem[FGST2]{FF3}
     B. L. Feigin, A.M. Ga{i}nutdinov, A.M. Semikhatov, and I. Yu Tipunin, 
     Modular group representations and fusion in logarithmic conformal field theories and in the quantum group center,
     {Comm. Math. Phys.} \textbf{265} (2006), 47–93.
\bibitem[FGST3]{FF4}
     B. L. Feigin, A.M. Ga{i}nutdinov, A.M. Semikhatov, and I. Yu Tipunin, 
     Kazhdan-Lusztig correspondence for the representation category of the triplet W-algebra in logarithmic CFT,
     {Theor. Math. Phys}. 148 (2006) 1210-1235; \textit{Teor. Mat. Fiz}. \textbf{148} (2006), 398-427.
\bibitem[FGST4]{FF22}
     B.L. Feigin, A.M. Gainutdinov, A.M. Semikhatov, and I. Yu Tipunin, 
     Kazhdan-Lusztig-dual quantum group for logarithmic extensions of Virasoro minimal models,
     {J. Math. Phys.} \textbf{48}(3), (2007).     
\bibitem[Fe]{Felder}
     G. Felder,
     BRST approach to minimal models,
     {Nuc. Phy. B.} \textbf{317} (1989), 215-236.
\bibitem[FS]{FS}
     G. Felder and R. Silvotti, 
     Conformal blocks of minimal models on a Riemann surface,
     {Commun. Math. Phys.} \textbf{144} (1992) 17-40.      
\bibitem[FFHST]{FJ}
     J. Fjelstad, J. Fuchs, S. Hwang, A. M. Semikhatov, and I. Yu. Tipunin,
     Logarithmic conformal field theories via logarithmic deformation,
     {Nuclear Phys. B} \textbf{633} (2002), 379-413.
%追加
\bibitem[Fo]{Fo}
        P. J. Forrester, 
        \textit{Log-gases and random matrices (LMS-34)}. 
        Princeton University Press 2010.     

%追加下
\bibitem[FHL]{FHL}
I. Frenkel, Y. Z. Huang, and J. Lepowsky, 
\textit{On axiomatic approaches to vertex operator algebras and modules},
Vol. 494. Amer. Math. Soc., 1993. 

%追加
\bibitem[FHST]{FHST}
     J. Fuchs, S. Hwang, A. M. Semikhatov, and I. Y. Tipunin, 
     Nonsemisimple fusion algebras and the Verlinde formula,
     {Commun. Math. Phys.} \textbf{247} (2004) 713–742.
      %\href{https://doi.org/10.1007/s00220-004-1058-y}{\textit{Commun. Math. Phys.}} \textbf{247} (2004) 713–742, \href{https://arxiv.org/abs/hep-th/0306274}{arXiv:hep-th/0306274}. 

\bibitem[GRW1]{GRW0}
     M. Gaberdiel, I. Runkel, and S. Wood, 
     Fusion rules and boundary conditions in the $c=0$ triplet model,
     {J. Phys. A:math. Theor}. \textbf{42}(32) (2009), 325403.
     % arXiv:0905.0916 [hep-th].
\bibitem[GRW2]{GRW}
     M. Gaberdiel, I. Runkel, and S. Wood, 
     A modular invariant bulk theory for the $c=0$ triplet model,
     {J.Phys. A:math. Theor}. \textbf{44} (2011), 015204.
     % arXiv:1008.0082v1.
\bibitem[GoK]{GK}
     M. Gorelik and V. Kac, 
     On complete reducibility for infinite-dimensional Lie algebras,
     {Adv. Math}. \textbf{226} (2011), no. 2, 1911-1972.
\bibitem[Hu]{H}
     Y. Z. Huang,
     Cofiniteness conditions, projective covers and the logarithmic tensor product theory,
     {J. Pure Appl. Algebra}, \textbf{213}(4), (2009), 458-475.
\bibitem[IK]{IK}
     K. Iohara and Y. Koga,
     \textit{Representation Theory of the Virasoro Algebra},
     Springer, Berlin, 2011.
\bibitem[Kau]{Ka}
     H. G. Kausch,
     Extended conformal algebras generated by multiplet of primary fields,
     {Phys. Lett. B} \textbf{259} (1991), 448-455.
      

\bibitem[KyR]{KR}
      K. Kyt\"{o}l\"{a} and D. Ridout, 
      On staggered indecomposable Virasoro modules,
      {J. Math. Phys.} \textbf{50} (2009), 123503.
      % arXiv:0905.0108[math-ph].



%追加2
\bibitem[MS1]{MS1}
R. McRae and V. Sopin, 
Fusion and (non)-rigidity of Virasoro Kac modules in logarithmic minimal models at $(p, q)$-central charge,
Phys. Scr., Vol. \textbf{99} (2024), no. 3, Paper No. 035233.

%追加
\bibitem[MS2]{MS}
R. McRae and V. Sopin, A tensor category construction of the $ W_ {p, q}$ triplet vertex operator algebra and applications,
{{SIGMA}}, \textbf{22} (2026), 012.
%\href{https://doi.org/10.3842/SIGMA.2026.012}{\textit{SIGMA}}, \textbf{22} (2026): 012.,
%\href{https://arxiv.org/abs/2508.18895}{arXiv:2508.18895[math.QA]}.

\bibitem[MY]{McRae}
      R. McRae and J. Yang,
      Structure of Virasoro tensor categories at central charge $13-6p-6p^{-1}$ for integers $p>1$,
      arXiv:2011.02170 (2020).
\bibitem[Mi]{Milas}
       A. Milas, 
      Fusion rings for degenerate minimal models, 
       {J. Algebra} \textbf{254} (2002), no. 2, 300-335.
\bibitem[NT]{NT}
      K. Nagatomo and A. Tsuchiya, 
      The Triplet Vertex Operator Algebra $W(p)$ and Restricted Quantum Group at Root of Unity, 
      {Adv. Stdu. in Pure Math., Exploring new Structures and Natural Constructions in Mathematical Physics, Amer. Math. Soc}. \textbf{61} (2011), 1–49.
%arXiv:0902.4607. 
%\bibitem[Na1]{Nakano0}   
%        H. Nakano, ``The category of modules of the triplet W-algebras associated to the Virasoro minimal models, 
%        {the Doctor Thesis, Mathematical Institute, Tohoku University} (2023).          
\bibitem[Na1]{Nakano}
      H. Nakano, 
      Explicit formulas of the logarithmic couplings of certain staggered Virasoro modules, 
      {Lett. Math. Phys.} \textbf{113}(2), (2023), 44.          
%追加2
\bibitem[Na2]{Na3}
H. Nakano, Fusion rules for the triplet $W$-algebra $W_{p+,p_-}$, arXiv: 2308.15954 (2023). 
\bibitem[Na3]{Na0}           
H. Nakano, Derivations on the triplet $W$-algebras with $\mathfrak{sl}_2 $-symmetry, arXiv:2604.25729 (2026).

        
\bibitem[Li]{Lin}
      X. Lin, Fusion rules of Virasoro vertex operator algebras, 
      {Proc. Amer. Math. Soc.} \textbf{143}(9), (2015), 3765-3776.
\bibitem[Ra]{Ra}
      J. Rasmussen, 
      W-extended logarithmic minimal models, 
      {Nucl. Phys. B} \textbf{807} (2009), 495. 
      %[0805.2991 [hep-th]].  

%追加
\bibitem[Se]{Se}
A. Selberg,  Bemerkninger om et Mulitiplet Integral, {Norsk Mat. Tidsskrift} \textbf{26} (1944), 71-78.  
       
 %追加下
\bibitem[Su]{S1}
      E. Sussman,
      The regularization of Dotsenko–Fateev integrals,
      {Letters in Mathematical Physics}, \textbf{113} (2023), 1-29.      
                  
\bibitem[TK]{TK}
      A. Tsuchiya and Y. Kanie, 
     Fock space representations of the Virasoro algebra - Intertwining operators,
     {Publ. RIMS, Kyoto Univ}. \textbf{22} (1986), 259-327. 
\bibitem[TW]{TW}
		A. Tsuchiya and S. Wood,
		On the extended W-algebra of type $sl_2$ at positive rational level,
		{Int. Math. Res. Not.}, Volume 2015, Issue 14, (2015), 5357-5435.

\bibitem[Wo]{W}
        S. Wood, Fusion Rules of the $\mathcal{W}_{p,q}$ Triplet Models,
        {J. Phys. A} \textbf{43} (2010), 045212.

%追加
\bibitem[Yam]{Yam}
Y. Yamada,
\textit{Introduction to Conformal Field Theory} (in Japanese),
Suributsuri Series, \textbf{1},
Baifukan, (2006).
%\bibitem[Yam]{Yam}
%Y. Yamada, \textit{Introduction to Conformal field theories(in Japanese)}, suributsuri-series, \textbf{1}, BaifuKan (2006).        
\bibitem[Yan]{Yanagida}
        S. Yanagida. 
        Norm of logarithmic primary of Virasoro algebra, 
        {Lett. Math. Phys.} \textbf{98}(2), (2011), 133-156.
\bibitem[Zh]{Zhu}
     Y. Zhu,
     Modular invariance of characters of of vertex operator algebras,
     \textit{J. Amer. Math. Soc.} \textbf{9} (1996), 237-302.
\end{thebibliography}
\end{document}